\documentclass[a4paper,11pt,twoside,reqno]{amsart}

\usepackage[utf8]{inputenc}
\usepackage{hyperref}
\usepackage{amssymb,amsthm}
\usepackage[margin=1in]{geometry}
\usepackage{slashed}

\newtheorem{Satz}{Theorem}[section]
\newtheorem{Prop}[Satz]{Proposition}
\newtheorem{Lem}[Satz]{Lemma}

\theoremstyle{definition}
\newtheorem{Dfn}[Satz]{Definition}
\newtheorem{Bem}[Satz]{Remark}

\newcommand{\tr}{\operatorname{Tr}}
\newcommand{\Ric}{\operatorname{Ric}}
\newcommand{\Id}{\operatorname{Id}}
\newcommand{\Scal}{\operatorname{Scal}}

\newcommand{\grad}{\operatorname{grad}}

\newcommand{\sff}{\mathrm{I\!I}}

\newcommand{\C}{{\mathbb{C}}}
\newcommand{\dv}{\,\mathrm{d}v}   % I think \, is enough extra spacing, \text{ } is a bit too wide

\newcommand{\R}{\mathbb{R}}
\allowdisplaybreaks[1]

\renewcommand{\epsilon}{\varepsilon}
\numberwithin{equation}{section}

\usepackage{booktabs}
\usepackage{empheq}
\usepackage{tensor}
\usepackage[inline]{enumitem}
\usepackage{mathrsfs}
\setenumerate{font=\upshape,label=(\roman*),itemsep=0.1cm}

\newcommand{\autdiff}{\delta}
\newcommand{\diff}{\mathrm{d}}
\newcommand{\Y}{\mathrm{Y}}
\newcommand{\id}{\mathrm{Id}}
\newcommand{\iprod}{\mathbin{\lrcorner}}
\newcommand{\Fer}{\mathscr{F}} %Fermion bundle
\newcommand{\Hig}{\mathscr{H}} %Higgs field bundle
\newcommand{\dirac}{\slashed{\mathrm{D}}{}}
\newcommand{\g}{\mathfrak{g}}

\newcommand{\Riem}{\operatorname{Riem}}
\newcommand{\comment}[1]{}
\newcommand{\Langle}{\langle\!\langle}
\newcommand{\Rangle}{\rangle\!\rangle}
\newcommand{\lLangle}{\left\langle\!\!\left\langle}
\newcommand{\rRangle}{\right\rangle\!\!\right\rangle}
\newcommand{\proj}{\mathrm{proj}}

\renewcommand{\Re}{\mathrm{Re}}
\renewcommand{\Im}{\mathrm{Im}}

\DeclareMathOperator{\Cl}{Cl}

\DeclareMathOperator{\SO}{\mathbf{SO}}
\DeclareMathOperator{\Spin}{\mathbf{Spin}}
\DeclareMathOperator{\End}{End}
\DeclareMathOperator{\Aut}{Aut}
\DeclareMathOperator{\Ad}{Ad}

\DeclareMathOperator{\Alt}{Alt}
\DeclareMathOperator{\comp}{comp}

\usepackage{color}

\title[Critical points of the Standard Model on expanding spacetimes I: Four dimensions]{Global critical points of the Standard Model on four-dimensional spacetimes of expanding type}
\date{\today}

\author{Volker Branding}
\address{University of Vienna, Faculty of Mathematics\\
Oskar-Morgenstern-Platz 1, 1090 Vienna, Austria\\}
\email{volker.branding@univie.ac.at}

\author{Marko Sobak}
\address{University of Vienna, Faculty of Mathematics\\
Oskar-Morgenstern-Platz 1, 1090 Vienna, Austria\\}
\email{marko.sobak@univie.ac.at}

\thanks{The authors gratefully acknowledge the support of the Austrian Science Fund (FWF) through the projects "Geometric Analysis of Biwave Maps" (DOI: 10.55776/P34853) and 
"The Standard Model as a Geometric Variational Problem" (DOI: 10.55776/P36862).
}

\makeatletter

\def\@testdef #1#2#3{%
  \def\reserved@a{#3}\expandafter \ifx \csname #1@#2\endcsname
 \reserved@a  \else
\typeout{^^Jlabel #2 changed:^^J%
\meaning\reserved@a^^J%
\expandafter\meaning\csname #1@#2\endcsname^^J}%
\@tempswatrue \fi}

\makeatother

\begin{document}

\begin{abstract}
The Standard Model of elementary particle physics is one of the most successful models of contemporary theoretical physics being in full agreement with experiments. However, its mathematical structure deserves further investigations both from a geometric and an analytic point of view.

The aim of this manuscript is to provide a mathematically well-defined and self-contained description of the Standard Model in terms of gauge theory and differential geometry on globally hyperbolic manifolds.
Within this setup we then prove the existence of a global solution for the Euler-Lagrange equations of the Standard Model (with the conformal Higgs potential) under the assumptions that the globally hyperbolic manifold is a four-dimensional spacetime of expanding type and small initial data.
This is achieved by establishing a gauge-invariant energy estimate which is of independent mathematical interest.
\end{abstract}

\maketitle

\tableofcontents

\section{Introduction}
The Standard Model of elementary particle physics is one of the greatest successes in modern theoretical physics which has been awarded with a number of Nobel prizes for physics. 
The mathematical formulation of the Standard Model employs the language of gauge theory and is thus naturally connected to differential geometry. When physicists employ the Standard Model 
to make predictions for elementary particle physics
they most often choose flat four-dimensional Minkowski space to work on. While the calculations performed by physicists are in perfect agreement with experimental data they can not be carried out, at least with the mathematical tools available at present,
in a rigorous mathematical fashion. On the other hand, the Standard Model is well-defined as a classical, that is non-quantized, mathematical field theory such that it can be investigated from the point of view of the geometric calculus of variations which is the path that we are following throughout this manuscript. Moreover, in order to couple the Standard Model to General Relativity, it seems necessary to investigate the Standard Model in the most general geometric setup available which are globally hyperbolic manifolds.

The Standard Model comes with different sectors: Yang-Mills, Dirac and Higgs. Each of them has a sound description in terms of fibre bundles over a semi-Riemannian manifold.
There are of course a significant number (and too numerous to cite) of textbooks that provide the mathematical background of each of the sectors of the Standard Model, but for a modern mathematical introduction we refer to the recent book of Hamilton \cite{MR3837560}.
However, it seems that references studying the analytic properties of the entire coupled Standard Model in the curved setting are quite rare, see \S \ref{sec-past} for a brief historic introduction.
Therefore, in this manuscript we provide a mathematically well-defined and self-contained description of the Standard Model in terms of gauge theory and differential geometry on globally hyperbolic manifolds, focusing in particular on the Euler-Lagrange equations, their properties, and the existence of global solutions.

\subsection{Setting}

Let us briefly discuss the setting.
Details are purposefully omitted in this section for brevity, but more complete descriptions are provided in the future sections.

Let $(M,h)$ be a four-dimensional spin spacetime, and $P \to M$ a principal $G$-bundle for a compact Lie group $G$ equipped with an $\Ad$-invariant inner product on its Lie algebra.%
\footnote{
We would like to point out that in physics the gauge group 
used in the Standard Model is \((\mathbf{SU}(3)\times \mathbf{SU}(2)\times \mathbf{U}(1)) / \mathbb{Z}^6\),
while for the mathematical analysis presented in this manuscript any compact Lie group can be employed.
}
The Euler-Lagrange equations associated with the (conformal) Standard Model Lagrangian are given by 
\begin{empheq}[left=\empheqlbrace]{align*}
    \diff_\omega^\ast  F_\omega + \Re \langle \nabla_{\omega}\Phi \otimes \rho_*\Phi \rangle - \tfrac12 \Im \langle \id\cdot\Psi \otimes \chi_*\Psi \rangle &= 0 , \\[0.2cm]
    \Box_{\omega} \Phi - \tfrac16 \Scal\Phi - \lambda|\Phi|^2 \Phi  - \langle \Psi, i\Y^-  \Psi \rangle &= 0 , \\[0.2cm]
    \dirac_{\omega} \Psi + \Y_\Phi\Psi &= 0,
\end{empheq}
with the unknowns being the connection $\omega \in \Omega^1(P,\g)$, the Higgs field $\Phi \in \Gamma(\Hig)$ (for a vector bundle $\Hig$ associated to $P$), and the Dirac field $\Psi \in \Gamma(\Fer_+)$ (for a bundle $\Fer_+$ of chiral spinors twisted by a vector bundle associated to $P$), see \S \ref{sec-standard-model} for all the details.
A solution $(\omega, \Phi, \Psi)$ of these equations will be referred to as a \emph{conformal Standard Model triplet} (we call it conformal due to the choice of Higgs potential, cf.\ (\ref{eq-conformal-potential})).

We study the conformal Standard Model equations as a Cauchy problem, so it is natural to assume that the background spacetime $(M,h)$ is globally hyperbolic.
Topologically, the spacetime is then foliated $M \cong \R \times \Sigma \cong \pi^\ast \Sigma$, where $\pi : M \to \Sigma$ is the projection.
One can naturally consider spatial tensors, i.e.\ sections of $(\pi^\ast T\Sigma)^{\otimes \ell} \otimes (\pi^\ast T^\ast\Sigma)^{\otimes k}$.
In a similar vein, any principal $G$-bundle $P \to M \cong I \times \Sigma$ can be foliated in the sense that $P \cong \pi^\ast P^\Sigma \cong \R \times P^\Sigma$ are diffeomorphic as principal $G$-bundles with natural projections and actions, where $P^\Sigma = \iota_0^\ast P \to \Sigma$ is the spatial principal bundle defined via the embedding $\iota_0 : \Sigma \xhookrightarrow{} \{ 0 \} \times \Sigma$, see \S \ref{sec-spatial-principal-bundle} or \cite[\S 1]{MR1604914} for details.
Any connection $\omega$ on $P$ then has the form $\omega = \diff\tau \otimes \alpha + \omega^\Sigma$, where $\omega^\Sigma$ is the pullback of a one-parameter family of connections on $P^\Sigma$, and $\alpha$ is a section of the adjoint bundle $\Ad P \to M$.
We also decompose the curvature form into its electric part $E_\omega$ and the magnetic part $B_\omega$, so that $F_\omega = \diff t \wedge E_\omega + B_\omega$.

The metric on the spacetime $(M,h)$ can without loss of generality be put into the form \cite{MR2163568}
\begin{equation*}
    h_{(t,x)} = - N(t,x)^2 \, \diff t \otimes \diff t + (g_{t})_x, \qquad (t,x) \in \R \times \Sigma,
\end{equation*}
where $(g_t)_{t \in \R}$ is a a smooth family of Riemannian metrics on $\Sigma$.
We will further restrict ourselves to spacetimes of \emph{expanding type},
for which there exists a smooth positive function $s : \R \to \R$ depending only on $t$, with $1/s \in L^1[0,\infty)$, and such that the conformal metrics 
\begin{equation*}
    \tilde h = (Ns)^{-2} h
\end{equation*}
have uniformly bounded geometry of a given order $k$, see Definition \ref{def-expanding-spacetime} for the technical details. Here, $s$ can be viewed as the expansion factor of the spacetime.
The conformally transformed spacetime $(M, \tilde h)$ is then Gaussian foliated by Cauchy hypersurfaces of bounded geometry.
The assumption of $\tilde g$ having uniformly bounded geometry enables us to use global intrinsic Sobolev spaces, as it ensures the validity of the usual Sobolev apparatus \cite{MR2343536}. 
In particular, given a Hermitian vector bundle $V \to M$ with a connection, we can define the $H^k$-norm of $V$-valued spatial tensors via
\begin{equation*}
    \Vert \xi \Vert_{H^k}^2 = \sum_{\ell=0}^k \int_\Sigma |\widetilde{D}^\ell \xi|_{\tilde g}^2 \, \dv_{\tilde g},
\end{equation*}
where $\widetilde{D}$ is the connection induced by the Levi-Civita connection of $(\Sigma,\tilde g)$ and the connection on $V$ (one could of course also define the $H^k$-norm also with respect to the original metric $g$; the relation between the two will be studied in \S \ref{sec-global-original}).
We denote by $C^0(I, H^k(V))$ the space of measurable sections $\xi$ for which the $H^k$-norm is uniformly bounded for $t \in I$.
One can also define an intrinsic Sobolev topology on the space of connections on $P^\Sigma$ \cite{MR1175322}, and the resulting space can be viewed as the topological sum
\begin{equation*}
    \mathcal{H}^k(P^\Sigma) = \coprod_{i \in \mathscr{I}} \left( \sigma_i + H^k_{\sigma_i}(T^\ast\Sigma \otimes \Ad P^\Sigma) \right)
\end{equation*}
for some index set $\mathscr{I}$ and smooth connections $\sigma_i$ of bounded geometry of order $k$ on $P^\Sigma$ (with respect to $\tilde g$), and where the subscript $\sigma_i$ in $H^k_{\sigma_i}$ emphasizes that the connection on $\Ad P^\Sigma$ is the one induced by $\sigma_i$, cf.\ also \S \ref{sec-sobolev-connections} for additional details.
At spacetime level, we then consider the space of uniformly Sobolev (spatial) connections on $I\times\Sigma \subset M$ given by
\begin{equation*}
    C^0\left(I, \mathcal{H}^k(P^\Sigma)\right) = \coprod_{i \in \mathscr{I}}
    \left( \omega_i + C^0(I, H^k_{\omega_i}(\pi^\ast T^\ast \Sigma \otimes \Ad P)) \right),
\end{equation*}
where $\omega_i = \Pi^\ast \sigma_i$ and $\Pi : P \to P^\Sigma$ denotes the natural projection.
Note that connections $\omega \in C^0(I, \mathcal{H}^k(P^\Sigma))$ are not the most general connections on $P$ since they have no temporal component (with respect to the splitting $\R \times P^\Sigma \cong P$).
However, this does not restrict the generality of the solutions, since connections with temporal coefficients of appropriate regularity can always be reached by a bundle automorphism, cf.\ Theorem \ref{thm-main} below.

\subsection{Main result}

Roughly speaking, our main result concerns the (future) global existence of conformal Standard Model triplets assuming small initial data.  
It can be formulated as follows.

\begin{Satz}\label{thm-main}
    Let $k > 3/2$ be an integer.
    Let $(M, h)$ be a four-dimensional spin spacetime of expanding type with scale factor $s$, and let $P \cong \pi^\ast P^\Sigma \to M$ be a principal $G$-bundle for a compact Lie group $G$.
    Fix an initial connection
    \begin{equation*}
         \sigma^0+\eta^0 \in \mathcal{H}^{k+1}(P^\Sigma), 
    \end{equation*}
    where $\eta^0 \in H^{k+1}_{\sigma^0}(T^\ast \Sigma \otimes \Ad P^\Sigma)$ and $\sigma^0$ is a smooth connection of bounded geometry of order $k+1$ on $P^\Sigma$, as well as initial data
    \begin{equation*}
        \begin{array}{lll}
            E^0 \in H^{k+1}_{\sigma^0}(T^\ast \Sigma \otimes \Ad P^\Sigma), \quad 
            \Phi^0 \in H^{k+1}_{\sigma^0}(\iota_0^\ast \Hig), \quad
            \dot\Phi^0 \in H^{k}_{\sigma^0}(\iota_0^\ast \Hig), \quad
            \Psi^0 \in H^{k+1}_{\sigma^0}(\iota_0^\ast \Fer_+),  
        \end{array}
    \end{equation*}
    satisfying the constraint
    \begin{equation*}
        - D_{\sigma^0}^\ast E^0 + \langle[\eta^0 \otimes E^0]\rangle + \Re \langle \dot\Phi^0 \otimes \rho_*\Phi^0 \rangle - \frac12\Im \langle \partial_t \cdot\Psi^0 \otimes \chi_*\Psi^0 \rangle = 0.
    \end{equation*}
    Then there is an $\varepsilon > 0$ such that if the initial data are small in the sense that
    \begin{equation*}
        \Vert F_{\sigma^0+\eta^0} \Vert_{H^{k+1}_{\sigma^0}} + \Vert \eta^0 \Vert_{H^{k+1}_{\sigma^0}} + \Vert E^0 \Vert_{H^{k+1}_{\sigma^0}} + \Vert \Phi^0 \Vert_{H^{k+1}_{\sigma^0}} + \Vert \dot\Phi^0 \Vert_{H^{k}_{\sigma^0}} + \Vert \Psi^0 \Vert_{H^{k+1}_{\sigma^0}} < \varepsilon, 
    \end{equation*} 
    then there exists a unique triplet $(\eta,\Phi,\Psi) \in \Gamma(\pi^\ast T^\ast \Sigma \otimes \Ad P) \times \Gamma(\Hig) \times \Gamma(\Fer_+)$ defined on $[0, \infty) \times \Sigma$ (future globally) such that 
    \begin{equation*}
        \omega = \omega^0 + \eta = \Pi^\ast \sigma^0 + \eta
    \end{equation*}
    defines a connection on $P$, and $(\omega, \Phi, \Psi)$ is a conformal Standard Model triplet on $(M, h)$ satisfying the initial conditions
    \begin{equation*}
        \iota_0^\ast \eta = \eta^0, \qquad 
        \iota_0^\ast E_\omega = E^0, \qquad 
        \iota_0^\ast \Phi = \Phi^0, \qquad 
        \iota_0^\ast \frac{\nabla_\omega \Phi}{\diff t} = \dot\Phi^0, \qquad
        \iota_0^\ast \Psi = \Psi^0,
    \end{equation*}
    with regularities
    \begin{align*}
        \omega  &\in C^0\left([0,\infty), \mathcal{H}^{k+1}(P^\Sigma)\right),\\
        E_\omega &\in C^0\left([0,\infty), H^{k+1}_{\omega^0}(\pi^\ast T^\ast \Sigma \otimes \Ad P)\right) \cap C^1\left([0,\infty),  H^{k}_{\omega^0}(\pi^\ast T^\ast \Sigma \otimes \Ad P)\right),\\
        B_\omega &\in C^0\left([0,\infty), H^{k+1}_{\omega^0}(\pi^\ast \Lambda^2 \Sigma \otimes \Ad P)\right) \cap C^1\left([0,\infty), H^{k}_{\omega^0}(\pi^\ast \Lambda^2 \Sigma \otimes \Ad P)\right),\\
        \Phi &\in C^0\left([0,\infty), H^{k+1}_{\omega^0}(\Hig)\right) \cap C^1\left([0,\infty), H^{k}_{\omega^0}(\Hig)\right),\\
        \Psi &\in C^0\left([0,\infty), H^{k+1}_{\omega^0}(\Fer_+)\right) \cap C^1\left([0,\infty), H^{k}_{\omega^0}(\Fer_+)\right),
    \end{align*}
    and $E_\omega$ and $\Phi$ decay uniformly at the rate $s^{-1}$ while $\Psi$ decays uniformly at the rate $s^{-\frac32}$, as $t \to \infty$.
    Furthermore, if
    \begin{equation*}
        \alpha \in C^0\left([0,\infty), H^{k+2}_{\omega^0}(\Ad P)\right),
    \end{equation*} 
    then there exists a unique bundle automorphism $f = f_\alpha : P \to P$ which restricts to the identity map of $P^\Sigma$ at $t=0$, and is such that the transformed triplet $f(\omega, \Phi, \Psi)$ is a conformal Standard Model triplet with $(f\omega)(\partial_t) = \alpha$ and the same regularity and decay properties as $(\omega,\Phi,\Psi)$. 
\end{Satz}

This is achieved via conformal techniques, and in particular the main ingredient of the proof is a gauge-invariant higher order energy estimate.

\begin{Bem}
    Here, $\varepsilon$ depends on the following fixed quantities: the integer $k$, the geometric data of $(M,\tilde h)$, the integral $\left\Vert 1/s \right\Vert_{L^1[0,\infty)}$ and value $s(0)$, and the Standard Model couplings $\lambda, \rho_*, \chi_*, \Y$ (cf.\ \S \ref{sec-standard-model}).
\end{Bem}

\begin{Bem}
    The final part of the theorem shows that we do not lose generality by working with connections with no temporal coefficient modulo bundle automorphisms, provided that we restrict temporal coefficients to the appropriate space.
    One may even take this one step further and show that the solution $(\omega, \Phi, \Psi)$ produced by the theorem is actually unique up to bundle automorphisms adapted to the Sobolev setting in the sense of Eichhorn and Heber \cite{MR1458665}, cf.\ also \S \ref{sec-sobolev-connections}. We do not discuss the full configuration quotient space here so as to not further overload the contents of the paper, but we plan to study this in a future work as this is of independent interest.
\end{Bem}

\subsection{Past and present}
\label{sec-past}

Let us give an overview on previous results that 
provide existence results for certain sectors of the Standard Model and 
are connected to the main theorem established in this manuscript. Here, we present the results in chronological order.
\begin{enumerate}
    \item In \cite{MR649158}, \cite{MR649159} Eardley and Moncrief established an existence result for the Yang-Mills-Higgs equations on four-dimensional Minkowski space.
    \item In the case of four-dimensional Riemannian manifolds Parker 
    studied the mathematical structure of the Standard Model in \cite{MR677998}.
    \item The seminal work of Choquet-Bruhat and  Christodoulou \cite{MR654209} provides an existence result for the Yang-Mills-Higgs-Dirac system in four-dimensional Minkowski space assuming small initial data in certain weighted Sobolev spaces. To obtain their result the authors exploit the conformal structure of the Euler-Lagrange equation
    by mapping four-dimensional Minkowski space onto the Einstein cylinder.
    \item Chru\'sciel and Shatah \cite{MR1604914} established an existence result for the Yang-Mills sector of the Standard Model on four-dimensional globally hyperbolic manifolds.
    \item Psarelli \cite{MR2131047} investigated the massive Dirac-Maxwell system on four-dimensional Minkowski space and proved an existence result assuming small initial data.
    \item More recently, Ginoux and Müller used the method of conformal extension of a globally hyperbolic manifold to prove an existence result for the massless Dirac-Maxwell system \cite{MR3846239} assuming that the initial data is small in a weighted Sobolev space. 
    \item In \cite{MR4226448} Dong, LeFloch and Wyatt studied the global evolution of the \(U(1)\) Higgs boson on four-dimensional Minkowski space.
    \item In another recent article \cite{MR4038549} Taujanskas proved the large data decay of Yang-Mills-Higgs fields on four-dimensional Minkowski and de Sitter spacetimes.
\end{enumerate}

The present manuscript contributes to this list by proving the global existence on four-dimensional expanding spacetimes, stated in Theorem \ref{thm-main}.
Aside from this, we focus on presenting the theory in a geometric manner, working with intrinsic quantities.
We would like to put an emphasis on the latter, since the common approach to Yang-Mills theory on curved spacetimes is to choose a convenient gauge in which the equations have desirable analytic properties, whereas in this manuscript we focus on a completely geometric and intrinsic gauge invariant ansatz.
Thus, the strategy followed in this article may be of interest for both geometers interested in studying the analytic properties of the Standard Model system, and analysts interested in learning about the gauge invariant approaches to the theory.

We would also like to point out that, besides the Standard Model, other models that arise in theoretical physics can also be successfully treated on spacetimes of expanding type, such as the equation for wave maps \cite{MR1776183}, and Dirac-wave
maps with curvature term \cite{MR3830277} which are a mathematical version of the supersymmetric non-linear
sigma model of quantum field theory. 
Moreover, the stability of the Milne universe, considered as a Kaluza-Klein reduced spacetime, was established in \cite{MR3951701}.

\subsection{Organization}

The article is organized as follows.
In \S \ref{sec-conventions}, we briefly review some preliminaries necessary for the understanding of this work and fix our conventions.
In \S \ref{sec-standard-model}, we describe the structure of the Standard Model Lagrangian, present the Euler-Lagrange equations, and also discuss some basic properties of the equations.
In \S \ref{sec-gaussian-foliations}, we discuss in more detail spacetimes of expanding type and objects related to them.
In \S \ref{sec-energy-estimates}, we define the notion of natural geometric total energy for the Standard Model and derive an a priori energy estimate.
Finally, in \S \ref{sec-main-result}, we provide a proof of our main result.

\section{Conventions and preliminaries}
\label{sec-conventions}

As the structure of Standard Model couples many different elements of differential geometry and gauge theory, we would like to recall some basic notions and more importantly fix the conventions that we use throughout the article.

\subsection{Spacetime conventions}
We fix an oriented and time-oriented Lorentzian spin manifold $(M,h)$ of dimension four. 
We assume the spacelike convention for the metric $h$, i.e.\ we work with signature $(-+++)$.
The Einstein summation convention will be used, and Greek indices will always represent the full spacetime indices $\{0,1,2,3\}$ while Latin indices will only run through the spatial indices $\{1,2,3\}$.

Since we will eventually wish to work with spinors, it will prove to be more convenient to work in frames rather than in coordinates.
We denote a positively oriented vielbein, i.e.\ a local section of the oriented semi-orthonormal frame bundle $\SO^+(TM)$, by $e = (e_\mu)$, so that $\eta_{\mu\nu} = h(e_\mu, e_\nu)$. 
We also denote the dual frame by $e^\mu$.

Next, let us fix our conventions for differential forms.
Wedge products of forms are defined inductively with the convention%
\footnote{
Some authors prefer to define wedge products without the binomial factor on the right-hand side.
}
\begin{equation}\label{eq-wedge-alt}
    \omega \wedge \theta = \binom{k+\ell}{k} \, \Alt(\omega \otimes \theta),
\end{equation}
for $\omega \in \Omega^k(M), \theta \in \Omega^\ell(M)$, where the alternator map $\Alt : \otimes^k T^\ast M \to \otimes^k T^\ast M$ is given by 
\begin{equation*}
    \Alt(\xi)(X_1,\ldots,X_k) = \frac{1}{k!}\sum_{\sigma \in S_k} \text{sgn}(\sigma)\,\xi(X_{\sigma(1)},\ldots, X_{\sigma(k)}),
\end{equation*}
where $S_k$ is the symmetric group of $k$ elements and $\xi \in \Gamma(\otimes^k T^\ast M)$.
For a $k$-form $\omega \in \Omega^k(M)$, its fully antisymmetric components are $\omega_{\mu_1\ldots\mu_k} = \omega(e_{\mu_1},\ldots,e_{\mu_k})$, and we have
\begin{equation}\label{eq-forms-tensors}
    \omega = \frac{1}{k!} \, \omega_{\mu_1\ldots\mu_k} \, e^{\mu_1} \wedge \cdots \wedge e^{\mu_k}
    =
    \omega_{\mu_1\ldots\mu_k} \, e^{\mu_1} \otimes \cdots \otimes e^{\mu_k},
\end{equation}
where the right-hand side also describes the embedding $\Lambda^kM \xhookrightarrow{} \otimes^k T^\ast M$.
We define the inner product on forms $\Omega^k(M)$ as usual by
\begin{equation*}
   \langle e^{\mu_1} \wedge \cdots \wedge e^{\mu_k}, \, e^{\nu_1} \wedge \cdots \wedge e^{\nu_k} \rangle = \det (\eta^{\mu_a\nu_b}). 
\end{equation*}
Note that, with this choice of inner product, the embedding $\Lambda^k M \hookrightarrow \otimes^k T^\ast M$ is an isometry provided that the inner product on the tensor bundle $\otimes^k T^\ast M$ is defined with an additional factor $1/k!$, explicitly%
\footnote{
Even if one defines wedge products without any additional factor on the right-hand side of (\ref{eq-wedge-alt}), one would have to rescale the inner product on the tensor bundle $\otimes^k T^\ast M$ with a factor $k!$ to make $\Lambda^k M \hookrightarrow \otimes^k T^\ast M$ an isometry.
}
\begin{equation*}
    \langle e^{\mu_1} \otimes \cdots \otimes e^{\mu_k}, \, e^{\nu_1} \otimes \cdots \otimes e^{\nu_k} \rangle = \frac{1}{k!} \, \eta^{\mu_1\nu_1} \cdots \eta^{\mu_k\nu_k}.
\end{equation*}

The connection forms $\tensor{\Gamma}{_\nu^\lambda}$ associated to the Levi-Civita connection $\nabla$ are defined by
\begin{equation}
\label{eq-christoffel-def}
    \tensor{\Gamma}{_\nu^\lambda}(e_\mu) = \tensor{\Gamma}{_{\mu\nu}^\lambda} = e^\lambda(\nabla_{e_\mu}e_\nu),
    \qquad
    \tensor{\Gamma}{_{\mu\nu\lambda}} = h(\nabla_{e_\mu}e_\nu, e_\lambda).
\end{equation}
The compatibility of $\nabla$ with the metric $h$ is equivalent to the skew-symmetry $\tensor{\Gamma}{_\nu^\lambda} = -\tensor{\Gamma}{^\lambda_\nu}$, while the torsion-freedom implies
\begin{equation*}
    \diff e^\lambda = -\tensor{\Gamma}{_\nu^\lambda}\wedge e^\nu,
    \qquad
    [e_\mu,e_\nu] = (\tensor{\Gamma}{_{\mu\nu}^\lambda} - \tensor{\Gamma}{_{\nu\mu}^\lambda}) \, e_\lambda.
\end{equation*}

For the Riemann tensor, we use the conventions%
\footnote{Many authors prefer to flip the last two indices so that
$R_{\mu\nu\lambda\rho} = h(R^{TM}(e_\mu, e_\nu)e_\rho, e_\lambda)$ but we do not follow this tradition.}
\begin{equation}\label{eq-riemann-convention}
    \Riem(X,Y) = [\nabla_X, \nabla_Y] - \nabla_{[X,Y]},
    \qquad
    R_{\mu\nu\lambda\rho} = h(\Riem(e_\mu, e_\nu)e_\lambda, e_\rho).
\end{equation}
We recall the (skew) symmetries
\begin{equation*}
    R_{\mu\nu\lambda\rho} = -R_{\nu\mu\lambda\rho} = -R_{\mu\nu\rho\lambda} = R_{\lambda\rho\mu\nu}.
\end{equation*}
The Ricci tensor is given by
\begin{equation*}
    \Ric = R_{\mu\nu} \, e^\mu \otimes e^\nu, \qquad R_{\mu\nu} = \tensor{R}{_\lambda_\mu_\nu^\lambda}.
\end{equation*}
For the wave operator corresponding to a given connection $\nabla$ on some vector bundle, we use the sign convention 
\begin{equation*}
    \Box = -\nabla^*\nabla = \tr \nabla^2.
\end{equation*}

Finally, throughout the paper we will employ the so-called $\ast$-notation \cite[\S 2.1]{MR2265040}, to simplify the presentation of certain estimates, in particular when the exact formulas for the terms are less important than the general structure of the terms.
More precisely, given tensor fields $A$ and $B$, we denote by $A \ast B$ any tensor field obtained by raising or lowering or contracting (using the metric $h$) any number of indices of the tensor $A \otimes B$, or switching the order of the indices.
For example, if $\theta$ is a 2-covariant tensor field, we can write
\begin{equation*}
    R_{\mu\rho}\tensor{\theta}{^\rho^\lambda}\tensor{\theta}{_\lambda_\nu} \, e^\mu \otimes e^\nu = {\Ric} \ast \theta \ast \theta. 
\end{equation*}
Covariant derivatives then distribute over stars since the metric is parallel.
Of course, the notation is ambiguous in the sense that equality signs lose their meaning, but this is immaterial for the purpose of doing estimates up to constants.

\subsection{Spin geometry}

We briefly recall some elements of spin geometry. For a more extensive review, we recommend \cite{MR1031992} for the Riemannian setting and \cite{MR701244} for the Lorentzian setting.

Let $\gamma : \R^{1,3} \to \Cl(\R^{1,3}, \eta)$ be the Clifford algebra of Minkowski space.
We consider the Pauli matrices
\begin{equation*}
    \sigma_1 = \begin{bmatrix}
        0 & 1\\
        1 & 0
    \end{bmatrix},
    \quad
    \sigma_2 = \begin{bmatrix}
        0 & -i\\
        i & 0
    \end{bmatrix},
    \quad
    \sigma_3 = \begin{bmatrix}
        1 & 0\\
        0 & -1
    \end{bmatrix},
\end{equation*}
and use the Weyl representation $\nu : \Cl(\R^{1,3}, \eta) \to \mathbf{GL}(4)$, such that
\begin{equation*}
    \gamma_0 = (\nu\circ\gamma)(e_0) = \begin{bmatrix}
        0 & I\\ 
        I & 0
    \end{bmatrix},
    \quad
    \gamma_k = (\nu\circ\gamma)(e_k) = \begin{bmatrix}
        0 & -\sigma_k\\
        \sigma_k & 0
    \end{bmatrix}.
\end{equation*}
Elements of the vector space $\Sigma = \C^4$ of the representation $\nu$ are called \textit{(Dirac) spinors}.
We also define the volume element (also known as the fifth gamma matrix $\gamma_5$)
\begin{equation*}
    \omega = i\gamma_0\gamma_1\gamma_2\gamma_3 = \begin{bmatrix}
        I & 0\\
        0 & -I
    \end{bmatrix},
\end{equation*}
as well as the projections $\pi_\pm = \tfrac12 (I \pm \omega)$.
Then any spinor $\psi \in \Sigma = \C^4$ decomposes as
\begin{equation*}
\psi = (\psi_+, \psi_-), \quad\text{where}\quad \psi_\pm \in \C^2 \cong \Sigma_\pm = \pi_\pm(\C^4).
\end{equation*}
The spaces $\Sigma_\pm$ are referred to as the chiral subspaces of $\Sigma$, and elements thereof are called \textit{chiral (Weyl) spinors}.

There is a natural $\Spin^+(1,3)$-invariant inner product on $\Sigma$ given explicitly by
\begin{equation}\label{eq-spinor-inner-product}
    \langle \psi, \phi \rangle = \psi^\dagger \gamma_0\phi = \psi_-^\dagger\phi_+ + \psi_+^\dagger \phi_-,
\end{equation}
where the dagger denotes the conjugate transpose.%
\footnote{
    In the physics literature, one sometimes defines the \textit{Dirac conjugate} of $\psi$ as $\overline\psi = \psi^\dagger \gamma_0$, in which case one can also write $\langle \psi,\phi \rangle = \overline\psi\phi$, although throughout this article we will use the usual inner product notation.
}
One easily verifies that Clifford multiplication of spinors by vectors is symmetric with respect to this inner product.
Note that the chiral subspaces $\Sigma_\pm$ are both null, and in particular the inner product is not positive-definite.%
\footnote{We would like to emphasize that the symmetry of Clifford multiplication and the indefiniteness of the natural geometric inner product of spinors are the main two features that distinguish the Lorentzian setting from the Riemannian setting, at least from an analytic perspective.}

Let $\Spin^+(TM) \to M$ be the principal spin bundle with corresponding spin structure $\Lambda : \Spin^+(TM) \to \SO^+(TM)$.
We denote the spinor bundle by
$$\Sigma M = \mathbf{Spin}^+(TM) \times_\nu \C^4,$$ 
where $\nu : \Spin^+(1,3) \to \mathbf{GL}(4)$ is the spinor representation, i.e.\ the restriction of the Weyl representation above to $\Spin^+(1,3)$. 
Sections of $\Sigma M$ are called \textit{spinor fields} (or most of the time simply \textit{spinors}). 
If $\epsilon : U \subset M \to \Spin^+(TM)$ is a section and $e = \Lambda \circ \epsilon$ is the associated vielbein, we can write a $\Psi \in \Gamma(\Sigma M)$ as $\Psi = [\epsilon, \psi]$, where $\psi : U \to \C^4$.
Clifford multiplication between tangent vectors and spinors then operates as
\begin{equation*}
    e_\mu \cdot \Psi = [\epsilon, \gamma_\mu\psi].
\end{equation*}
We can also Clifford multiply spinors by covectors via the musical isomorphism, which then extends to $k$-forms via the standard (vector space) isomorphism between the Clifford algebra and the exterior algebra. Explicitly, if we take a form $\omega$  with fully antisymmetric components as in (\ref{eq-forms-tensors}), then 
\begin{equation*}
    \omega \cdot \Psi = \frac{1}{k!} \, \omega^{\mu_1\ldots\mu_k} \, e_{\mu_1} \cdots e_{\mu_k} \cdot \Psi.
\end{equation*}

The inner product (\ref{eq-spinor-inner-product}) is invariant under the action of $\Spin^+(1,3)$ and hence induces a bundle metric on $\Sigma M$, such that the chiral subbundles $\Sigma_\pm M$ are null.
The Levi-Civita connection on $\SO^+(TM)$ induces via the spin structure a connection on $\Spin^+(TM)$, which then induces a covariant derivative $\nabla = \nabla^{\Sigma M}$ on the associated spinor bundle $\Sigma M$ which is compatible with the bundle metric.
Explicitly, in a section
\begin{equation*}
    \nabla^{\Sigma M}\Psi = [\epsilon, \nabla^{\Sigma M}\psi], \quad 
    \nabla^{\Sigma M} \psi = \diff \psi + \sigma \psi,
\end{equation*}
where
\begin{equation*}
    \sigma\psi = \frac14 \Gamma^{\nu\lambda} \otimes (\gamma_\nu \gamma_\lambda \psi) 
    = \frac14 \Gamma_{\mu\nu\lambda} \, e^\mu \otimes (\gamma^\nu\gamma^\lambda \psi),
\end{equation*}
and the connection symbols $\Gamma$ are defined as in (\ref{eq-christoffel-def}).
We recall that the curvature tensor of the spinor connection is related to the Riemann tensor via
\begin{equation}\label{eq-riemann-vs-spinor}
    R^{\Sigma M}(X,Y)\Psi = -\frac14 (\Riem(X,Y)e_\alpha) \cdot e^\alpha \cdot \Psi,
\end{equation}
cf.\ \cite[Theorem 4.15]{MR1031992}.
The \textit{Dirac operator} is the map
\begin{equation*}
    \dirac : \Gamma(\Sigma M) \to \Gamma(\Sigma M), \qquad \dirac\Psi = \eta^{\mu\nu} e_\mu \cdot (\nabla^{\Sigma M} \Psi)(e_\nu).
\end{equation*}
Note that \(i\dirac\) is self-adjoint with respect to the \(L^2\)-norm as we are working on a Lorentzian manifold.

\subsection{Bundle-valued forms}
\label{sec-bundle-valued-forms}

Let $E \to M$ be a vector bundle with a bundle metric and a compatible connection $\nabla^E$.
We will denote by $\Omega^k(M,E)$ the bundle of $E$-valued $k$-forms, i.e. sections of $\Lambda^kM \otimes E$.
This then also induces a \textit{twisted exterior derivative} on $E$-valued forms by
\begin{equation}\label{eq-twisted-diff}
    \diff^E (\theta \otimes \Phi) = \diff \theta \otimes \Phi + (-1)^k \, \theta \wedge \nabla^E \Phi,
\end{equation}
where $\theta \in \Omega^k(M)$ and $\Phi \in \Gamma(E)$.
We can also view elements of $\Omega^k(M, E)$ as twisted covariant tensor fields, i.e.\ sections of $\otimes^k T^\ast M \otimes \Ad P$, cf.\ (\ref{eq-forms-tensors}).
On the latter bundle we have a natural affine connection induced by the Levi-Civita connection and the connection on $E$.
Denoting this connection also by $\nabla^E$, we have the relation
\begin{equation}\label{eq-diff-vs-cov}
    \diff^E \eta = (k+1) \, \text{Alt}(\nabla^E \eta), \qquad \eta \in \Omega^k(M, E),
\end{equation}
where the alternator on the right-hand side only acts on the $\otimes^k T^\ast M$ factor.
Throughout the article, we will always write all indices to the right, i.e.\ if $\omega$ is an $E$-valued one-form, then
\begin{equation*}
    (\nabla^E \omega)_{\mu\nu} = (\nabla^E\omega)(e_\mu,e_\nu)
\end{equation*}
is the covariant derivative of $\omega$ in the direction of $e_\mu$ evaluated at $e_\nu$.
Though this quantity would perhaps more commonly be denoted by $\nabla_\mu \omega_\nu$ 
(cf.\ also Footnote \footref{foot-gauge-covariant-derivative} on p.\ \pageref{foot-gauge-covariant-derivative}), 
such notation can be ambiguous when working with bundle-valued quantities.
This notation also better supports the picture of a connection being a map taking sections of $E$ to sections of $T^\ast M \otimes E$.

\subsection{Principal bundles}
Let us also briefly recall some elements of gauge theory.
For a more extensive review of the topic, we recommend the standard textbooks \cite{MR701244,MR3837560}.
Let $G$ be a compact Lie group, equipped with an Ad-invariant positive-definite inner product on the Lie algebra $\g$, and let $G \to P \xrightarrow\pi M$ be a principal bundle.

We recall that a \emph{connection} on $P$ is a $\g$-valued one-form $\omega \in \Omega^1(P,\g)$ which is of $\Ad$-type in the sense that  $r_g^\ast \omega = \Ad_{g^{-1}}\omega$ for each $g\in G$, and $\omega(X^\ast) = X$ for each $X\in \g$, where $X^\ast \in \Gamma(TP)$ denotes the fundamental vertical vector field of $X$.
The connection can in general not be viewed as a global object on $M$, but the difference of two connections can be viewed as an $\Ad P$-valued one-form on $M$, where $\Ad P = P \times_{\Ad} \g$ is the associated \textit{adjoint bundle}.
The space of connections is thus an affine space over $\Omega^1(M,\Ad P)$.

Let $\rho : G \to \mathbf{GL}(V)$ be a representation of $G$ on the vector space $V$ with associated vector bundle $E = P \times_\rho V \to M$.
A connection $\omega \in \Omega^1(P,\mathfrak g)$ induces a \textit{twisted covariant derivative} $\nabla_\omega$ on $E$ given explicitly by
\begin{equation*}
    \nabla_\omega \Phi = [s, \mathrm \nabla_\omega \phi],
    \quad
    (\nabla_\omega \phi)(e_\mu) = (\diff \phi)(e_\mu) + \rho_*((s^\ast\omega)(e_\mu))\phi,
\end{equation*}
where $\Phi = [s,\phi] \in \Gamma(E)$ for a local section $s : U \subset M \to P$ and  $\phi : U \to V$.
This also induces a twisted exterior derivative $\diff_\omega$ on $\Omega^k(M,E)$ by (\ref{eq-twisted-diff}).

The curvature two-form of the connection $\omega$ is defined as
\begin{equation*}
    F_\omega = \diff \omega + \frac12 [\omega \wedge \omega].
\end{equation*}
We can also view $F_\omega$ as an element of $\Omega^2(M, \Ad P)$ and write
\begin{equation*}
    F_\omega = e^\mu\wedge e^\nu \otimes \tfrac12 F_{\mu\nu} = e^\mu \otimes e^\nu \otimes F_{\mu\nu}, \quad F_{\mu\nu} = -F_{\nu\mu} \in \Gamma(\Ad P).
\end{equation*}
Unless otherwise specified, we will view the curvature form as an element of $\Omega^2(M,\Ad P)$ as above.
The curvature form satisfies the Bianchi identity
\begin{equation*}
    \diff _\omega F_\omega = 0,
\end{equation*}
which can by (\ref{eq-diff-vs-cov}) also be written as
\begin{equation*}
    (\nabla_\omega F_\omega)_{\mu\nu\lambda}
    + (\nabla_\omega F_\omega)_{\nu\lambda\mu}
    + (\nabla_\omega F_\omega)_{\lambda\mu\nu} = 0,
\end{equation*}
where $\nabla_\omega$ denotes the induced connection on $\otimes^k T^\ast M \otimes \Ad P$.%

\footnote{
\label{foot-gauge-covariant-derivative}
In the physics literature one usually defines the (local) \textit{gauge covariant derivative}
$\mathcal{D}_\mu = \nabla_\mu + [A_\mu, \cdot],$
in terms of which
$(\nabla_\omega F_\omega)_{\lambda\mu\nu} = \mathcal{D}_\lambda F_{\mu\nu}$, cf.\ also \S \ref{sec-bundle-valued-forms}.
}   

If $\eta \in \Omega^1(M,\Ad P)$, then
\begin{equation}\label{eq-curvature-sum-intro}
    F_{\omega + \eta} = F_\omega + \diff_\omega \eta + \frac12 [\eta \wedge \eta].
\end{equation}
We recall that the square of the twisted exterior differential on any associated bundle $E = P \times_\rho V \to M$ is not generally zero, it satisfies
\begin{equation*}
    \diff_\omega \diff_\omega \Phi = \rho_*(F_\omega) \Phi = e^\mu \wedge e^\nu \otimes \tfrac12 \rho_*(F_{\mu\nu})\Phi,
\end{equation*}
for $\Phi \in \Gamma(E)$.
In particular, for the curvature form we also have
\begin{equation}\label{eq-curvature-divergence-free}
    \diff_\omega^* \diff_\omega^* F_\omega = -\star^{-1} \diff_\omega \diff_\omega \star F_\omega = -\star^{-1} [F_\omega \wedge \star F_\omega] = 0.
\end{equation}
Here, $\star$ denotes the Hodge star defined on forms by
\begin{equation*}
    \theta \wedge \star \,\eta = \langle \theta, \eta \rangle \dv_h,
\end{equation*}
for $\theta, \eta \in \Omega^k(M)$, where $\dv_h$ denotes the volume form of $(M,h)$.

\subsection{Bundle automorphisms}
\label{sec-bundle-aut}

We consider the group $\Aut P$ of bundle automorphisms of $P$, i.e.\ diffeomorphisms $f: P \to P$ which are $G$-equivariant so that $f(p\cdot h) = f(p)\cdot h$ for each $h\in G$.
We recall that automorphisms can be identified with maps $g : P \to G$ such that $g(p\cdot h)=h^{-1} \cdot g(p)\cdot h$ for each $h\in G$, via $f(p) = p\cdot g(p)$.

Each automorphism induces left actions on the space of connections $\omega \in \Omega^1(P,\g)$ and the  space of sections $\xi = [p, v] \in \Gamma(E)$ of any associated bundle $E = P \times_\rho V$ via
\begin{equation*}
    f \omega = \omega \circ \diff f^{-1} = (f^{-1})^\ast \omega, \qquad f\xi = [f(p), v] = [p, \rho(g(p)) v],
\end{equation*}
and the action naturally extends also to $E$-valued forms and tensors.
Note that twisted exterior derivatives transform under the action of $f$ as $\diff_{f\omega} (f\xi) = f\diff_\omega\xi$.
The map 
\begin{equation*}
    \autdiff_\omega f = f\omega - \omega : TP \to \g
\end{equation*}
is said to be the \emph{variation of $f$ with respect to $\omega$}, cf.\ \cite[\S 3]{MR1458665}.
Note that $\autdiff_\omega f$ is a horizontal form of $\Ad$-type and thus can also be identified with an element of $\Ad P$-valued one-form on $M$.
In terms of the correspondence $f(p) = p \cdot g(p)$, we have
\begin{equation}\label{eq-aut-derivative-formula}
    \autdiff_\omega f = \left(\Ad_g-\Id\right)\omega + (g^{-1})^\ast \mu_G 
    = ((g^{-1})^\ast \mu_G) \circ \proj^H_\omega
    = - \diff R_{g^{-1}} \circ \diff g \circ \proj^H_\omega, 
\end{equation}
where $\mu_G$ denotes the Maurer-Cartan form of $G$, and $\proj^H_\omega$ denotes the horizontal projection induced by $\omega$, i.e.\ onto $H = \ker \omega$.
If $\eta$ is an $\Ad P$-valued form on $M$, then we have
\begin{equation}\label{eq-aut-derivative-sum}
    \autdiff_{\omega + \eta} \, f = \autdiff_\omega f + \left(f-\Id\right)\eta.
\end{equation}
More generally if $\xi$ is a section of any associated bundle $E = P \times_\rho W$, we have $\nabla_{f\omega} (f\xi) = f\nabla_{\omega} \xi$, or
\begin{equation}\label{eq-section-aut-derivative}
    \nabla_\omega (f\xi) = f \nabla_\omega\xi - \rho_*(\autdiff_\omega f) (f\xi)
    = f(\nabla_\omega\xi + \rho_*(\autdiff_\omega f^{-1})\xi),
\end{equation}
where we use the fact that $\autdiff_\omega f^{-1} = -f^{-1} \autdiff_\omega f$.
These formulae also extend trivially to $E$-valued tensors and forms.
Since we can view $\autdiff_\omega f \in \Omega^1(M,\Ad P)$, we can also consider higher order derivatives of $f$ using the twisted exterior derivative $\diff_\omega: \Omega^k(M,\Ad P) \to \Omega^{k+1}(M,\Ad P)$. 
In particular, since $fF_\omega = F_{f\omega} = F_{\omega + \autdiff_\omega f}$, the identity (\ref{eq-curvature-sum-intro}) gives
\begin{equation*}
    \diff_\omega \autdiff_\omega f = (f - \Id) \, F_{\omega} - \frac12 [\autdiff_\omega f \wedge \autdiff_\omega f].
\end{equation*}
In terms of covariant derivatives, this gives the commutation law
\begin{equation}\label{eq-aut-derivative-commutator}
    (\nabla_\omega \, \autdiff_\omega f)(X,Y) = (\nabla_{\omega} \, \autdiff_\omega f)(Y,X) + (f - \Id)\,F_\omega(X,Y) - [\autdiff_\omega f(X), \autdiff_\omega f(Y)],
\end{equation}
for all $X,Y \in \Gamma(TM)$, by the relation (\ref{eq-diff-vs-cov}).

\section{Geometric structure of the Standard Model}
\label{sec-standard-model}

Throughout the section, we fix an oriented and time-oriented four-dimensional spin Lorentzian manifold $(M,h)$.

\subsection{Bosons (Yang-Mills-Higgs sector)}
\textit{Bosons} are the so-called \textit{force-carrier particles}, as they mediate forces between all other particles.
These forces include the electromagnetic force, the weak interaction, and the strong interaction. 
Let us describe how elementary bosons are modelled mathematically.

\subsubsection{Yang-Mills energy density}
Let $G$ be a compact Lie group equipped with an Ad-invariant positive-definite inner product on its Lie algebra $\g$.
We would like to note that for the actual Standard Model, the structure group is \cite{MR3837560}
\begin{equation*}
    G = (\mathbf{SU}(3) \times \mathbf{SU}(2) \times \mathbf{U}(1))/\mathbb{Z}^6,
\end{equation*}
but the results of this paper apply for any compact Lie group as above.

Let $G \to P \xrightarrow\pi M$ be a principal bundle equipped with a connection $\omega \in \Omega^1(P,\g)$. 
The connection $\omega$ encodes information about bosons.
The bosonic kinetic energy density is given by the \textit{Yang-Mills density} 
\begin{equation*}
    \omega \mapsto -|F_\omega|^2,
\end{equation*}
where the norm on $\Omega^2(M,\Ad P)$ is the natural one induced by the spacetime metric $h$ and the $\Ad$-invariant inner product on $\g$.

\subsubsection{Bosonic mass and Higgs fields}
We now wish to supplement this density by a potential (mass) term, which should formally appear as $m^2 |\omega|^2$, where $m\in\R_+$.
However, the latter quantity is not gauge invariant.
Nevertheless, one can provide mass to bosons via the celebrated \textit{Higgs sector}, which we now describe.

Let $W$ be a complex linear space equipped with a Hermitian inner product and let $\rho : G \to \mathbf{U}(W)$ be a complex representation. The associated bundle $\Hig = P \times_\rho W$ is called the \textit{Higgs field bundle},%
\footnote{There exists also the notion of a \textit{Higgs bundle}, which should not be confused with the definition given here.}
and sections thereof are called \textit{Higgs fields}.
The connection $\omega$ induces a covariant derivative on $\Hig$ which shall be denoted by $\nabla_{\omega}$.
Since the representation $\rho$ is unitary, the inner product on $W$ also induces a metric on the Higgs field bundle $\Hig$.
The \textit{Higgs density} is given by 
\begin{equation*}
    (\omega, \Phi) \mapsto -\left(|\nabla_{\omega} \Phi|^2 + U(\Phi)\right),
\end{equation*}
where $U: \Hig \to \R$ is a smooth $G$-invariant function.
For the Standard Model, one uses the so-called \textit{Mexican hat potential}
\begin{equation}\label{eq-mexican-hat-potential}
    U(\Phi) = - \mu |\Phi|^2 + \frac{\lambda}{2} |\Phi|^4
\end{equation}
where $\mu,\lambda \geq 0$ are constants.
An alternative choice of potential is the \textit{conformal potential}
\begin{equation}\label{eq-conformal-potential}
    U(\Phi) = \frac16 \Scal |\Phi|^2 + \frac{\lambda}{2} |\Phi|^4,
\end{equation}
where $\Scal = \Scal_h$ is the scalar curvature of $(M,h)$ and $\lambda$ is a real constant (or even a scalar function on $M$).
Throughout the first part of the article we will keep working with a general potential $U$, but we will later work only with the conformal potential.

\subsubsection{Bosons summarized}
The total Lagrangian density for bosons is the Yang-Mills-Higgs density 
\begin{equation*}
    (\omega, \Phi) \mapsto -\left(|F_\omega|^2 + |\nabla_{\omega} \Phi|^2 + U(\Phi)\right),
\end{equation*}
where $\omega \in \Omega^1(P,\g)$ is a connection form, $\Phi\in \Gamma(\Hig)$ is a Higgs field, and $U: \Hig \to \R$ is a smooth $G$-invariant function.

\subsection{Fermions (Dirac-Yukawa sector)}
\textit{Fermions} are the so-called \textit{matter particles}.
They can further be subdivided into \textit{quarks} and \textit{leptons}.
Mathematically, elementary fermionic particles are modelled as twisted chiral spinors, which we now describe.

\subsubsection{Twisted chiral spinors and the Dirac energy density}
Let $V$ be a complex linear space equipped with a Hermitian inner product.
Let $\chi : G \to \mathbf{U}(V)$ be a complex representation and define the associated bundle $\mathscr{S} = P \times_\chi V$. We consider the twisted spinor bundle 
\begin{equation*}
\Fer = \Sigma M \otimes \mathscr{S}.
\end{equation*}
The connection $\omega$ induces a covariant derivative on $\mathscr{S}$, which together with the spinor covariant derivative induces a \textit{twisted spinor covariant derivative} which will be denoted by $\nabla_\omega$.
The \textit{twisted Dirac operator} is the map
\begin{equation*}
    \dirac_{\omega} : \Gamma(\Fer) \to \Gamma(\Fer), \qquad 
    \dirac_{\omega}\Psi = \eta^{\mu\nu} e_\mu \cdot (\nabla_{\omega} \Psi)(e_\nu). 
\end{equation*}
Here, Clifford multiplication ignores the coefficient bundle $\mathscr{S}$ and acts only on the spinorial factor.

Now, to be able to consider the full chiral structure, we assume that the representation $\chi$ splits, so that $V = V_+ \oplus V_-$ is an orthogonal decomposition and $\chi = \chi_+ \oplus \chi_-$ where $\chi_\pm : G \to \mathbf{U}(V_\pm)$ are representations (in general non-isomorphic).
We also define the associated bundles $\mathscr{S}_\pm = P \times_{\chi_\pm} V_\pm$, as well as the \textit{twisted chiral/mixed spinor bundles} respectively as 
\begin{equation*}
    \Fer_+ = (\Sigma_+M \otimes \mathscr{S}_+) \oplus (\Sigma_-M \otimes \mathscr{S}_-), \qquad \Fer_- = (\Sigma_+M \otimes \mathscr{S}_-) \oplus (\Sigma_-M \otimes \mathscr{S}_+).
\end{equation*}
Sections $\Psi \in \Gamma(\Fer_+)$ are called \textit{twisted chiral spinors} and encode information about fermions.
The twisted Dirac operator restricts to a mapping
\begin{equation*}
    \dirac_{\omega} : \Gamma(\Fer_\pm) \to \Gamma(\Fer_\mp).
\end{equation*}
The fermionic kinetic energy density is given by the \textit{Dirac density} 
\begin{equation*}
    (\omega, \Psi) \mapsto -\Re \langle \Psi, i\dirac_{\omega} \Psi \rangle.
\end{equation*}
Again, we would like to point out that 
\(i\dirac_\omega\) is self-adjoint with respect to the \(L^2\)-norm.

\subsubsection{Fermionic mass via the Yukawa coupling}
Note that $\Fer = \Fer_+ \oplus \Fer_-$ is a null decomposition, since $\Sigma_\pm M$ are null and $\mathscr{S}_\pm$ are orthogonal.
In particular, we cannot supplement the energy density with a mass term of the form $m^2 |\Psi|^2$ as such a term is identically zero for $\Psi \in \Gamma(\Fer_+)$. 
The solution to this problem is again provided by the Higgs field, via the so-called Yukawa coupling, which we now proceed to describe.

A \textit{Yukawa map} is an $\R$-linear map $\Y : W \to \mathfrak{u}(V), \, w \mapsto \Y_w$, such that
\begin{enumerate}
    \item $\Y_w (V_\pm) \subset V_\mp$ for each $w\in W$,
    \item $\Y$ is $G$-equivariant in the sense that 
    \begin{equation}\label{eq-yukawa-equivariant}
        \Y_{\rho(g)w}(\chi(g)v) = \chi(g) \Y_w v
    \end{equation}
    for each $g\in G$, $w\in W$ and $v \in V$. 
\end{enumerate}

A Yukawa map thus has block form (with respect to the decomposition $V = V_+ \oplus V_-$)
\begin{equation*}
    \Y_w = \begin{bmatrix}
        0 & -Z_w^\dagger\\
        Z_w & 0
    \end{bmatrix},
    \qquad
    Z_w : V_+ \to V_-.
\end{equation*}

Every Yukawa map induces an $\R$-tensorial map
\begin{equation*}
    \Y : \Gamma(\Hig) \to \Gamma(\End_{\mathfrak{u}}(\mathscr{S})), \quad\text{via}\quad  \Y_{[s,\phi]} \, [s, v] = [s, \Y_\phi v],
\end{equation*}
where $s : U \subset M \to P$ is a local section, and $\End_{\mathfrak{u}}(\mathscr{S})$ denotes the bundle of skew-Hermitian endomorphisms of $\mathscr{S}$.
Note that $\Y$ is well-defined at bundle level in this way owing to the equivariance condition (\ref{eq-yukawa-equivariant}).

We also extend $\Y_\Phi$ to an endomorphism of $\Fer$ by ignoring the spinorial factor.
Then $\Y_\Phi$ can also be viewed as a mapping
\begin{equation*}
    \Y_\Phi : \Gamma(\Fer_\pm) \to \Gamma(\Fer_\mp).
\end{equation*}

The \textit{Yukawa coupling} (fermionic mass) is given by 
\begin{equation*}
    (\Phi, \Psi) \mapsto -\Re \langle \Psi, i\Y_\Phi \Psi \rangle,
\end{equation*}
for a given fixed Yukawa map $\Y$ which can be viewed as a parameter of the theory.

The equivariance property (\ref{eq-yukawa-equivariant}) infinitesimally reads
\begin{equation}\label{eq-yukawa-equivariant-infinitesimal-linear}
    [\chi_*(X), \Y_w] = \Y_{\rho_*(X)w}, \qquad X\in \g,
\end{equation}
where the commutator on the left-hand side is taken on $\mathfrak{u}(V).$
At the bundle level this reads
\begin{equation}\label{eq-yukawa-equivariant-infinitesimal-bundle}
    \chi_*(\xi) \circ \Y_\Phi - \Y_\Phi \circ \chi_*(\xi) = \Y_{\rho_*(\xi)\Phi}, \qquad \xi \in \Gamma(\Ad P).
\end{equation}
We also observe the product rules%
\footnote{
Here we cannot use the above mentioned bundle-level equivariance (\ref{eq-yukawa-equivariant-infinitesimal-bundle}) because the connection $\omega$ is not an $\Ad P$-valued form, so instead one has to check this locally in a section and use (\ref{eq-yukawa-equivariant-infinitesimal-linear}) directly.
}
\begin{align} \nonumber
    \nabla_\omega(\Y_\Phi\Psi) &= \Y_\Phi (\nabla_\omega\Psi) + \Y_{\nabla_{\omega} \Phi}\Psi,
    \\[0.1cm] \label{eq-dirac-yukawa-commutator}
    \dirac_{\omega} (\Y_\Phi\Psi) &= \Y_\Phi (\dirac_{\omega}\Psi) + \Y_{\nabla_\omega\Phi} \cdot \Psi
    = \Y_\Phi (\dirac_{\omega}\Psi) + \eta^{\mu\nu} e_\mu \cdot \Y_{(\nabla_{\omega}\Phi)_\nu}\Psi,
\end{align}
for $\Phi \in \Gamma(\Hig)$ and $\Psi \in \Gamma(\Fer)$.

\subsubsection{Fermions summarized}
For fixed $\omega$ and $\Phi$, the total Lagrangian density for fermions is the Dirac-Yukawa density 
\begin{equation*}
    \Psi \mapsto -\Re \langle \Psi, i(\dirac_{\omega} + \Y_\Phi) \Psi \rangle,
\end{equation*}
where $\Psi \in \Gamma(\Fer_+)$.

\subsection{Lagrangian and Euler-Lagrange equations}
With these notations, the full Lagrangian of the Standard Model is given by
\begin{equation}\label{eq-full-lagrangian}
    L: (\omega,\Phi,\Psi) \mapsto 
    -\left(\vert F_\omega \vert^2  + \vert \nabla_{\omega} \Phi \vert^2 + U(\Phi) + 
    \Re\langle \Psi, i(\dirac_{\omega} + \Y_\Phi) \Psi \rangle\right),
\end{equation}
where $\omega \in \Omega^1(P,\g)$ is a connection, $\Phi \in \Gamma(\Hig)$, $\Psi \in \Gamma(\Fer_+)$, and $U: \Hig \to \R$ is a smooth $G$-invariant function.

\begin{Bem}
    The relative signs in front of the terms in the Lagrangian (\ref{eq-full-lagrangian}) are important both from a physics and a mathematics perspective. In particular, the relative signs between the Yang-Mills and Higgs terms cannot be fixed without loss of generality. The signs also depend on the choice of signature for the Lorentzian metric. Indeed, under a signature inversion $\eta\mapsto -\eta$ of the Minkowski metric, the Higgs kinetic term also flips sign, whereas the Yang-Mills term remains invariant (since the former is the norm of a one-form, and the latter of a two-form). In contrast, the sign(s) in front of the Dirac-Yukawa term are more flexible, as one can always multiply the Yukawa map and/or the inner product on the spinor bundle by a negative constant without leading to inconsistencies.
\end{Bem}

\begin{Satz}\label{thm-el-eqs}
The Euler-Lagrange equations associated to the Lagrangian {\normalfont (\ref{eq-full-lagrangian})} are given by
\begin{subequations}
\begin{empheq}[left=\empheqlbrace]{align}
    \diff_\omega^\ast  F_\omega + \Re \langle \nabla_{\omega}\Phi \otimes \rho_*\Phi \rangle - \tfrac12 \Im \langle \id\cdot\Psi \otimes \chi_*\Psi \rangle &= 0 , \label{eq-yang-mills}\\[0.2cm]
    \Box_{\omega} \Phi - \tfrac12 \grad U_\Phi - \langle \Psi, i\Y^-  \Psi \rangle &= 0 , \label{eq-higgs}\\[0.2cm]
    \dirac_{\omega} \Psi + \Y_\Phi\Psi &= 0, \label{eq-dirac}
\end{empheq}
\end{subequations}
where $\omega \in \Omega^1(P,\g)$ is a connection, $\Phi \in \Gamma(\Hig)$, $\Psi \in \Gamma(\Fer_+)$, and $U: \Hig \to \R$ is a smooth $G$-invariant function.
Here, the currents can locally explicitly be written as
\begin{align*}
    \Re \langle \nabla_{\omega}\Phi \otimes \rho_*\Phi \rangle
        &= \Re \langle (\nabla_{\omega}\Phi)_\mu, \, \rho_*(\xi^a)\Phi \rangle \, e^\mu \otimes \xi_a,\\[0.1cm]
    \Im \langle \id\cdot\Psi \otimes \chi_*\Psi \rangle
        &= \Im \langle e_\mu \cdot \Psi, \, \chi_*(\xi^a)\Psi \rangle\, e^\mu \otimes \xi_a,\\[0.1cm]
    \langle \Psi, i\Y^-  \Psi \rangle
        &= \tfrac{1}{2} \langle \Psi, (i\Y_{W_k} - \Y_{iW_k}) \Psi \rangle \, W_k,
\end{align*}
where $e_\mu$ is a semi-orthonormal frame for $\SO^+(TM)$, $\xi_a$ is an orthonormal frame for $\Ad P$ and $W_k$ is an orthonormal frame for $\Hig$.
\end{Satz}

\begin{proof}
    See Appendix \ref{appendix-el-eqs-derivation}.
\end{proof}

\begin{Bem}\label{rem-gradient-antilinear}
Note that $\grad U_\Phi$ is the bundle gradient of $U$ evaluated at $\Phi$, satisfying 
\begin{equation*}
    \Re \langle \alpha, \grad U_\Phi \rangle = \diff U_\Phi(\alpha),
\end{equation*}
for all $\alpha \in \Gamma(\Hig)$,
and $\Y^-$ is the complex antilinear part of $\Y$, satisfying
\begin{equation}\label{eq-yukawa-complex-antilinear}
    2\Re \langle \Phi, \langle \Psi, i\Y^-  \Psi \rangle \rangle = \langle \Psi, i\Y_\Phi  \Psi \rangle
\end{equation}
for all $\Phi \in \Gamma(\Hig)$ and $\Psi \in \Gamma(\Fer_+)$.
\end{Bem}

\subsection{Euler-Lagrange equations as waves}
We now wish to demonstrate that the Euler-Lagrange equations (\ref{eq-yang-mills}--\ref{eq-dirac}) imply that the variables $(F_\omega,\Phi,\Psi)$ satisfy a suitable system of wave equations.

\begin{Satz}\label{thm-wave-eqs}
    If $(\omega, \Phi, \Psi)$ satisfy the system {\normalfont (\ref{eq-yang-mills}--\ref{eq-dirac})}, then they also satisfy the wave system
    \begin{subequations}
    \begin{empheq}[left=\empheqlbrace]{align}
        &\Box_\omega F_\omega 
        = \langle [F_\omega \wedge F_\omega] \rangle + \mathcal RF_\omega + \Im \langle \id \cdot \Psi \wedge \chi_* \nabla_\omega \Psi \rangle  \label{eq-yang-mills-wave}\\[0.1cm]
        &\hspace{1.4cm} + \Re \langle \rho_*(F_\omega) \Phi \otimes \rho_*\Phi \rangle - \Re \langle \nabla_\omega \Phi \wedge  \rho_*\nabla_\omega\Phi \rangle,\nonumber\\[0.2cm]
        &\Box_\omega \Phi = \tfrac12 \grad U_\Phi + \langle \Psi, i\Y^-  \Psi \rangle , \label{eq-higgs-wave}\\[0.2cm]
        &\Box_\omega \Psi =\tfrac14 \Scal \Psi + \chi_*(F_\omega) \cdot \Psi + \Y_{\nabla_\omega\Phi} \cdot \Psi - \Y_\Phi^2 \Psi, \label{eq-dirac-wave}
    \end{empheq}
    \end{subequations}
    where $\Box_\omega = -\nabla_\omega^\ast \nabla_\omega$ with $\nabla_\omega$ denoting the twisted connection of the corresponding bundle (respectively in the equations, these are $\Lambda^2M \otimes \Ad P$, $\,\Hig$, and $\Fer_+$), and we use the notations
    \begin{align*}
        \langle [F_\omega \wedge F_\omega] \rangle &= e^\mu \wedge e^\nu \otimes  [\tensor{F}{^\alpha_\mu}, F_{\alpha\nu}], \\[0.1cm]
        \mathcal R F_\omega &= e^\mu \wedge e^\nu \otimes \tfrac12 (R_{\mu\alpha} \tensor{F}{^\alpha_\nu} + R_{\nu\alpha} \tensor{F} {_\mu^\alpha}-2R_{\alpha\mu\nu\beta}F^{\alpha\beta}), \\[0.1cm]
        \Re \langle \nabla_\omega \Phi \wedge  \rho_*\nabla_\omega\Phi \rangle
        &=
        e^\mu \wedge e^\nu \otimes \langle (\nabla_{\omega} \Phi)_\mu,\, \rho_*(\xi^a)(\nabla_{\omega} \Phi)_\nu  \rangle\} \, \xi_a, \\[0.1cm]
        \Im \langle \id \cdot \Psi \wedge \chi_* \nabla_{\omega} \Psi \rangle
        &= e^\mu \wedge e^\nu \otimes \tfrac12\Im \{\langle e_\mu \cdot \Psi, \, \chi_*(\xi^a)(\nabla_{\omega}\Psi)_\nu \rangle - \langle e_\nu \cdot \Psi, \, \chi_*(\xi^a)(\nabla_{\omega}\Psi)_\mu \rangle \} \xi_a.
    \end{align*}
\end{Satz}

\begin{proof}
The Higgs equation (\ref{eq-higgs}) is already a wave equation so we only analyze the Yang-Mills equation and the Dirac equation.

If the Dirac equation (\ref{eq-dirac}) is satisfied, we can apply $\dirac_{\omega}$ once more and use the commutation law (\ref{eq-dirac-yukawa-commutator}) to find
\begin{equation*}
    \dirac_{\omega}^2 \Psi = \Y_\Phi^2\Psi - \Y_{\nabla_{\omega}\Phi} \cdot \Psi,
\end{equation*}
and by the Lichnerowicz-Weitzenböck formula (cf.\ Proposition \ref{prop-dirac-lichnerowicz-weitzenbock}) this is equivalent to
\begin{equation*}
    \Box_{\omega} \Psi = \frac14 \Scal\Psi + \slashed{\mathcal R}{}^S_\omega \Psi + \Y_{\nabla_{\omega}\Phi} \cdot \Psi - \Y_\Phi^2 \Psi,
\end{equation*}
where $\slashed{\mathcal R}{}^S_\omega \xi = \chi_*(F_\omega) \cdot \xi$ by Lemma \ref{lem-weitzenbock-curvature-explicit} (ii).

To turn the Yang-Mills equation (\ref{eq-yang-mills}) into a wave equation we employ the Weitzenböck formula for the twisted exterior derivative (cf.\ Proposition \ref{prop-twisted-weitzenbock})
\begin{equation*}
    \diff_\omega^\ast  \diff_\omega^{} + \diff_\omega^{} \diff_\omega^\ast 
    = -\Box_\omega^{\Lambda^k M\otimes \Ad P} + \mathcal{R}_\omega^{\Lambda^k M\otimes \Ad P},
\end{equation*}
where the twisted Weitzenböck curvature operator is given by
\begin{equation*}
    \mathcal{R}_\omega^{\Lambda^k M\otimes \Ad P}\xi = \eta^{\nu\alpha} \, e^\mu \wedge e_\nu \iprod (R_\omega^{\Lambda^k M \otimes \Ad P}(e_\alpha, e_\mu)\xi).
\end{equation*}
More precisely, we can apply $\diff_\omega$ to (\ref{eq-yang-mills}) and $\diff_\omega^\ast $ to the Bianchi identity to get
\begin{align*}
    \Box_\omega F_\omega 
    &= \mathcal R_\omega^{\Lambda^2 M \otimes \Ad P} F_\omega - \diff_\omega^{}\diff_\omega^*F_\omega - \diff_\omega^* \diff_\omega^{} F_\omega \\
    &= \mathcal R_\omega^{\Lambda^2 M \otimes \Ad P} F_\omega + \diff_\omega\Re \langle \nabla_{\omega} \Phi \otimes \rho_*\Phi \rangle - \tfrac12 \diff_\omega \Im \langle \id \cdot \Psi \otimes \chi_* \Psi \rangle,
\end{align*}
where $\Box_\omega = \Box_\omega^{\Lambda^2M \otimes \Ad P}$.
By Lemma \ref{lem-weitzenbock-curvature-explicit}, we explicitly have
\begin{align*}
    \mathcal R^{\Lambda^2 M} F_\omega &=  e^\mu \wedge e^\nu \otimes \frac12 (-2R_{\alpha\mu\nu\beta}F^{\alpha\beta} + R_{\mu\alpha} \tensor{F}{^\alpha_\nu} + R_{\nu\alpha} \tensor{F} {_\mu^\alpha}),\\[0.1cm]
    \mathcal R_\omega^{\Ad P} F_\omega &= e^\mu \wedge e^\nu \otimes [\tensor{F}{^\alpha_\mu}, F_{\alpha\nu}],
\end{align*}
and thus $\mathcal R_\omega^{\Lambda^2 M \otimes \Ad P}F_\omega = \mathcal R F_\omega +  \langle [F_\omega \wedge F_\omega] \rangle$.
Furthermore, using Lemma \ref{lemma-current-differential} we can calculate
\begin{align*}
    \diff_\omega\Re \langle \nabla_{\omega} \Phi \otimes \rho_*\Phi \rangle
    &= \Re \langle \rho_*(F_\omega) \Phi \otimes \rho_*\Phi \rangle
    - \Re \langle \nabla_{\omega} \Phi \wedge  \rho_*\nabla_\omega\Phi \rangle\\[0.1cm]
    &=  e^\mu \wedge e^\nu \otimes \frac12\Re\{\langle \rho_*(F_{\mu\nu})\Phi,\, \rho_*(\xi^a)\Phi \rangle
    - 2\langle (\nabla_{\omega} \Phi)_\mu,\, \rho_*(\xi^a)(\nabla_{\omega} \Phi)_\nu  \rangle\} \, \xi_a,
\end{align*}
and (also using the symmetry of Clifford multiplication)
\begin{align*}
    \frac12 \diff_\omega \Im \langle \id \cdot \Psi \otimes \chi_* \Psi \rangle
    &=  -\Im \langle \id \cdot \Psi \wedge \chi_* \nabla_{\omega} \Psi \rangle\\[0.1cm]
    &= -e^\mu \wedge e^\nu \otimes \frac12\Im \{\langle e_\mu \cdot \Psi, \, \chi_*(\xi^a)(\nabla_{\omega}\Psi)_\nu \rangle - \langle e_\nu \cdot \Psi, \, \chi_*(\xi^a)(\nabla_{\omega}\Psi)_\mu \rangle \} \xi_a.
\end{align*}
Hence, the wave equation for $F_\omega$ reads
\begin{align*}
    \Box_\omega F_\omega 
    &= \mathcal R_\omega^{\Lambda^2 M \otimes \Ad P} F_\omega + \Re \langle \rho_*(F_\omega) \Phi \otimes \rho_*\Phi \rangle
    - \Re \langle \nabla_{\omega} \Phi \wedge  \rho_*\nabla_{\omega}\Phi \rangle + \Im \langle \id \cdot \Psi \wedge \chi_* \nabla_{\omega} \Psi \rangle,
\end{align*}
as desired.
\end{proof}

\subsection{Divergence freedom of currents}

Note that the curvature form of any connection satisfies the identity $\diff_\omega^*\diff_\omega^*F_\omega = 0$, cf.\ (\ref{eq-curvature-divergence-free}) so the currents of the Yang-Mills equation (\ref{eq-yang-mills}) should also have vanishing twisted divergence if the system is to be coherent.

\begin{Prop}\label{prop-currents-div-free}
    Consider the currents of the Yang-Mills equation
    \begin{equation*}
        \mathfrak{J}[\omega,\Phi,\Psi] = - \Re \langle \nabla_{\omega}\Phi \otimes \rho_*\Phi \rangle + \tfrac12 \Im \langle \id\cdot\Psi \otimes \chi_*\Psi \rangle.
    \end{equation*}
    If the Higgs equation {\normalfont (\ref{eq-higgs})} and the Dirac equation {\normalfont (\ref{eq-dirac})} are satisfied, then 
    \begin{equation*}
        \diff_\omega^\ast\, \mathfrak{J}[\omega,\Phi,\Psi] = 0.
    \end{equation*}
\end{Prop}

\begin{proof}
By Lemma \ref{lemma-current-differential} and the fact that the representations $\rho$ and $\lambda$ are unitary, we get
\begin{align}
    \label{eq-higgs-current-divergence}
    \diff_\omega^\ast  \Re \langle \nabla_{\omega} \Phi \otimes \rho_* \Phi \rangle &= -\Re \langle \Box_{\omega} \Phi,\, \rho_* \Phi \rangle,\\[0.1cm]
    \label{eq-dirac-current-divergence}
    \diff_\omega^\ast  \Im \langle \id \cdot \Psi \otimes \chi_*\Psi \rangle &= -2\Im \langle \dirac_{\omega} \Psi, \chi_* \Psi \rangle. 
\end{align}
Now observe that
\begin{equation*}
    \Re\langle \grad U_\Phi, \rho_*(\xi)\Phi \rangle = \diff U_\Phi(\rho_*(\xi)\Phi) 
    = \frac{\diff}{\diff t}\bigg|_{t=0} [U(\rho(\exp(t \xi))\Phi) - U(\Phi)] = 0
\end{equation*}
for each $\xi \in \Gamma(\Ad P)$, since $U$ is $G$-invariant.
On the other hand, using (\ref{eq-yukawa-equivariant-infinitesimal-bundle}, \ref{eq-yukawa-complex-antilinear}), we see that 
\begin{equation*}
    \Re\langle \langle \Psi, i\Y^- \Psi \rangle, \rho_*\Phi \rangle
    = \frac12 \langle \Psi, i\Y_{\rho_*\Phi} \Psi \rangle
    = -\Im \langle \Y_\Phi \Psi, \chi_*\Psi \rangle,
\end{equation*}
since both $\Y_\Phi$ and $\chi_*$ are unitary.
Thus the Higgs and Dirac equations (\ref{eq-higgs}, \ref{eq-dirac}) together with (\ref{eq-higgs-current-divergence}, \ref{eq-dirac-current-divergence}) imply that
\begin{equation*}
    \diff_\omega^\ast\, \mathfrak{J}[\omega,\Phi,\Psi] 
    =
    \diff_\omega^\ast  \Re \langle \nabla_{\omega} \Phi \otimes \rho_* \Phi \rangle - \tfrac12 \, \diff_\omega^\ast  \Im \langle \id \cdot \Psi \otimes \chi_*\Psi \rangle
    =
    \Im \langle \Y_\Phi \Psi + \dirac_{\omega}\Psi, \chi_*\Psi \rangle = 0,
\end{equation*}
as desired.
\end{proof}

\subsection{Choice of Higgs potential and conformality} 

For the Higgs potential, one usually chooses $U(\Phi) = u(|\Phi|^2) = u(x,|\Phi|_x^2)$ where $u:M\times\R\to\R$.
The Higgs equation is then
\begin{align*}
    0 &= \Box_{\omega} \Phi - u'(|\Phi|^2)\Phi - \langle \Psi, i\Y^-  \Psi \rangle,
\end{align*}
where prime denotes differentiation with respect to the second factor and we suppress $x$.
In particular, the choice $u(x,a) = - \mu a + \tfrac{\lambda}{2}a^2$ yields the Mexican hat potential (\ref{eq-mexican-hat-potential}), and the equation is
\begin{equation*}
    \Box_{\omega} \Phi + \mu \Phi - \lambda|\Phi|^2 \Phi - \langle \Psi, i\Y^-  \Psi \rangle = 0.
\end{equation*}
On the other hand $u(x, a) = \tfrac16 \Scal_h(x) a + \tfrac{\lambda}{2}a^2$ yields the conformal Higgs potential and the Higgs equation is 
\begin{equation*}
    \Box_{\omega} \Phi - \frac16\Scal_h \Phi - \lambda|\Phi|^2 \Phi - \langle \Psi, i\Y^-  \Psi \rangle = 0.
\end{equation*}

\begin{Dfn}
    A triplet $(\omega,\Phi,\Psi)$, where $\omega \in \Omega^1(P,\g)$ is a connection, $\Phi \in \Gamma(\Hig)$, and $\Psi \in \Gamma(\Fer_+)$, is called:
    \begin{enumerate}
        \item a \textit{Standard Model triplet} if it satisfies the system (\ref{eq-yang-mills}--\ref{eq-dirac}) with respect to the Mexican hat potential (\ref{eq-mexican-hat-potential}) for some $\mu, \lambda \geq 0$;
        \item a \textit{conformal Standard Model triplet} if it satisfies the system (\ref{eq-yang-mills}--\ref{eq-dirac}) with respect to the conformal Higgs potential (\ref{eq-conformal-potential}) for some $\lambda \in \R$;
    \end{enumerate}
\end{Dfn}

The reason for this nomenclature is the following result.

\begin{table}
\centering
\renewcommand{\arraystretch}{1.5}
\begin{tabular}{rc}
    \toprule
    Metric & $\tilde h = \exp(2f) h = \Omega^2 h$ \\ \midrule
    Spin structure & $\widetilde\Lambda = \Omega^{-1} \Lambda$ \\ \midrule
    Clifford multiplication by vectors & $X \;{\tilde\cdot}\; \psi = \Omega\, X\cdot\psi$ \\ \midrule
    Clifford multiplication by $k$-forms & $\theta \;{\tilde\cdot}\; \psi = \Omega^{-k} \,\theta\cdot\psi$ \\ \midrule
    Volume form & $\dv_{\tilde h} = \Omega^n \dv_h$ \\ \midrule
    Hodge star on $k$-forms & ${\widetilde\star} = \Omega^{n-2k} \, \star$ \\ \midrule 
    Levi-Civita connection & $\widetilde\nabla X = \nabla X + \diff f \otimes X + \diff f(X) \, \id - X^\flat \otimes \grad f$\\ \midrule
    Induced connection on one-forms & $\widetilde\nabla \theta = \nabla \theta - \theta \otimes \diff f - \diff f \otimes \theta + \langle \diff f, \theta \rangle h$\\ \midrule
    Formal adjoint of Levi-Civita & $\widetilde{\nabla}^{*}\theta = \Omega^{-2}\left(\nabla^*\theta - (n-2)\langle\diff f, \theta\rangle\right)$ \\ \midrule
    Wave operator on scalar functions & $\widetilde\Box \phi = \Omega^{-2} \left( \Box \phi + (n-2) \langle \diff f, \diff \phi \rangle \right)$ \\ \midrule
    Scalar curvature & $\widetilde\Scal = \Omega^{-2} \left[\Scal -  2(n-1)\Box f - (n-1)(n-2)|\diff f|^2\right]$ \\ \midrule
    Spin connection & $\widetilde\nabla^{\Sigma M} \psi = \nabla^{\Sigma M} \psi + \tfrac14 (\diff f \cdot \id - \id \cdot \diff f) \cdot \psi$ \\ \midrule
    Dirac operator & $\widetilde{\dirac}\psi = \Omega^{-1}  \left( \dirac\psi + \tfrac12(n-1) \, \diff f \cdot\psi \right)$ \\
    \bottomrule
\end{tabular}
\vspace{0.4cm}
\caption{Various metric operators under a conformal transformation on an $n$-dimensional manifold. The quantities on the right hand side are always taken with respect to $h$.}
\label{tab-conformal-transform}
\end{table}

\begin{Prop}
    \label{prop-conformal-trans}
    Consider a conformal transformation
    \begin{equation*}
        \tilde h = \Omega^2 h, \qquad \widetilde\Phi = \Omega^{-1}\Phi, \qquad \widetilde\Psi = \Omega^{-\frac32} \Psi.
    \end{equation*}
    Then $(\omega,\Phi,\Psi)$ satisfies the equations {\normalfont (\ref{eq-yang-mills}--\ref{eq-dirac})} with respect to $h$ if and only if $(\omega, \widetilde\Phi, \widetilde\Psi)$ satisfies {\normalfont (\ref{eq-yang-mills},\ref{eq-dirac})} with respect to $\tilde h$, while {\normalfont (\ref{eq-higgs})} should be replaced by
    \begin{equation*}
        \widetilde{\Box}_\omega \widetilde\Phi - \langle \widetilde\Psi, i\Y^- \widetilde\Psi \rangle - \tfrac16 \,\widetilde{\Scal}\,\widetilde\Phi - \Omega^{-2} \left( \tfrac{1}{2}\,\Omega^{-1} \grad U_{\Omega\widetilde\Phi} - \tfrac16 \Scal \widetilde\Phi \right) = 0.
    \end{equation*}
    In particular, the equations are conformally invariant if $U$ is the conformal Higgs potential {\normalfont (\ref{eq-conformal-potential})}.
\end{Prop}

\begin{proof}
This follows from Table \ref{tab-conformal-transform}, cf.\ also the discussion in \cite[\S 5]{MR654209}.
\end{proof}

\section{Spacetimes of expanding type and Gaussian foliations}
\label{sec-gaussian-foliations}

The goal of this section is to discuss some constructions and objects on spacetimes of expanding type.
We start with the explicit definition.

\begin{Dfn}
    \label{def-expanding-spacetime}
    Let $(M,h)$ be a four-dimensional globally hyperbolic Lorentzian manifold,
    \begin{equation*}
        M = \R \times \Sigma, \qquad h_{(t,x)} = - N(t,x)^2 \, \diff t \otimes \diff t + (g_{t})_x,
    \end{equation*}
    such that $N$ is a uniformly lower and upper bounded smooth function with
    \begin{equation*}
        \sup_{t\geq0}\Vert N(t,\cdot) \Vert_{C^k} < \infty, \qquad  \sup_{t\geq0}\Vert \partial_t N(t,\cdot) \Vert_{C^{k-1}} < \infty,
    \end{equation*}
    for all integers $k \geq 1$, and $(g_t)_{t \in \R}$ is a a smooth family of Riemannian metrics on $\Sigma$.%
    \footnote{We assume as usual that $\Sigma$ has a sufficiently well-behaved topology, in the sense that it is paracompact and Hausdorff.}
    Then we say that $(M,h)$ is of \emph{expanding type} (with scale factor $s$ and order $k$) if there exists a smooth positive function $s : \R \to \R$ with $1/s \in (L^1 \cap L^\infty)([0,\infty))$ such that the conformal metrics
    \begin{equation*}
        \tilde g_t := \frac{1}{N(t,\cdot)^2\,s(t)^2} \, g_t
    \end{equation*}
    on $\Sigma$ have:
    \begin{itemize}
        \item[(i)] uniformly bounded intrinsic geometry, in the sense that
        \begin{equation*}
            \inf_{t\geq 0}\text{inj}_{(\Sigma,\tilde g_t)} > 0, \qquad \sup_{t\geq 0} \Vert {\Riem_{\tilde g_t}} \Vert_{C^k(\tilde g_t)} < \infty, 
        \end{equation*}
        where the former denotes the injectivity radius, and we also have
        \begin{equation*}
            \sup_{t\geq 0} \Vert \tilde g_t - \tilde g_0 \Vert_{C^{k+2}(\tilde g_0)} < \infty, \qquad \sup_{t\geq 0} \Vert \tilde g_t - \tilde g_0 \Vert_{C^{k+2}(\tilde g_t)} < \infty,
        \end{equation*}
        \item[(ii)] uniformly bounded extrinsic geometry, in the sense that
         \begin{equation*}
            \sup_{t\geq 0}\Vert \widetilde\sff \Vert_{C^{k}(\tilde g_t)} < \infty, \qquad 
            \sup_{t\geq 0} \,s(t)\left\Vert \frac{\nabla\widetilde\sff}{\diff t} \right\Vert_{C^{k-1}(\tilde g_t)} < \infty,
        \end{equation*}
        where $\widetilde\sff = \widetilde\sff_t$ denotes the second fundamental form of $\{t\} \times \Sigma$ in $(M, \tilde h)$.
    \end{itemize}   
\end{Dfn}

The assumptions in condition (i) of the definition are natural when studying spaces of metrics on noncompact manifolds \cite{MR1331100}. 
Note that the $C^k(\tilde g_t)$-norms of tensors on $\Sigma$ are then uniformly equivalent for all $t\geq0$. 

A prototypical example of a spacetime of expanding type is the de Sitter spacetime
\begin{equation*}
    \R \times \mathbb{S}^3, \qquad h = -\diff t\otimes \diff t + a^2\cosh^2 (t/a) \, g_{\mathbb{S}^3}, \quad a > 0
\end{equation*}
which clearly satisfies the conditions with scale factor $s(t)=a\cosh(t/a)$.
More generally, one can consider any Robertson-Walker spacetime 
\begin{equation*}
    \R \times \Sigma, \qquad h = -\diff t\otimes \diff t + s(t)^2 \, g_{\Sigma}, \qquad (\Sigma, g_\Sigma) \text{ three-dimensional space form},
\end{equation*}
whose scale factor $s$ is such that $1/s$ is integrable.
Even more generally, one can replace the space form $(\Sigma, g_\Sigma)$ by any three-dimensional Riemannian manifold of bounded geometry. 
(In all of these examples, the second fundamental form on the conformally transformed side vanishes identically $\widetilde\sff \equiv 0$, while the Riemann tensor is constant with respect to $t$ and has bounded $C^k$-norms for all $k$.)

\begin{Bem}
    The spacetimes considered above are all examples of so-called \emph{expanding spacetimes}.
    Expanding spacetimes are more generally defined by requiring that the mean curvature of the slices $\Sigma$ in $M$ is positive \cite[Definition 3]{ring24}.
    We would like to observe that, since $1/s$ is assumed to be integrable, $s$ must be unbounded as $t \to \infty$, and in particular the volume of any compact spatial set is also unbounded as $t \to \infty$. 
    Note, however, that the integrability of $1/s$ does not imply that $s$ increases near infinity. 
    Even if we assume that $s$ increases, it appears that there is still no way of controlling the mean curvature without making additional assumptions on the lapse function $N$ and the spatial metric $g$.
    On the other hand, given an arbitrary (but sufficiently fast) expanding spacetime, it is also not clear how one would choose the function $s$ to fit the definition above.
    Our notion of spacetimes of expanding type and the one of mean curvature expanding spacetimes therefore seem to intersect in a large class, but do not really coincide, and for this reason we refrain from referring to these as expanding spacetimes but instead say they are of expanding type.
\end{Bem}

\begin{center}
    \emph{Unless otherwise specified, we hereafter work exclusively with the conformally transformed spacetime $(M, \tilde h)$, and drop the tildes for notational simplicity.}
\end{center}

We can simplify matters further and define a new temporal coordinate $\tau$ via 
\begin{equation*}
    \diff\tau = \frac{1}{s}\, \diff t, \qquad \tau(0) = 0.
\end{equation*}
With respect to this coordinate, we can view $M$ as $I \times \Sigma$, where $I \subset \R$ is a connected interval with
\begin{equation*}
    T:= \sup I = \int_0^\infty \frac{\diff t }{s(t)} < \infty,
\end{equation*}
since $1/s \in L^1(0,\infty)$ by definition (note that $\inf I$ is negative but may or may not be finite).
The metric becomes
\begin{equation*}
    h = -\diff\tau \otimes\diff\tau + g_\tau,
\end{equation*}
and the vector field $\partial_\tau$ is the future-directed timelike unit normal to the hypersurfaces $\{\tau\} \times \Sigma$.
This type of coordinate $\tau$ is often called a \emph{Gaussian time coordinate} \cite{MR1828850} and we say that the corresponding spacetime is \emph{Gaussian foliated}.
Definition \ref{def-expanding-spacetime} thus ensures that the original spacetime is conformal to a Gaussian foliated spacetime whose spatial slices have uniformly bounded geometry.

One of the main strong features of Gaussian foliated spacetimes is the fact that
\begin{equation}\label{eq-e0-parallel}
    \nabla_{\partial_\tau} \partial_\tau = 0,
\end{equation}
so that the curves $\tau \mapsto (\tau,x)$ are geodesics parametrised by arclength for each $x \in \Sigma$.
To see this, choose local coordinates $x^k$ on $\Sigma$, so that $(\tau,x)$ forms a coordinate system on $I \times \Sigma$.
From the compatibility of $\nabla$ with the metric $h$ and torsion freedom, we get
\begin{equation}\label{eq-e0-parallel-proof}
    h(\nabla_{\partial_\tau} \partial_\tau, \partial_{x^k}) = -h(\partial_\tau, \nabla_{\partial_\tau}\partial_{x^k}) = -h(\partial_\tau, \nabla_{\partial_{x^k}}\partial_\tau) = -\tfrac12 \partial_{x^k} h(\partial_\tau,\partial_\tau) = 0,
\end{equation}
and similarly from metric compatibility $h(\nabla_{\partial_\tau} \partial_\tau, \partial_\tau) = 0$ (cf.\ \cite[\S 4]{MR2121740}).

\subsection{Spatial tensors}
\label{sec-spatial-tensors}
The tangent bundle of $M$ admits an orthogonal decomposition
\begin{equation}\label{eq-tangent-bundle-decomposition}
    TM = \R \partial_\tau \oplus \pi^\ast T\Sigma,
\end{equation}
where $\pi : M \to \Sigma$ is the projection $(\tau,x) \mapsto x$ and $\pi^\ast T\Sigma \cong I \times T\Sigma \to M$ is the pullback bundle which we shall refer to as the
\emph{spatial tangent bundle}, and its sections as \emph{spatial vector fields}.
The decomposition (\ref{eq-tangent-bundle-decomposition}) induces naturally a decomposition of all tensor bundles, allowing us to consider \emph{spatial tensors}, i.e.\ sections of $\pi^\ast(T^\ast \Sigma)^{\otimes k} \otimes \pi^\ast(T\Sigma)^{\otimes \ell}$.

The bundle metric $h$ on $TM$ restricts to a bundle metric on $\pi^\ast T\Sigma$ (explicitly this is $g$).
The Levi-Civita connection $\nabla$ on $(M,h)$ induces a metric-compatible connection on the spatial tangent bundle, explicitly given by the projection $\nabla^{\pi^\ast T\Sigma} = \pi_\ast \nabla$.
We note that for $X,Y \in \Gamma(\pi^\ast (T\Sigma))$,
\begin{equation*}
    \nabla^{\pi^\ast T\Sigma}_{\partial_\tau} X = \nabla_{\partial_\tau} X,
    \qquad 
    \nabla^{\pi^\ast T\Sigma}_X \, Y = D_X Y,
\end{equation*}
where the former follows since $\nabla_{\partial_\tau}X$ is purely spatial by the calculation in (\ref{eq-e0-parallel-proof}), and in the latter $D=D^{\pi^\ast T\Sigma}$ corresponds to the pullback of the Levi-Civita connection of the Riemannian slice $\Sigma_\tau$. 
In fact, we can differentiate $Y$ in the direction of $X$ either with respect to $\nabla$ or $\nabla^{\pi^\ast T\Sigma}$, and the result differs only by a normal term, i.e.\
\begin{equation}\label{eq-gauss-formula}
    \nabla_X Y = -\sff(X,Y)\, \partial_\tau + D_X Y,
\end{equation}
where $\sff(X,Y) = h(\partial_\tau, \nabla_XY)$ is the second fundamental form.
We note that $\sff$ can be computed in coordinates by
\begin{equation*}
    \sff(\partial_{x^i}, \partial_{x^j}) = -\frac12 \frac{\partial}{\partial\tau} g(\partial_{x^i},\partial_{x^i}) , 
\end{equation*}
although we will not use this formula since we will work with frames.

More generally, assume that $V \to M$ is a vector bundle equipped with a bundle metric and a connection
\begin{equation*}
    \nabla^V : \Gamma(V) \to \Gamma(T^\ast M \otimes V).
\end{equation*}
Then, via (\ref{eq-tangent-bundle-decomposition}), we can decompose the connection into temporal and spatial derivatives%
\footnote{
If we set in particular $V=\pi^\ast T\Sigma$, then $D^V$ here corresponds to the aforementioned pullback of the Levi-Civita connection of the Riemannian slices. However, if we set $V = TM$, then $D^V$ is merely the restriction of the ambient connection $\nabla$ to spatial tangent vectors. These are as already mentioned not the same, but related by (\ref{eq-gauss-formula}). \label{foot-LC-vs-spatial}
}
\begin{equation*}
    \frac{\nabla^V}{\diff\tau} := \nabla^V_{\partial_\tau} : \Gamma(V) \to \Gamma(V), 
    \qquad 
    D^V : \Gamma(V) \to \Gamma(\pi^\ast (T^\ast \Sigma) \otimes V).
\end{equation*}
Using the connection $\nabla^{\pi^\ast T\Sigma}$, these extend naturally via the Leibniz rule to $V$-valued spatial tensors (i.e.\ sections of $\pi^\ast(T^\ast \Sigma)^{\otimes k} \otimes \pi^\ast(T\Sigma)^{\otimes \ell}\otimes V$),
so e.g.\ for a one-form $\xi \in \Gamma(\pi^\ast (T^\ast \Sigma) \otimes V)$ and $X,Y \in \pi^\ast T\Sigma$, we have
\begin{equation*}
    \left( \frac{\nabla^V \xi}{\diff\tau} \right)(X) = \frac{\nabla^V}{\diff\tau} (\xi(X)) - \xi(\nabla_{\partial_\tau} X),
    \qquad
    (D^V_X \xi)(Y) = D^V_X(\xi(Y)) - \xi(D_XY).
\end{equation*}
The connection wave operator acting on sections of $V$ satisfies
\begin{equation*}
    \Box^V = -(\nabla^V)^\ast \nabla^V  = -\frac{(\nabla^V)^2}{\diff\tau^2} + 3H \frac{\nabla^V}{\diff\tau}  - (D^V)^\ast D^V,
\end{equation*}
where $H = \frac13\tr_g(\sff)$ is the mean curvature of $\Sigma$.

\subsection{Adapted frames}
\label{sec-adapted-frame}
As previously in the paper, it will be convenient to work with frames rather than coordinates.
Setting $e_0 = \partial_\tau$ we obtain a unit timelike normal to $\Sigma$, so we only need to choose orthonormal spatial vector fields to complete the frame.
If $e_i \in \Gamma(\pi^\ast T\Sigma)$ are orthonormal, then using metric compatibility we can also compute (independently of the choice of orthonormal spatial vector fields)
\begin{equation*}
    h(\nabla_{e_i} e_0, e_0) = 0, \qquad
    h(\nabla_{e_i} e_0, e_k) = - h(e_0, \nabla_{e_i} e_k) = -\sff(e_i, e_k) = -\sff_{ik},
\end{equation*}
which implies that
\begin{equation*}
    \nabla_{e_i} e_0 = -\tensor{\sff}{_i^k} e_k.
\end{equation*}
Thus, the only remaining unknown connection symbols are the spatial components of $\nabla_{e_0} e_i$ (the temporal component is always zero by metric compatibility).
In fact, we can construct a particularly nice adapted frame for which $\nabla_{e_0} e_i=0$, simply by choosing an orthonormal frame $(e_i)_{\tau_0}$ on the tangent space $T\Sigma_{\tau_0}$ for fixed $\tau_0$, and parallel transporting it along the curves $\tau \mapsto (\tau, x)$ with tangent vector $\partial_\tau = e_0$ (cf.\ \cite[Lemma 3.1]{MR3830277}).
We thus obtain a semi-orthonormal frame $e_\mu$ for $(M,h)$, for which the connection operates as
\begin{equation*}
    \nabla_{e_0} e_0 = \nabla_{e_0} e_i = 0, \qquad \nabla_{e_i} e_0 = -\tensor{\sff}{_i^k} e_k, \qquad \nabla_{e_i} e_j = -\sff_{ij} \, e_0 + D_{e_i} e_j,
\end{equation*}
where $D$ is again the Levi-Civita connection on $\Sigma_\tau$.

Since this construction only depends on the choice of spatial frame, we see that we can cover the entire spacetime by a collection of such frames.
Throughout the remainder of this article, we work with such adapted semi-orthonormal frames, unless explicitly specified otherwise. 

For later convenience, we also note that the induced coframe satisfies
\begin{equation}\label{eq-connection-coframe}
    \nabla_{e_0} e^0 = \nabla_{e_0} e^i = 0, \qquad \nabla_{e_i} e^0 = \sff_{ik} e^k, \qquad \nabla_{e_i} e^j = \tensor{\sff}{_i^j} e^0 + D_{e_i} e^j,
\end{equation}
as well as
\begin{equation*}
    \diff e^0 = 0, \qquad \diff e^k = -\tensor{\sff}{_i^k}\, e^0 \wedge e^i - e^k(D_{e_i} e_j) \, e^i \wedge e^j.
\end{equation*}

\subsection{Curvature tensor}

Let us also express the curvature tensor of $(M,h)$ in terms of the
second fundamental form and the curvature tensor of $(\Sigma, g_\tau)$.

\begin{Lem}\label{lem-curvature-generalized-cyl}
    The non-zero components (up to symmetries) of the Riemann curvature tensor of $(M,h)$ are given by
    \begin{align*}
        R(h)_{k0i0} &= -\left(\tfrac{\nabla\sff}{\diff\tau}\right)_{ki} + \tensor{\sff}{_k^\ell} \sff_{\ell i}, \\[0.1cm]
        R(h)_{kij0} &= (D\sff)_{kij} - (D\sff)_{ijk},\\[0.1cm]
        R(h)_{ijk\ell} &= R(g)_{ijk\ell} + \sff_{ik}\sff_{j\ell} - \sff_{jk}\sff_{i\ell}.
    \end{align*}
    The Ricci tensor satisfies
    \begin{align*}
        R(h)_{00} &= \tr\left(\tfrac{\nabla\sff}{\diff\tau}\right) - 2|\sff|^2,\\[0.1cm]
        R(h)_{0k} &= 3(\diff H)_k + (D^\ast \sff)_k,\\[0.1cm]
        R(h)_{ik} &= R(g)_{ik}  - \left(\tfrac{\nabla\sff}{\diff\tau}\right)_{ik} + 2\tensor{\sff}{_i^\ell} \sff_{\ell k} - 3H \sff_{ik}. 
    \end{align*}
\end{Lem}

\begin{Bem}
    Recall that we use the convention (\ref{eq-riemann-convention}) for the Riemann tensor.
    Note that these formulae hold in any semi-orthonormal frame with $e_0 = \partial_\tau$.
    The reader can compare this result with \cite[Appendix B.2]{MR2527641} and \cite[Proposition 4.1]{MR2121740} (although one should note that certain signs are different here than in the Riemannian setting, due to the Lorentzian character of the normal field $e_0$). 
\end{Bem}

\subsection{Spatial principal bundle}\label{sec-spatial-principal-bundle}

Let $P \to M$ be a principal $G$-bundle. 
Since we have the foliation $M = I \times \Sigma$, we can also foliate the principal bundle appropriately, cf.\ \cite[\S 1]{MR1604914}.
%,
To see this, fix $\tau_0 \in I$ and let $P^\Sigma = \iota^\ast P$ be the pullback bundle under the embedding $\iota : \Sigma \to M, \; x\mapsto (\tau_0,x)$.
Then $P^\Sigma$ is a principal $G$-bundle over $\Sigma$ with fibers $P^\Sigma_x = P_{(\tau_0,x)}$.
Now pull $P^\Sigma$ back to $M$ via the projection $\pi : M \to \Sigma$ to obtain a principal $G$-bundle $\pi^\ast P^\Sigma \cong I \times P^\Sigma \to M \cong I \times \Sigma$, equipped with the trivial bundle projection and principal $G$-action
\begin{equation*}
    \pi^P(\tau,p) = (\tau,\pi^{P^\Sigma}(p)), \qquad (\tau,p)g = (\tau,pg),
\end{equation*}
where $(\tau,p) \in I \times P^\Sigma$.
Then $P \cong \pi^\ast P^\Sigma$ are diffeomorphic as principal $G$-bundles.
To see this, equip $P$ with any smooth connection%
\footnote{A connection always exists \cite[Satz 3.4]{BAUM}. Here, the choice of connection is immaterial and is used merely for the construction of the bundle diffeomorphism.}
and denote the parallel transport of an element $p_x \in P_{x}^\Sigma = P_{(\tau_0,x)}$ along the curve $\tau \mapsto (\tau,x)$ by $\Phi_\tau(p_x) \in P_{(\tau,x)}$.
Then, the mapping
\begin{equation*}
    \pi^\ast P^\Sigma \to P, \qquad (\tau,p) \mapsto \Phi_\tau(p)
\end{equation*}
is a principal bundle diffeomorphism (compatibility with projections holds trivially and $G$-equivariance follows from the well-known properties of parallel transport, cf.\ \cite[\S 5.8]{MR3837560}).
Hence, we can identify $P \cong \pi^\ast P^\Sigma \cong I \times P^\Sigma$ and we will henceforth always make this identification without mentioning it explicitly.

Although the present article focuses on the more global aspects of the theory, let us still say a few words about local gauges on $P \cong \pi^\ast P^\Sigma$. 
Sections of $\pi^\ast P^\Sigma$ are of all the form 
$$(\tau, x)\mapsto (\tau, s^\Sigma(\tau,x)),$$ 
where $s^\Sigma(\tau,\cdot) : U \subset \Sigma \to P^\Sigma$ is a one-parameter family of local sections of $P^\Sigma$.
In fact, it turns out to be most convenient to keep $s^\Sigma$ constant with respect to $\tau$, in which case the corresponding section of $\pi^\ast P^\Sigma$ will be called \textit{adapted}.
An atlas of sections of $P^\Sigma$ then naturally induces an atlas of adapted sections of $P$, so one can without loss of generality work with adapted gauges.

A connection $\omega$ on $P$ is always of the form
\begin{equation*}
    \omega_{} = \diff \tau \otimes \alpha + \Pi^\ast \sigma,
\end{equation*}
where $\Pi:P\to P^\Sigma$ is the natural projection, $\alpha \in \Gamma(\Ad P)$ and $\sigma=\sigma_\tau$ is a smooth one-parameter family of connections on $P^\Sigma$. 
We note that if $\sigma$ is independent of $\tau$ (so just a fixed connection on $P^\Sigma$), then for any adapted gauge the spatial part of $s^\ast \omega$ is independent of $\tau$ and depends only on $\sigma$ and the choice of gauge $s^\Sigma$ on $P^\Sigma$, while the temporal part is given simply by $\alpha$.%
\footnote{For non-adapted gauges, the temporal part would also depend on $\sigma$ and $s^\Sigma$, i.e.\ $(s^*\omega)(\partial_\tau) = \alpha+\sigma(\diff s^\Sigma(\partial_\tau)).$}

\section{Standard Model energy estimates}
\label{sec-energy-estimates}

We now wish to define the $k$-th total energy of the coupled system (\ref{eq-yang-mills-wave}--\ref{eq-dirac-wave}) on Gaussian foliated spacetimes $(M,h)$. 
At first, we will define the standard Sobolev spaces on vector bundles and also the Sobolev spaces of connections.
However, we need to be a bit careful as $F_\omega$ and $\Psi$ are sections respectively of $\Lambda^2 M \otimes \Ad P$ and $\Fer_+$, whose natural geometric inner product is not positive-definite, so we will also discuss each sector separately in more detail.

\subsection{Sobolev spaces on vector bundles}
\label{sec-sobolev-vector-bundle}

Assuming that a vector bundle $V \to M$ comes with a positive-definite (Hermitian) inner product and a (not necessarily metric compatible) connection, we can define a positive-definite $H^k$-norm at $\tau\in I$ of a given spacelike compactly supported section $\xi\in\Gamma(V)$ as
\begin{equation}\label{eq-Hk-norm}
    \Vert \xi \Vert_{H^k(g_\tau)}^2 = \sum_{\ell=0}^k \Vert D^\ell \xi \Vert_{L^2(g_\tau)}^2 = \sum_{\ell=0}^k \int_{\Sigma} |D^\ell \xi|^2 \dv_{g_\tau},
\end{equation}
where $D = D^V$ is as defined in \S \ref{sec-spatial-tensors}.
The analogous definition applies also to functions.
Note that the definition of the Sobolev norm depends on the choice of the connection $D$, which will be somewhat important later.
We also define the space $C^\ell(I, H^k(V))$ for $\ell<k$ as the closure of the space of spacelike compactly supported sections $\xi\in\Gamma(V)$ under the norm
\begin{equation*}
    \sup_{\tau \in I} \left\Vert \frac{\nabla^\ell \xi}{\diff\tau^\ell} \right\Vert_{H^k(g_\tau)} < \infty.
\end{equation*}
We recall that a Riemannian manifold $(\Sigma,g)$ is said to have \emph{bounded geometry} of order $k$ if its injectivity radius is lower bounded by a positive constant, and the covariant derivatives up to order $k$ of the Riemann tensor are uniformly bounded (i.e.\ with respect to the supremum norm). Similarly, a connection $D^V$ on a vector bundle $V$ is said to have bounded geometry of order $k$ if the derivatives up to order $k$ of the corresponding curvature tensor $R^V$ are uniformly bounded. These assumptions are sufficient to unlock the usual Sobolev inequalities (up to order $k$) on Riemannian manifolds (compact or non-compact) and vector bundles over them \cite{MR2343536, MR1481970}.
In particular, under the assumptions of bounded geometry, the above defined $H^k$-space coincides with the space of measurable sections for which the norm is finite \cite{MR2343536}, see also \cite{MR3126616}.

\subsection{Sobolev spaces of connections on principal bundles}
\label{sec-sobolev-connections}

Connections on principal bundles over compact Riemannian manifolds were studied in the seminal paper of Uhlenbeck \cite{MR648356}.
The more involved non-compact setting was described in the works of Eichhorn and Heber \cite{MR1175322, MR1458665}.
Let us describe their results.
As already hinted in the previous section, the appropriate assumption for working with Sobolev spaces is that of bounded geometry, i.e.\ we consider the space
\begin{equation*}
    \mathcal{C}^k(P^\Sigma) = \left\{ \sigma \in \Omega^1(P^\Sigma,\g) \,\mid\, \sigma \text{ smooth connection}, \; \sum_{\ell=0}^k\sup_\Sigma |D_\sigma^\ell F_\sigma| < \infty \right\}
\end{equation*}
for an $m$-dimensional Riemannian manifold $(\Sigma,g)$ of bounded geometry.
For a connection $\sigma \in \mathcal{C}^k(P^\Sigma)$ and Sobolev spaces on associated vector bundles as defined in the previous \S \ref{sec-sobolev-vector-bundle}, the usual Sobolev apparatus therefore holds true (up to order $k$).
One would now also like to define the notion of Sobolev regular connections.
It was shown in \cite{MR1175322} that the natural intrinsic way of doing this is to equip the space $\mathcal{C}^k(P^\Sigma)$ with the uniform topology induced by the fundamental system given by
\begin{equation*}
    U_\delta = \left\{ (\sigma, \sigma') \in \mathcal{C}^k(P^\Sigma) \times \mathcal{C}^k(P^\Sigma) \,\mid\, \sup_\Sigma | \sigma - \sigma' |_{H^k_\sigma} < \infty \right\},
\end{equation*}
where we view $\sigma - \sigma'$ as an $\Ad P^\Sigma$-valued one-form on $\Sigma$ and the $H^k_\sigma$ norm is taken with respect to the connection induced by $\sigma$.
The technical details of the matter are somewhat involved, but for the purposes of this work it will be sufficient to note that the closure of this topological space $\mathcal{H}^{k}(P^\Sigma) := \overline{\mathcal{C}^k(P^\Sigma)}$ with $k > m/2+1$ can be written as a topological sum of components
\begin{equation}\label{eq-sobolev-connections-riemannian}
    \mathcal{H}^{k}(P^\Sigma) = \coprod_{i \in \mathscr I} \comp(\sigma_i), \qquad  \comp(\sigma_i) = \sigma_i + H^{k}_{\sigma_i}(T^\ast \Sigma \otimes \Ad P^\Sigma)
\end{equation}
for some connections $\sigma_i \in \mathcal{C}^k(P^\Sigma)$ of bounded geometry and some index set $\mathscr I$.
Note that, in the compact case, the index set $\mathscr I$ has only a single element, and one recovers the classical definition of Uhlenbeck \cite{MR648356}.
Elements of $\mathcal{H}^{k}(P^\Sigma)$ are then what should be considered as $H^k$-connections.

At spacetime level, we wish to consider connections $\omega$ on $P \cong \pi^\ast P^\Sigma$ such that $\iota_\tau^\ast \omega$ defines an element of $\mathcal{H}^{k}(P^\Sigma)$ uniformly for each $\tau \in I$, restricted to a single topological component.
More precisely, we define
\begin{equation*}
    C^0(I, \mathcal{H}^{k}(P^\Sigma)) = \coprod_{i\in\mathscr{I}} \left( \omega_i + C^0(I, H^k_{\omega_i}(\pi^\ast T^\ast \Sigma \otimes \Ad P)) \right),
\end{equation*}
where $\omega_i = \Pi^\ast \sigma_i$, while $\mathscr{I}$ and $\sigma_i \in \mathcal{C}^k(P^\Sigma)$ are as in (\ref{eq-sobolev-connections-riemannian}).
Note that this definition is independent of the choice of $\sigma_i$ in (\ref{eq-sobolev-connections-riemannian}).
Indeed, if we choose another smooth connection of bounded geometry $\tilde\sigma_i \in \comp(\sigma_i) \subset \mathcal{H}^{k}(P^\Sigma)$, then $\theta_i = \sigma_i-\tilde\sigma_i \in H^k_{\sigma_i}(T^\ast \Sigma \otimes \Ad P^\Sigma)$ and $\pi^\ast \theta_i \in C^0(I, H^k_{\omega_i}(\pi^\ast T^\ast \Sigma \otimes \Ad P))$ trivially since the spatial metrics $g_\tau$ have uniformly bounded geometry, so that the $H^k_{\sigma_i}$- and $H^k_{\tilde\sigma_i}$-norms are equivalent.

\subsection{Higgs field energy}
The natural geometric inner product on the Higgs field bundle $\Hig$ is already positive-definite, so we take the $k$-th total energy of $\Phi$ to be
\begin{equation}\label{eq-k-total-energy}
    \mathfrak{E}_0(\Phi) = \Vert \Phi \Vert_{L^2}^2,
    \qquad
    \mathfrak{E}_{k+1}(\Phi; \omega) 
    = \left\Vert \frac{\nabla_\omega\Phi}{\diff\tau} \right\Vert_{H^k_\omega}^2 + \Vert \Phi \Vert_{H^{k+1}_\omega}^2,
\end{equation}
where the $H^k_\omega$ norms are taken with respect to the affine connection on $\Hig$ induced by $\omega$.

\subsection{Yang-Mills energy}\label{sec-yang-mills-energy-definition}
We decompose the curvature form $F_\omega$ into its electric and magnetic parts
\begin{align*}
    E_\omega &= e_0 \iprod F_\omega = e^i \otimes F_{0i} \in \Gamma(\pi^\ast T^\ast \Sigma \otimes \Ad P),\\[0.1cm]
    B_\omega &= F_\omega - e^0 \wedge E_\omega = e^i \wedge e^j \otimes \tfrac12 F_{ij} \in \Gamma(\pi^\ast \Lambda^2 \Sigma \otimes \Ad P),
\end{align*}
in terms of which we have the orthogonal splitting $F_\omega = e^0 \wedge E_\omega + B_\omega$.
Note that this splitting does not depend on the choice of frame, since $e_0 = \partial_\tau$ is the unique globally defined future-directed unit timelike vector.
Since $E_\omega$ and $B_\omega$ are now sections of Hermitian vector bundles over $M$, we can define the total positive-definite energy of $F_\omega$ as 
\begin{equation*}
    \mathfrak{E}^+_k(F_\omega; \omega) := \mathfrak{E}_k(E_\omega; \omega) + \mathfrak{E}_k(B_\omega; \omega),    
\end{equation*}
where on the right-hand side the energies are defined analogously as in (\ref{eq-k-total-energy}).%
\footnote{
Here it is really important to view $E_\omega$ and $B_\omega$ as sections of $\pi^\ast T^\ast \Sigma \otimes \Ad P$ and $\pi^\ast \Lambda^2 \Sigma \otimes \Ad P$ rather than $T^\ast M \otimes \Ad P$ and $\Lambda^2 M \otimes \Ad P$, since otherwise their spatial derivatives would also have a timelike component 
(cf.\ Footnote \ref{foot-LC-vs-spatial} on p.\ \pageref{foot-LC-vs-spatial}) 
and hence the higher energies would not be positive-definite.
}

\subsection{Dirac energy}
For spinorial quantities, we again face the problem that the geometric inner product is indefinite, and in particular we have $|\Psi|^2 = 0$ for all $\Psi \in \Gamma(\Fer_+)$.
We circumvent this issue by redefining the inner product by multiplying one of the factors by the unit timelike vector $e_0$.
To avoid confusion, we will conventionally denote the geometric inner product and corresponding norm by single lines as before, and the positive-definite one by double lines, i.e.\ we set
\begin{equation}\label{eq-spinor-positive-inner-prod}
    \Langle \Psi_1, \Psi_2 \Rangle := \langle \Psi_1, e_0 \cdot \Psi_2 \rangle, \qquad \Psi_1,\Psi_2 \in \Gamma(\Fer_+).
\end{equation}
Note that upon fixing a local spin frame, the spinorial factor of $\Langle\cdot,\cdot\Rangle$ just corresponds to the standard Hermitian inner product on $\C^4$.
The inner product $\Langle \cdot, \cdot \Rangle$ on $\Fer$ then induces an $L^2$-inner product of $\Psi$ on $M$ and we can define the positive-definite energy $\mathfrak{E}^+_k(\Psi; \omega)$ by the same formulae (\ref{eq-Hk-norm}--\ref{eq-k-total-energy}) but with respect to $\Langle \cdot, \cdot \Rangle$.
Due to the modification of the inner product, we also lose certain geometric properties (most importantly compatibility of connections with the inner products).

\subsection{Total energy of the Standard Model}
\label{sec-sm-energy}

We now wish to define the $k$-th total energy of the system as the sum of the above defined energies.
However, one expects that the contributions coming from the different sectors can also have different regularities, so one should also take care of the orders of the energies.
We therefore define the total energy of the Standard Model system as
\begin{equation}\label{eq-total-energy-sm}
    \mathfrak{E}_{k_1, k_2, k_3} := \mathfrak{E}^+_{k_1} (F_\omega; \omega) + \mathfrak{E}_{k_2}(\Phi; \omega) + \mathfrak{E}^+_{k_3}(\Psi; \omega),
\end{equation}
where $k_1,k_2,k_3$ are positive integers.%
\footnote{In fact, we will only need the case $k_1 = k_2 = k_3$ later, but we wanted to emphasize that there are other interactions between regularities for which the energy can be estimated.}
Note that this energy is a geometric invariant in the sense that if $f : P \to P$ is a bundle automorphism, then
\begin{equation*}
    \mathfrak{E}^+_{k_1} (F_\omega; \omega) + \mathfrak{E}_{k_2}(\Phi; \omega) + \mathfrak{E}^+_{k_3}(\Psi; \omega) =
    \mathfrak{E}^+_{k_1} (fF_\omega; f\omega) + \mathfrak{E}_{k_2}(f\Phi; f\omega) + \mathfrak{E}^+_{k_3}(f\Psi; f\omega).
\end{equation*}

The main goal now is to show the following a priori energy estimate.

\begin{Satz}\label{thm-total-energy-estimate}
    Let $k > 3/2$ be an integer.
    Let $(M,h)$ be a Gaussian foliated spacetime, and let $P\to M$ be a principal bundle equipped with a connection $\omega \in C^0([0,\infty), \mathcal{C}^k(P^\Sigma))$ of uniformly bounded geometry of order $k+1$.    
    Suppose that $(\omega,\Phi,\Psi)$ is a conformal Standard Model triplet.
    Then, for $0 \leq a \leq c \leq 1$, the total energy 
    $$\mathfrak{E}_{k+a, k+1, k+c} = \mathfrak{E}^+_{k+a} (F_\omega; \omega) + \mathfrak{E}_{k+1}(\Phi; \omega) + \mathfrak{E}^+_{k+c}(\Psi; \omega)$$ 
    is uniformly bounded for $\tau \in (0,T)$ if it is sufficiently small initially, and the bound depends only on $k$, $T$,
    $$\sup_{\tau\in(0,T)}\Vert {\Riem_g} \Vert_{C^{k}}, \quad \sup_{\tau\in(0,T)}\Vert \sff \Vert_{C^{k+1}}, \quad \sup_{\tau\in(0,T)}\left\Vert \frac{\nabla\sff}{\diff\tau} \right\Vert_{C^k},$$
    and the Standard Model couplings $\lambda, \rho_*, \chi_*, \Y$.
\end{Satz}

To show this, we need to estimate the energy of each term separately, as follows.

\begin{Prop}\label{prop-total-energy-inequality}
    Suppose that $(\omega, \Phi, \Psi)$ is a conformal Standard Model triplet.
    If $k > 3/2$ and $a,b,c \in \{0,1\}$, then
    \begin{align}
        \label{eq-energy-estimate-yang-mills-total}
        &\frac{\diff}{\diff\tau} \mathfrak{E}^+_{k+a}(F_\omega; \omega) 
        \lesssim 
        C(k,g,\sff) \,\mathfrak{E}_{k+a, \circ, \circ}
        + \mathfrak{E}_{k+a, k+1, k+a}^{\frac32}
        + C(k,\sff) \,\mathfrak{E}_{k+a, k, \circ}^2,
        \\[0.1cm]
        \label{eq-energy-estimate-higgs-total}
        &\frac{\diff}{\diff\tau} \mathfrak{E}_{k+b}(\Phi; \omega) 
        \lesssim 
        C(k,g,\sff) \, \mathfrak{E}_{\circ, k+b, \circ}  + \mathfrak{E}_{k+b-1,k+b,k}^{\frac32} + |\lambda| \mathfrak{E}_{\circ,k+b,\circ}^2,
        \\[0.1cm]
        \label{eq-energy-estimate-dirac-total}
        &\frac{\diff}{\diff\tau} \mathfrak{E}^+_{k+c}(\Psi; \omega)
        \lesssim
        C(k,g,\sff) \, \mathfrak{E}_{\circ, \circ, k+c} + C(k,\sff) \, \mathfrak{E}_{k+c-1,k+c,k+c}^{\frac32} + \mathfrak{E}_{\circ,k+c,k+c}^2,
    \end{align}
    where $\circ$ is used to denote parameters that can be chosen freely.
\end{Prop}

Proposition \ref{prop-total-energy-inequality} will be proven in the remainder of the section, but first let us deduce Theorem \ref{thm-total-energy-estimate} from it.

\begin{proof}[Proof of Theorem \ref{thm-total-energy-estimate}]
    Let us denote $\mathfrak{E} = \mathfrak{E}_{k+a, k+1, k+c}$ for brevity.
    The estimates (\ref{eq-energy-estimate-yang-mills-total}--\ref{eq-energy-estimate-dirac-total}) imply that
    \begin{equation*}
        \frac{\diff}{\diff\tau} \mathfrak{E} \leq C(k,g,\sff) \, \left( \mathfrak{E} + \mathfrak{E}^{\frac32} + \mathfrak{E}^2 \right),
    \end{equation*}
    where we need to set $b=1$ in (\ref{eq-energy-estimate-higgs-total}) because the second term in (\ref{eq-energy-estimate-yang-mills-total}) contains a Higgs energy of order $k+1$ (regardless of whether $a=0$ or $a=1$) and we need $a \leq c$ since the same term contains a Dirac energy of order $k+a$.
    As long as $\mathfrak{E} \leq 1$, we therefore have the inequality
    \begin{equation*}
        \frac{\diff}{\diff\tau} \mathfrak{E} \leq C(k,g,\sff) \, \mathfrak{E},
    \end{equation*}
    and hence also
    \begin{equation*}
        \mathfrak{E}(\tau) \leq \mathfrak{E}(0) \, \exp\big(C(k,g,\sff)\tau\big), \quad \tau \in (0,T).
    \end{equation*}
    Thus if we assume $\mathfrak{E}(0) < \exp\big(-C(k,g,\sff)T\big)$, the formula above shows that $\mathfrak{E}$ stays in the region $\mathfrak{E} < 1$ for all $\tau \in (0,T)$.
\end{proof}

\subsection{Some preliminary estimates}

In this section we discuss some estimates that will be used in the derivation of Proposition \ref{prop-total-energy-inequality}.

\subsubsection{Evolution of Sobolev norms}

We will often want to evolve the $H^k$-norm of a bundle section.
To that end, we note the following.

\begin{Lem}\label{lem-Hk-norm-evolution}
    Let $(M,h)$ be a Gaussian foliated spacetime and let $V \to M$ be a Hermitian vector bundle with a compatible connection $\nabla^V$.
    If $\xi\in \Gamma(V)$ is spacelike compactly supported, then
    \begin{equation*}
        \frac{\partial}{\partial \tau} \Vert \xi \Vert_{H^k(V)}^2
        \lesssim 
        \left(\left\Vert \frac{\nabla \xi}{\diff \tau} \right\Vert_{H^k(V)} + \Vert R^V(\partial_\tau, \cdot)\xi \Vert_{H^{k-1}(V)} + \Vert \sff \Vert_{C^k} \Vert \xi\Vert_{H^k(V)} \right) \Vert \xi \Vert_{H^k(V)},
    \end{equation*}
    where the curvature term is understood to be zero if $k=0$.
\end{Lem}
\begin{proof}
    We calculate
    \begin{align*}
        \frac12 \frac{\partial}{\partial \tau} \Vert \xi \Vert_{H^k(V)}^2
        &= \frac12\sum_{\ell=0}^k \frac{\partial}{\partial \tau} \Vert D^\ell \xi \Vert_{L^2(V)}^2
        = \sum_{\ell=0}^k \int_\Sigma \left( \left\langle \frac{\nabla}{\diff \tau} D^\ell \xi, D^\ell \xi \right\rangle - \frac32 H |D^\ell\xi|^2 \right) \dv_g,
    \end{align*}
    where the second term comes from the evolution of the volume form and we recall that $H = \tfrac13 \tr_g\sff$ is the mean curvature.
    Hence,
    \begin{align*}
        \frac{\partial}{\partial \tau} \Vert \xi \Vert_{H^k(V)}^2
        &\lesssim \sum_{\ell=0}^k \left\Vert \frac{\nabla}{\diff \tau}  D^\ell \xi \right\Vert_{L^2(V)} \Vert \xi \Vert_{H^k(V)} + \Vert H \Vert_{C^0}\Vert \xi\Vert_{H^k(V)}^2\\
        &\lesssim \left\Vert \frac{\nabla\xi}{\diff \tau} \right\Vert_{H^k(V)} \Vert \xi \Vert_{H^k(V)} + \sum_{\ell=0}^k \left\Vert \left[\frac{\nabla}{\diff \tau},  D^\ell\right] \xi \right\Vert_{L^2(V)} \Vert \xi \Vert_{H^k(V)} + \Vert \sff \Vert_{C^0}\Vert \xi\Vert_{H^k(V)}^2.
    \end{align*}
    Now it only remains to deal with the commutator term.
    To this end we note that for $\xi \in \Gamma(T^\ast \Sigma^{\otimes p} \otimes V)$, $X \in \Gamma(T\Sigma)$ and $Y\in \Gamma(T\Sigma^{\otimes p})$, we have
    \begin{equation}\label{eq-commutator-temporal-spatial-explicit}
        \left(\left[\frac{\nabla}{\diff\tau}, D\right]\xi\right)(X, Y) = R^V(\partial_\tau, X) (\xi(Y)) + (D\xi)(\sff(X), Y) - \xi((D\sff)(X,Y)),
    \end{equation}
    where we interpret $$(D\sff)(X, Y_1 \otimes \ldots \otimes Y_p) = (D\sff)(X, Y_1) \otimes Y_2 \otimes \ldots\otimes Y_p + \ldots + Y_1 \otimes\ldots \otimes Y_{p-1}\otimes (D\sff)(X,Y_p).
    $$   
    For the higher powers, we use the general identity
    \begin{equation}\label{eq-commutator-power}
        [B,A^k]\xi=\sum_{\ell=0}^{k-1}A^{k-\ell-1}\big([B,A]A^{\ell}\xi\big), 
    \end{equation}
    which holds for arbitrary differential operators \(A,B\).
    This allows us to write 
    \begin{equation*}
        \left[\frac{\nabla}{\diff\tau}, D^\ell\right]\xi=
        \sum_{j=0}^{\ell}D^{\ell-j-1}\big(R^V(\partial_\tau,\cdot)\big)\ast
        D^{j}\xi
        +\sum_{j=0}^{\ell} D^{\ell-j}\sff\ast D^{j}\xi,
    \end{equation*}
    where we use the $\ast$-notation as described at the end of \S \ref{sec-conventions},
    and hence we can estimate
    \begin{equation}\label{eq-commutator-curvature-estimate}
        \sum_{\ell=0}^k\left\Vert\left[\frac{\nabla}{\diff\tau}, D^{\ell}\right]\xi\right\Vert_{L^2}
        \lesssim 
        \left\Vert R^V(\partial_\tau,\cdot)\xi \right\Vert_{H^{k-1}}
        +\left\Vert\sff\right\Vert_{C^{k}}\left\Vert\xi\right\Vert_{H^{k}},
    \end{equation}
    which gives the result.
\end{proof}

\subsubsection{General energy estimate}

We have defined the energy of a section of a Hermitian vector bundle in (\ref{eq-k-total-energy}).
Provided that the bundle connection is compatible with the inner product used in this definition, one can derive a general energy estimate as follows.

\begin{Prop}\label{prop-general-energy-estimate}
    Let $(M,h)$ be a Gaussian foliated spacetime and let $V \to M$ be a Hermitian vector bundle with a compatible connection $\nabla^V$.
    Then if $\xi \in \Gamma(V)$ is   spacelike compactly supported, we have
    \begin{align*}
        \frac{\diff}{\diff\tau} \mathfrak{E}_0(\xi) 
        &\lesssim \mathfrak{E}_0(\xi)^{\frac12}\mathfrak{E}_1(\xi)^{\frac12} + \Vert H \Vert_{C^0} \, \mathfrak{E}_0(\xi),
    \end{align*}
    and if $\xi \in C^1(I, H^{k+1}(V))$ for $k \geq 0$, then
    \begin{align*}
        \frac{\diff}{\diff\tau} \mathfrak{E}_{k+1}(\xi) 
        \lesssim&\; \mathfrak{E}_0(\xi)^{\frac12}\mathfrak{E}_1(\xi)^{\frac12} + \Vert {\Riem_g} \Vert_{C^{k}} \mathfrak{E}_{k}(\xi) + \Vert \sff \Vert_{C^{k+1}} \, \mathfrak{E}_{k+1}(\xi)
        \\[0.1cm]\nonumber
        &+ \left(\Vert \Box^V \xi \Vert_{H^{k}} 
        + \mathfrak{R}^V_k(\xi)
        \right) \mathfrak{E}_{k+1}(\xi)^{\frac12},
    \end{align*}
    where we denote
    \begin{equation*}
        \mathfrak{R}^V_k(\xi) = 
        \begin{cases}
            \left\Vert R^V(\partial_\tau, \cdot)\xi \right\Vert_{L^2}, 
            & \text{if } k=0,\\[0.25cm]
            \left\Vert R^V(\partial_\tau, \cdot)\tfrac{\nabla\xi}{\diff\tau} \right\Vert_{H^{k-1}}
            +
            \left\Vert R^V\xi \right\Vert_{H^{k}},
            +
            \left\Vert R^V(\partial_\tau, \cdot)\xi \right\Vert_{H^{k}}, 
            & \text{if } k\geq 1,
        \end{cases}
    \end{equation*}
    and $R^V(\partial_\tau, \cdot)\xi$ and $R^V\xi$ are viewed as $V$-valued spatial tensors, i.e.\ sections respectively of the bundles $\pi^\ast T^\ast\Sigma \otimes V$ and $\pi^\ast \Lambda^2 \Sigma \otimes V$.
\end{Prop}

\begin{Bem}
    For our purposes, it will be entirely sufficient to trivially estimate the first three terms on the right-hand side, so that
    \begin{equation}\label{eq-energy-estimate-general}
        \frac{\diff}{\diff\tau} \mathfrak{E}_{k+1}(\xi) 
        \lesssim \left( 1 + \Vert {\Riem_g} \Vert_{C^{k}} + \Vert \sff \Vert_{C^{k+1}} \right)\mathfrak{E}_{k+1}(\xi)
        + \left(\Vert \Box^V \xi \Vert_{H^{k}} 
        + \mathfrak{R}^V_k(\xi)
        \right) \mathfrak{E}_{k+1}(\xi)^{\frac12}.
    \end{equation}
\end{Bem}

\begin{proof}
For the zeroth energy we have the trivial bound
\begin{align*}
    \frac{\diff}{\diff\tau} \mathfrak{E}_0(\xi)
    &= 2\int_\Sigma \left\langle \frac{\nabla\xi}{\diff\tau}, \, \xi \right\rangle \dv_g
    - 3 \int_\Sigma H|\xi|^2 \dv_g
    \lesssim  \left\Vert\frac{\nabla\xi}{\diff\tau}\right\Vert_{L^2} \Vert \xi \Vert_{L^2} + \Vert H \Vert_{C^0} \Vert \xi \Vert_{L^2}^2,
\end{align*}
which implies the desired inequality.
To estimate the higher order energies, let us also define
\begin{equation}\label{eq-k-energy}
     \mathcal{E}_{k+1}(\xi)(\tau)
     = \left\Vert D^k\tfrac{\nabla\xi}{\diff\tau} \right\Vert_{L^2}^2 + \Vert D^{k+1}\xi \Vert_{L^2}^2
     = \int_{\Sigma} \left( \left|D^k\tfrac{\nabla\xi}{\diff\tau}\right|^2 + |D^{k+1}\xi|^2 \right) \dv_{g_\tau},
\end{equation}
so that we can write
\begin{equation*}
    \mathfrak{E}_{k+1}(\xi) 
    = \mathfrak{E}_0(\xi) + \sum_{\ell=1}^{k+1} \mathcal{E}_\ell(\xi).
\end{equation*}
A calculation shows (cf.\ \cite[Proposition 3.3]{MR3830277} but note the last term is missing)
\begin{align*}
    \frac{\diff}{\diff\tau} \mathcal{E}_{k+1}(\xi) &=
    -2\int_\Sigma \left\langle  D^k \Box\xi, \, D^k \tfrac{\nabla\xi}{\diff\tau} \right\rangle \dv_g 
    - 3\int_\Sigma \left(|D^{k+1}\xi|^2 - \left|D^k\tfrac{\nabla\xi}{\diff\tau}\right|^2 \right)H \dv_g\\
     &\quad + 2\int_\Sigma \left\langle \left[\tfrac{\nabla}{\diff\tau}, D^k\right] \tfrac{\nabla\xi}{\diff\tau} + [D^*D, D^k]\xi, \, D^k \tfrac{\nabla\xi}{\diff\tau} \right\rangle \dv_g  
     + 2\int_\Sigma \left\langle \left[\tfrac{\nabla}{\diff\tau}, D^{k+1}\right]\xi, \, D^{k+1}\xi \right\rangle \dv_g\\
     &\quad + 6\sum_{\ell=1}^k \binom{k}{\ell} \int_\Sigma \left\langle D^\ell H \otimes D^{k-\ell}\tfrac{\nabla\xi}{\diff\tau}, \, D^k\tfrac{\nabla\xi}{\diff\tau} \right\rangle \dv_g.
\end{align*}
Hölder's inequality permits the estimate
\begin{equation*}
    -\int_\Sigma \left\langle  D^k \Box\xi, \, D^k \tfrac{\nabla\xi}{\diff\tau} \right\rangle \dv_g 
    \leq \Vert D^k \Box \xi \Vert_{L^2} \left\Vert D^k\tfrac{\nabla\xi}{\diff\tau} \right\Vert_{L^2} 
    \leq \Vert D^k \Box \xi \Vert_{L^2} \, \mathcal{E}_{k+1}(\xi)^{\frac12},
\end{equation*}
which can similarly be applied to the terms involving commutators.
On the other hand we also have
\begin{equation*}
    \int_\Sigma \left(|D^{k+1}\xi|^2 - \left|D^k\tfrac{\nabla\xi}{\diff\tau}\right|^2 \right)H \dv_g
    \lesssim \Vert H \Vert_{C^0} \, \mathcal{E}_{k+1}(\xi)
\end{equation*}
and
\begin{align*}
    &\sum_{\ell=1}^k \binom{k}{\ell}  \int_\Sigma \left\langle D^\ell H \otimes D^{k-\ell}\tfrac{\nabla\xi}{\diff\tau}, \, D^k\tfrac{\nabla\xi}{\diff\tau} \right\rangle \dv_g\\
    &\lesssim \sum_{\ell=1}^k \left\Vert D^\ell H \otimes D^{k-\ell}\tfrac{\nabla\xi}{\diff\tau} \right\Vert_{L^2} \left\Vert D^k\frac{\nabla\xi}{\diff\tau} \right\Vert_{L^2} 
    \lesssim \Vert H \Vert_{C^k} \sum_{\ell=1}^{k} \mathcal{E}_{\ell}(\xi)^{\frac12} \mathcal{E}_{k+1}(\xi)^{\frac12}.
\end{align*}
We conclude that
\begin{align*}
    \frac{\diff}{\diff\tau} \mathcal{E}_{k+1}(\xi) 
    \lesssim&\;
    \bigg(
    \Vert D^k \Box \xi \Vert_{L^2} 
    + \left\Vert \left[\tfrac{\nabla}{\diff\tau}, D^k\right] \tfrac{\nabla\xi}{\diff\tau}\right\Vert_{L^2} 
    + \Vert [D^*D, D^k]\xi \Vert_{L^2} \\
    &+ \left\Vert \left[\tfrac{\nabla}{\diff\tau}, D^{k+1}\right]\xi \right\Vert_{L^2}
    + \Vert H \Vert_{C^k} \sum_{\ell=1}^{k} \mathcal{E}_{\ell}(\xi)^{\frac12}
    \bigg) \,
    \mathcal{E}_{k+1}(\xi)^{\frac12}
    + \Vert H \Vert_{C^0} \, \mathcal{E}_{k+1}(\xi)
\end{align*}
and thus
\begin{align*}
    \frac{\diff}{\diff\tau} &\mathfrak{E}_{k+1}(\xi) 
    \lesssim \mathfrak{E}_0(\xi)^{\frac12}\mathfrak{E}_1(\xi)^{\frac12} + \Vert H \Vert_{C^k} \, \mathfrak{E}_{k+1}(\xi)\\\nonumber
    &+\sum_{\ell=0}^k
    \bigg(
    \Vert D^\ell \Box \xi \Vert_{L^2} 
    + \left\Vert \left[\tfrac{\nabla}{\diff\tau}, D^\ell\right] \tfrac{\nabla\xi}{\diff\tau}\right\Vert_{L^2} 
    + \Vert [D^*D, D^\ell]\xi \Vert_{L^2} + \left\Vert \left[\tfrac{\nabla}{\diff\tau}, D^{\ell+1}\right]\xi \right\Vert_{L^2}
    \bigg) \,
    \mathfrak{E}_{\ell+1}(\xi)^{\frac12}.
\end{align*}

For the commutator terms, we note by (\ref{eq-commutator-curvature-estimate}) that
\begin{align*}
\left\Vert \left[\tfrac{\nabla}{\diff\tau}, D^k\right] \tfrac{\nabla\xi}{\diff\tau}\right\Vert_{L^2}&\lesssim 
\left\Vert D^{k-1} \left(R^V(\partial_\tau,\cdot)\tfrac{\nabla\xi}{\diff\tau}\right)\right\Vert_{L^2}
+\left\Vert\sff\right\Vert_{C^k}
\left\Vert \tfrac{\nabla\xi}{\diff\tau}\right\Vert_{H^{k}},\\
\left\Vert\left[\tfrac{\nabla}{\diff\tau}, D^{k+1}\right]\xi\right\Vert_{L^2}&\lesssim 
\left\Vert D^k(R^V(\partial_\tau,\cdot)\xi)\right\Vert_{L^2}
+\left\Vert\sff\right\Vert_{C^{k+1}}\left\Vert\xi\right\Vert_{H^{k+1}}.
\end{align*}
On the other hand, for $\xi \in \Gamma(T^\ast \Sigma^{\otimes k} \otimes V)$, $X \in \Gamma(T\Sigma)$ and $Y\in \Gamma(T\Sigma^{\otimes k})$ we also have
\begin{align*}
    &([D^*D, D]\xi)(X,Y)\\[0.1cm]
    &= 2\tr R^V(X,\cdot)((D\xi)(\cdot, Y)) - 2\tr(D\xi)(\cdot, \Riem_g(X,\cdot)Y) - (D\xi)(\Ric_g(X), Y)\\[0.1cm]
    &\quad + \tr(DR^V)(\cdot,X,\cdot)(\xi(Y)) + \xi(\tr(D\Riem_g)(\cdot,X,\cdot)Y),
\end{align*}
where the traces are taken with respect to $g_\tau$.
Using the identity (\ref{eq-commutator-power}) we can then write
    \begin{equation*}
    [D^*D, D^k]\xi = \sum_{\ell=0}^{k}D^{k-\ell}R^V\ast D^{\ell}\xi
    +\sum_{\ell=0}^{k}D^{k-\ell} \Riem_g\ast D^{\ell}\xi,
\end{equation*}
and we get 
\begin{equation*}
    \left\Vert [D^*D, D^k]\xi\right\Vert_{L^2}\lesssim 
    \left\Vert D^k(R^V\xi)\right\Vert_{L^2}
    +\left\Vert\Riem_g\right\Vert_{C^k}\left\Vert\xi\right\Vert_{H^k},
\end{equation*}
which together with the estimates above gives the result.
\end{proof}

Thus, to complete the energy estimate, one only needs to estimate the $H^k$-norm of $\Box^V \xi$, and the bundle curvature terms $\mathfrak{R}_k^V(\xi)$.
The norm of $\Box^V \xi$ will of course depend heavily on the equation that the section satisfies,
whereas the bundle curvature terms are fixed and directly computable once one specifies the bundle $V$.

The energy estimate (\ref{eq-energy-estimate-general}) can directly be applied to the Higgs energy $\mathfrak{E}_k(\Phi)$ and each of the terms separately from the Yang-Mills energy $\mathfrak{E}_k^+(F_\omega)=\mathfrak{E}_k(E_\omega)+\mathfrak{E}_k(B_\omega)$.
The Dirac energy $\mathfrak{E}_k^+(\Psi)$ is however defined with respect to an inner product that is no longer compatible with the spinor connection.
In particular, the evolution of $\mathfrak{E}_k^+(\Psi)$ picks up an additional term, since $e_0$ is not parallel with respect to spatial directions, in fact $\nabla_{e_i}e_0 = -\tensor{\sff}{_i^k}e_k$ (recall, however, that $\nabla_{e_0}e_0=0$).
More precisely, when commuting $D^*$ with Clifford multiplication by $e_0$, one gets a term of the form
\begin{equation*}
    \int_\Sigma \lLangle D_\omega^k \tfrac{\nabla_\omega\Psi}{\diff\tau}, \, e_0 \cdot \sff^{ij} e_i \cdot (e_j \iprod D_\omega^{k+1} \Psi) \rRangle \dv_g
    \lesssim 
    \Vert \sff \Vert_{C^0} \left( \left\Vert D_\omega^k \tfrac{\nabla_\omega\Psi}{\diff\tau} \right\Vert_{L^2}^2 + \Vert D_\omega^{k+1} \Psi \Vert_{L^2}^2 \right),
\end{equation*}
which already gets consumed in the right-hand side of the estimate (\ref{eq-energy-estimate-general}), and hence the estimate remains unchanged.

\subsubsection{Estimates involving products}

In estimating the norms, we will often use the following generalization of the Sobolev multiplication lemma.

\begin{Lem}\label{lem-sobolev}
    Let $(\Sigma^m, g)$ be a complete $m$-dimensional Riemannian manifold with bounded geometry, and $k$ a natural number satisfying $k>\frac{m}{2}$.
    For $N \geq 1$ and $1\leq i \leq N+1$, let $V_i\to M$ be Hermitian vector bundles, each equipped with a connection $D = D^{V_i}$.
    \begin{enumerate}
        \item[(i)]
            If $\xi_i \in H^k(V_i)$ and $\sum_{i=1}^{N} \ell_i\leq k$, then
            \begin{align*}
            \left\Vert D^{\ell_1}\xi_1 \ast \cdots \ast D^{\ell_N} \xi_N \right\Vert_{L^2} 
            \lesssim 
            \Vert \xi_1\Vert_{H^k} \cdots \Vert \xi_N\Vert_{H^k}.
            \end{align*}
        \item[(ii)]
        If also $\xi_{N+1}\in H^{k-1}(V_{N+1})$ and  $\sum_{i=1}^{N+1} \ell_i\leq k-1$, then 
        \begin{align*}
            \left\Vert D^{\ell_1}\xi_1 \ast \cdots \ast D^{\ell_N} \xi_N \ast D^{\ell_{N+1}} \xi_{N+1} \right\Vert_{L^2} 
            \lesssim  
            \Vert \xi_1\Vert_{H^k} \cdots \Vert \xi_N\Vert_{H^k} \Vert \xi_{N+1}\Vert_{H^{k-1}} .
        \end{align*}
    \end{enumerate}
    In both cases, $\lesssim$ denotes an inequality up to a constant depending on the Sobolev embedding on \((\Sigma,g)\),
    the numbers \(k\), and the contractions $\ast$.
\end{Lem}

\begin{proof}
See \cite[Lemma 3.6]{MR3830277}, \cite{MR1205820}, and \cite[\S VI.3, \S VI.16.3]{MR1016603}.
\end{proof}

\begin{Bem}\label{rem-star-notation}
    Here we use the $\ast$-notation as described at the end of \S \ref{sec-conventions}.
    We are also tacitly assuming that the bundles $V_i$ admit pairings such that the $\ast$-expressions on the left-hand side are well-defined and real-valued. For example, we will eventually wish to apply this result to $\Ad P$-valued forms for which the coefficients are contracted using the Lie commutator bracket, and also to contractions between $TM$- and $\Sigma M$-valued forms for which the coefficients are contracted using Clifford multiplication. The main point is that such contractions should be bounded in an appropriate sense, cf.\ sections that follow.
\end{Bem}

\subsubsection{Estimates involving Clifford multiplication}
\label{sec-estimates-clifford}

We will often have to estimate Clifford products between vectors and spinors, to which end we note that Clifford multiplication is only compatible with the ambient Levi-Civita connection $\nabla$ on $TM$, so that
\begin{equation*}
    D_\omega (X\cdot \Psi)(e_k) = \nabla_\omega (X\cdot \Psi)(e_k) = \nabla_{e_k} X \cdot \Psi + X\cdot(\nabla_\omega\Psi)_k = (D^{TM} X \cdot \Psi + X \cdot D_\omega\Psi)(e_k),
\end{equation*}
where $D^{TM}$ is the restriction of $\nabla$ to spatial tensors (cf.\ \S \ref{sec-spatial-tensors} and Footnote \ref{foot-LC-vs-spatial} on page \pageref{foot-LC-vs-spatial}).
Hence, we have the estimate
\begin{equation*}
    \Vert X \cdot \Psi \Vert_{H^k} \lesssim \Vert X \Vert_{C^k(TM)} \Vert \Psi \Vert_{H^k}
\end{equation*}
where
\begin{equation*}
    \Vert X \Vert_{C^k(TM)} = \sum_{\ell=0}^k \Vert (D^{TM})^{\ell} X \Vert_{C^0(TM)}
\end{equation*}
and the $C^0$-norms on the right-hand side are taken with respect to the Riemannian reference metric $\diff\tau \otimes \diff\tau + g_\tau$ on $TM$.
To estimate the $C^k(TM)$-norms, let
\begin{equation*}
    \Lambda = \tensor{\Lambda}{_k^i}  \, e^k \otimes e_i \in \Gamma(\pi^\ast T^\ast \Sigma \otimes \pi^\ast T\Sigma).
\end{equation*}
A simple calculation (cf.\ \S \ref{sec-adapted-frame}) then shows that
\begin{equation*}
    D^{TM} e_0 = -\tensor{\sff}{_k^i}  \, e^k \otimes e_i, \qquad D^{TM} \Lambda = - \tensor{\Lambda}{_k^m}\sff_{m\ell} \, e^k \otimes e^\ell \otimes e_0 + D\Lambda.
\end{equation*}
Now to estimate the $C^k(TM)$-norm of $e_0$ and $\Lambda$ one would have to calculate the higher derivatives $(D^{TM})^\ell$, which can in theory be done by "zig-zagging" between the formulae above, although there does not appear to be a simple closed formula for it.
However let us note that if we put
\begin{equation*}
    \mathcal{D}^p(\sff) = \sum_{q=0}^{p} \sum_{|a|=p-q} D^{a_1}\sff \ast \cdots \ast D^{a_q}\sff,
\end{equation*}
then schematically 
\begin{align*}
    (D^{TM})^\ell e_0 
    &=
    \sum_{p = 1}^{\ell-1} \mathcal{D}^p (\sff) \, e_0 + \sum_{q=0}^{\ell-1} \mathcal{D}^q(\sff) \ast \id^{T\Sigma} ,
    \\[0.1cm]
    (D^{TM})^\ell \Id^{T\Sigma}
    &=
    \sum_{p = 0}^{\ell-1} \mathcal{D}^p (\sff) \, e_0 + \sum_{q=1}^{\ell-1} \mathcal{D}^q(\sff) \ast \id^{T\Sigma} ,
\end{align*}
and we can trivially bound
\begin{equation*}
    \Vert \mathcal{D}^{p}(\sff) \Vert_{C^0} \lesssim \sum_{q=0}^{p} \Vert \sff \Vert_{C^q}^{p-q+1} \lesssim \Vert \sff \Vert_{C^p}(1 + \Vert \sff \Vert_{C^p}^p).
\end{equation*}
Therefore we have the rather weak but for our purposes sufficient estimates
\begin{equation*}
    \Vert e_0 \Vert_{C^k(TM)} 
    \lesssim  
    \Vert \sff \Vert_{C^{k-1}}(1+ \Vert \sff \Vert^{k-1}_{C^{k-1}}),
    \qquad
    \Vert \Lambda \Vert_{C^k(TM)} 
    \lesssim 
    \Vert \sff \Vert_{C^{k-1}}(1+ \Vert \sff \Vert^{k-1}_{C^{k-1}}) \Vert \Lambda \Vert_{C^k}.
\end{equation*}

\subsubsection{Estimates involving representations and Yukawa maps}

We will also often have to estimate terms which are acted upon by a representation of $\g$, in particular within the currents of (\ref{eq-yang-mills-wave}).
Let $\rho : G \to \mathbf{GL}(W)$ be a representation of $G$ and consider the associated vector bundle $E = P \times_\rho W$, equipped with the metric connection $\nabla_\omega$ induced by $\omega$.
If $\sigma \in \Gamma(E)$, then for $\rho_*\sigma \in \Gamma(E \otimes \Ad P)$ we have (cf.\ Lemma \ref{lemma-current-differential})
\begin{equation*}
   D_\omega (\rho_* \sigma) = \rho_* D_\omega \sigma,
\end{equation*}
so we can estimate
\begin{equation*}
    \Vert \rho_*\sigma \Vert_{H^k} \leq \Vert \rho_* \Vert   \Vert \sigma \Vert_{H^k},
\end{equation*}
where $$\Vert \rho_* \Vert = \max\{ \Vert \rho_*(X) \Vert \,\mid\, X \in \g, \, \Vert X \Vert = 1 \}$$ is the usual operator norm of the linear map $\rho_* : \g \to \End(W)$.
Analogously if $\xi \in \Gamma(\Ad P)$, then for $\rho_*(\xi) \in \Gamma(\End(E))$ we have
\begin{equation*}
    \Vert \rho_*(\xi) \Vert_{H^k} \leq \Vert \rho_* \Vert \Vert \xi \Vert_{H^k}.
\end{equation*}

For the Yukawa map, one should recall that we only have $\R$-linearity, but nevertheless its norm can be defined as
\begin{equation*}
    \Vert \Y \Vert = \sup \left\{ \frac{\Vert \Y_w \Vert}{\Vert w \Vert} \,\mid\, w \in W, \, w \not= 0 \right\} = \max\{ \Vert \Y_w \Vert \,\mid\, w \in W, \, \Vert w \Vert = 1 \},
\end{equation*}
where $\R$-linearity is sufficient to achieve the last equality, as it implies $\frac{1}{\Vert w \Vert}\Y_w = \Y_{\frac{w}{\Vert w \Vert}}$.
Thus we have
\begin{equation*}
    \Vert \Y_\Phi \Vert_{H^k} \leq \Vert \Y \Vert \Vert \Phi \Vert_{H^k}.
\end{equation*}
We note also that the complex antilinear part $\Y^-$ appearing in (\ref{eq-higgs}) satisfies
\begin{equation*}
    \Vert \Y^- \Vert \leq \Vert \Y \Vert. 
\end{equation*}

\subsection{Higgs energy estimate}
\label{sec-higgs-energy-estimates}

We start by considering the evolution of $\mathfrak{E}_k(\Phi) = \mathfrak{E}_k(\Phi; \omega)$, as this term is simpler to analyze than the other two, and
we use the general energy estimate (\ref{eq-energy-estimate-general}).
In particular, we will show below that for $k > 3/2$, we have
\begin{align}\label{eq-energy-estimate-higgs}
    \frac{\diff}{\diff\tau} \mathfrak{E}_{k+b}(\Phi) 
    \lesssim&\, \left( 1 + \Vert {\Riem_g} \Vert_{C^{k}} + \left\Vert \tfrac{\nabla\sff}{\diff\tau} \right\Vert_{C^k} + \Vert \sff \Vert_{C^{k+1}}^2 \right)\mathfrak{E}_{k+b}(\Phi)
    \\\nonumber
    &+ \mathfrak{E}^+_{k+b-1}(F_\omega)^{\frac12} \, \mathfrak{E}_{k+b}(\Phi) + \mathfrak{E}^+_{k}(\Psi) \, \mathfrak{E}_{k+b}(\Phi)^{\frac12} + |\lambda|\mathfrak{E}_{k+b}(\Phi)^2
\end{align}
which directly implies (\ref{eq-energy-estimate-higgs-total}).

\subsubsection{Contribution from the Euler-Lagrange equations}
We consider the system only with respect to the conformal Higgs potential, so that the Higgs equation reads
\begin{equation*}
    \Box_{\omega} \Phi - \tfrac16 \Scal_h \Phi - \lambda|\Phi|^2 \Phi - \langle \Psi, i\Y^-  \Psi \rangle = 0.
\end{equation*}
By Lemma \ref{lem-curvature-generalized-cyl} we can also decompose
\begin{equation}\label{eq-scalar-curvature}
    \Scal_h = -R(h)_{00} + g^{ij} R(h)_{ij} 
    = -2\tr\left(\tfrac{\nabla\sff}{\diff\tau}\right) + 6|\sff|^2 - 9H^2 + \Scal_g.
\end{equation}
We can estimate
\begin{align*}
    \Vert \Scal_h \Phi \Vert_{H^{k+b-1}} 
    &\lesssim
    \left( \left\Vert \tfrac{\nabla\sff}{\diff\tau} \right\Vert_{C^k} + \Vert \sff \Vert_{C^k}^2 + \Vert {\Scal_g} \Vert_{C^k}  \right) \Vert \Phi \Vert_{H^k},
    \\[0.2cm]
    \Vert \lambda|\Phi|^2 \Phi \Vert_{H^{k+b-1}}
    &\lesssim |\lambda| \sum_{|\ell|\leq k+b-1} \Vert D_\omega^{\ell_1} \Phi \ast D_\omega^{\ell_2} \Phi \ast D_\omega^{\ell_3} \Phi \Vert_{L^2} 
    \lesssim |\lambda|\Vert \Phi \Vert_{H^k}^3,
    \\
    \Vert  \langle \Psi, i\Y^-  \Psi \rangle \Vert_{H^{k+b-1}}
    &\lesssim
    \sum_{|\ell|\leq k+b-1} \Vert D_\omega^{\ell_1} \Psi \ast  \Y^- (D_\omega^{\ell_2} \Psi) \Vert_{L^2}
    \lesssim
    \Vert \Psi \Vert_{H^k}^2,
\end{align*}
where we also use Lemma \ref{lem-sobolev} in the last two estimates.
Combining all the terms, we get 
\begin{align}\label{eq-phi-wave-estimate}
    \Vert \Box_\omega \Phi \Vert_{H^{k+b-1}}
    &\lesssim
    \left( \left\Vert \frac{\nabla\sff}{\diff\tau} \right\Vert_{C^k} + \Vert \sff \Vert_{C^k}^2 + \Vert {\Scal_g} \Vert_{C^k}  \right) \mathfrak{E}_{k}(\Phi)^{\frac12} + |\lambda|\mathfrak{E}_{k}(\Phi)^{\frac32} + \mathfrak{E}^+_{k}(\Psi),
\end{align}
and thus also
\begin{equation*}
    \Vert \Box_\omega \Phi \Vert_{H^{k+b-1}} \mathfrak{E}_{k+b}(\Phi)^{\frac12} 
    \lesssim
    C(k,g,\sff) \, \mathfrak{E}_{k+b}(\Phi) + |\lambda|\mathfrak{E}_{k+b}(\Phi)^2 + \mathfrak{E}^+_{k}(\Psi) \, \mathfrak{E}_{k+b}(\Phi)^{\frac12}.
\end{equation*}

\subsubsection{Contribution from the bundle curvatures}
\label{sec-higgs-bundle-curvatures}

Next, we estimate the terms involving $\mathfrak{R}^\Hig_{k+b-1}(\Phi)$ from (\ref{eq-energy-estimate-general}), cf.\ also Proposition \ref{prop-general-energy-estimate}.
For the Higgs field bundle $\Hig$ we have $R^\Hig = \rho_*(F_\omega)$, and hence
\begin{equation*}
    R^\Hig(\partial_\tau, X) = \rho_*(E_\omega(X)), \quad R^\Hig(X,Y) = \rho_*(B_\omega(X,Y)).
\end{equation*}
We therefore get
\begin{align*}
    \left\Vert R^\Hig(\partial_\tau, \cdot) \frac{\nabla_\omega\Phi}{\diff\tau} \right\Vert_{H^{k+b-2}} 
    &\lesssim \Vert E_\omega \Vert_{H^{k+b-1}} \left\Vert \frac{\nabla_\omega\Phi}{\diff\tau}\right\Vert_{H^{k-1}},
    \\[0.2cm]
    \Vert R^\Hig\Phi \Vert_{H^{k+b-1}}
    &\lesssim \Vert B_\omega \Vert_{H^{k+b-1}} \Vert \Phi \Vert_{H^{k}},
    \\[0.2cm]
    \Vert R^\Hig(\partial_\tau, \cdot) \Phi \Vert_{H^{k+b-1}}
    &\lesssim \Vert E_\omega \Vert_{H^{k+b-1}} \Vert \Phi \Vert_{H^{k}},
\end{align*}
where in the first inequality we use Lemma \ref{lem-sobolev} (ii) and in the last two we use Lemma \ref{lem-sobolev} (i).
Hence,
\begin{equation*}
    \mathfrak{R}^\Hig_{k+b-1}(\Phi)\,\mathfrak{E}_{k+b}(\Phi)^{\frac12} 
    \lesssim \mathfrak{E}_{k+b-1}^+(F_\omega)^{\frac12} \, \mathfrak{E}_{k+b}(\Phi),
\end{equation*}
which completes the estimate (\ref{eq-energy-estimate-higgs}).

\subsection{Yang-Mills energy estimate}
\label{sec-yang-mills-energy-estimates}
Next, we treat the Yang-Mills energy $\mathfrak{E}^+_k(F_\omega) = \mathfrak{E}^+_k(F_\omega; \omega)$, which is technically more involved and contains some more constraining (the wave equation for $F$ is quadratic in first derivatives of $\Phi$) terms than the previous two.

\subsubsection{Contribution from the Euler-Lagrange equations}

We recall from (\ref{eq-yang-mills-wave}) that $F_\omega$ satisfies the wave equation
\begin{align*}
    \Box_\omega F_\omega 
    =&\, \langle [F_\omega \wedge F_\omega] \rangle + \mathcal RF_\omega + \Im \langle \id \cdot \Psi \wedge \chi_* \nabla_\omega \Psi \rangle\\[0.1cm]
    &+ \Re \langle \rho_*(F_\omega) \Phi \otimes \rho_*\Phi \rangle - \Re \langle \nabla_\omega \Phi \wedge  \rho_*\nabla_\omega\Phi \rangle.
\end{align*}
We need to split this equation to get expressions for $\Box_\omega E_\omega$ and $\Box_\omega B_\omega$, where $\Box_\omega$ needs to be taken with respect to the bundles $\pi^\ast T^\ast \Sigma \otimes \Ad P$ and $\pi^\ast \Lambda^2 \Sigma \otimes \Ad P$ respectively.

\begin{Prop}\label{prop-E-B-wave-eqs}
    Consider the decomposition $F_\omega = e^0 \wedge E_\omega + B_\omega$,
    where $E_\omega$ and $B_\omega$ are considered as spatial forms, i.e.\ sections of $\pi^\ast T^\ast \Sigma \otimes \Ad P$ and $\pi^\ast \Lambda^2 \Sigma \otimes \Ad P$ respectively.
    Then the electric part satisfies the wave equation
    \begin{align*}
        (\Box_\omega E_\omega)_i 
        =&\; 2[E^k, B_{ik}] + R(g)_{ik} E^k + \sff_{ik}\sff^{k\ell}E_\ell - 3H \sff_{ik}E^k + \left(\tfrac{\nabla\sff}{\diff\tau}\right)_{ik} E^k - \tr\left(\tfrac{\nabla\sff}{\diff\tau}\right) E_i \\[0.1cm]
        &  -2\tensor{\sff}{^k^\ell} (D_\omega B_\omega)_{k\ell i} + 2\tensor{(D^\ast \sff)}{^k} B_{ki} + 2 (D\sff)_{\ell ki}B^{k\ell} + 3(\diff H)_k \tensor{B}{^k_i} \\[0.1cm]
        &+ \Im \langle e_0 \cdot \Psi, \, \chi_*(D_\omega\Psi)_i \rangle - \Im\langle e_i \cdot \Psi, \, \chi_* \tfrac{\nabla_{\omega}\Psi}{\diff\tau} \rangle \\[0.1cm]
        &+ \Re \langle \rho_*(E_i)\Phi,\, \rho_*\Phi \rangle
        - 2\Re \langle \tfrac{\nabla_{\omega} \Phi}{\diff\tau},\, \rho_*(D_\omega \Phi)_i  \rangle,
    \end{align*}
    where $\Box_\omega = \Box_\omega^{\pi^\ast T^\ast \Sigma \otimes \Ad P}$,
    while the magnetic part satisfies
    \begin{align*}
        (\Box_\omega B_\omega)_{ij} 
        =&\; 2[\tensor{B}{^k_i}, B_{kj}] - 2[E_i, E_j]  - 2\tensor{\sff}{_i^k} (D_\omega E_\omega)_{kj} + 2\tensor{\sff}{_j^k} (D_\omega E_\omega)_{ki} \\[0.1cm]
        & + 2\left( (D\sff)_{ijk} - (D\sff)_{jik} \right) E^k - 3(\diff H)_i E_j + 3(\diff H)_j E_i \\[0.1cm]
        &+ R(g)_{ik}\tensor{B}{^k_j} + R(g)_{jk}\tensor{B}{_i^k} - 2R(g)_{kij\ell} B^{k\ell} - \left(\tfrac{\nabla\sff}{\diff\tau}\right)_{ik}\tensor{B}{^k_j} -  \left(\tfrac{\nabla\sff}{\diff\tau}\right)_{jk}\tensor{B}{_i^k} \\[0.1cm]
        &- 3H (\sff_{ik} \tensor{B}{^k_j} + \sff_{jk} \tensor{B}{_i^k}) -2 \sff_{kj} \sff_{i\ell} B^{k\ell} + \tensor{\sff}{_i^\ell} \sff_{\ell k}\tensor{B}{^k_j} + \tensor{\sff}{_j^\ell} \sff_{\ell k} \tensor{B}{_i^k}  \\[0.1cm]
        &+ \Im \langle e_i \cdot \Psi, \, \chi_*(D_\omega\Psi)_j \rangle - \Im\langle e_j \cdot \Psi, \, \chi_* (D_\omega\Psi)_j \rangle\\[0.1cm]
        &+ \Re \langle \rho_*(B_{ij})\Phi,\, \rho_*\Phi \rangle
        - 2 \Re \langle (D_\omega \Phi)_i,\, \rho_*(D_\omega \Phi)_j  \rangle,
    \end{align*}
    where $\Box_\omega = \Box_\omega^{\pi^\ast \Lambda^2 \Sigma \,\otimes \Ad P}$.
\end{Prop}

\begin{proof}
    We defer the proof to Appendix \ref{appendix-curvature-splitting}.
\end{proof}

\begin{Bem}
    For the purpose of doing estimates, it is convenient to gather terms of similar types and write these equations more schematically as
    \begin{align}\label{eq-E-wave-schematic}
        \Box_\omega E_\omega =&\, E_\omega \ast B_\omega + {\Riem_g} \ast E_\omega
        \\[0.1cm] \nonumber
        &+ \sff \ast \sff \ast E_\omega + \tfrac{\nabla\sff}{\diff\tau} \ast E_\omega + \sff \ast D_\omega B_\omega + D\sff \ast B_\omega
        \\[0.1cm] \nonumber
        &+ \Im \langle e_0 \cdot \Psi, \, \chi_* D_\omega\Psi \rangle - \Im\langle \id^{T\Sigma} \cdot \Psi, \, \chi_* \tfrac{\nabla_{\omega}\Psi}{\diff\tau} \rangle 
        \\[0.1cm] \nonumber
        &+ \Re \langle \rho_*(E_\omega)\Phi,\, \rho_*\Phi \rangle
        - 2\Re \langle \tfrac{\nabla_{\omega} \Phi}{\diff\tau},\, \rho_*D_\omega \Phi  \rangle,
    \end{align}
    and
    \begin{align}\label{eq-B-wave-schematic}
        \Box_\omega B_\omega =&\, B_\omega \ast B_\omega + E_\omega \ast E_\omega  + {\Riem_g} \ast B_\omega 
        \\[0.1cm] \nonumber
        &+ \sff \ast D_\omega E_\omega + D\sff \ast E_\omega + \tfrac{\nabla\sff}{\diff\tau} \ast B_\omega + \sff \ast \sff \ast B_\omega 
        \\[0.1cm] \nonumber
        &+ \Im \langle \id^{T\Sigma} \cdot \Psi \wedge \chi_*D_\omega\Psi \rangle + \Re \langle \rho_*(B_\omega)\Phi,\, \rho_*\Phi \rangle
        - \Re \langle D_\omega \Phi \wedge \rho_*D_\omega \Phi  \rangle,
    \end{align}
    where we use the $\ast$-notation as discussed in Remark \ref{rem-star-notation}.
\end{Bem}

Now we can estimate the norms $\Vert \Box_\omega E_\omega \Vert_{H^{k+a-1}}$ and $\Vert \Box_\omega B_\omega \Vert_{H^{k+a-1}}$ using (\ref{eq-E-wave-schematic}, \ref{eq-B-wave-schematic}).
We have that 
\begin{align*}
    \Vert E_\omega \ast B_\omega \Vert_{H^{k+a-1}}
    \lesssim \Vert E_\omega \Vert_{H^k} \Vert B_\omega \Vert_{H^k}.
\end{align*}
On the other hand for the linear terms in $E_\omega$ and $B_\omega$, we trivially have
\begin{align*}
    \left\Vert {\Riem_g} \ast E_\omega + \sff \ast \sff \ast E_\omega + \frac{\nabla\sff}{\diff\tau} \ast E_\omega \right\Vert_{H^{k+a-1}}
    &\lesssim \left(\Vert {\Riem_g} \Vert_{C^k} + \Vert \sff \Vert_{C^k}^2 + \left\Vert \frac{\nabla\sff}{\diff\tau} \right\Vert_{C^k} \right) \Vert E_\omega \Vert_{H^{k+a-1}}, \\[0.2cm]
    \Vert \sff \ast D_\omega B_\omega + D\sff \ast B_\omega \Vert_{H^{k+a-1}}
    &\lesssim \Vert \sff \Vert_{C^{k+a}}  \Vert B_\omega \Vert_{H^{k+a}}.
\end{align*}
The terms arising from the currents are a bit more subtle.
By the multinomial theorem, we can estimate the $\Phi$-current terms as
\begin{align*}
    \Vert \Re \langle \rho_*(E_\omega)\Phi, \rho_*\Phi \rangle \Vert_{H^{k+a-1}}
    &\lesssim 
    \sum_{|\ell| \leq k+a-1} \Vert \rho_*(D_\omega^{\ell_1} E_\omega) \ast D_\omega^{\ell_2}\Phi \ast \rho_* D_\omega^{\ell_3} \Phi\Vert_{L^2} 
    \lesssim
    \Vert E_\omega \Vert_{H^k} \Vert \Phi \Vert_{H^k}^2,
    \\
    \left\Vert \Re \left\langle \frac{\nabla_{\omega} \Phi}{\diff\tau},\, \rho_*D_\omega \Phi  \right\rangle \right\Vert_{H^{k+a-1}}
    &\lesssim 
    \sum_{|\ell| \leq k+a-1} \left\Vert\left(D_\omega^{\ell_1}\frac{\nabla_{\omega} \Phi}{\diff\tau}\right) \ast  (\rho_*D_\omega^{\ell_2}D_\omega \Phi)\right\Vert_{L^2} 
    \lesssim
    \left\Vert \frac{\nabla_\omega \Phi}{\diff\tau} \right\Vert_{H^k} \Vert \Phi \Vert_{H^{k+1}},
\end{align*}
and in the latter estimate we cannot get a Higgs energy order lower than $k+1$. 

For the $\Psi$-current term we also need to deal with Clifford multiplication, which we have discussed in \S \ref{sec-estimates-clifford}.
Thus
\begin{align*}
    \Vert \Im \langle e_0 \cdot \Psi, \, \chi_* D_\omega\Psi \rangle \Vert_{H^{k+a-1}}
    &\lesssim \sum_{|\ell|\leq k+a-1} \left\Vert (D^{TM})^{\ell_1} e_0 \ast D_\omega^{\ell_2} \Psi \ast \chi_* D_\omega^{\ell_3} D_\omega \Psi \right\Vert_{L^2}
    \\[0.1cm]
    &\lesssim \Vert e_0 \Vert_{C^k(TM)} \Vert \Psi \Vert_{H^{k+a}}^2,
    \\[0.2cm]
    \left\Vert \Im \left\langle \id^{T\Sigma} \cdot \Psi, \, \chi_* \frac{\nabla_{\omega}\Psi}{\diff\tau} \right\rangle \right\Vert_{H^{k+a-1}}
    &\lesssim
    \Vert {\Id^{T\Sigma}} \Vert_{C^k(TM)} \Vert \Psi \Vert_{H^{k+a}} \left\Vert \frac{\nabla_\omega\Psi}{\diff\tau} \right\Vert_{H^{k+a-1}}.
\end{align*}
Combining all these estimates, we get
\begin{align*}
    \Vert \Box_\omega E_\omega \Vert_{H^{k+a-1}} 
    &\lesssim
    \mathfrak{E}_{k}(E_\omega)^{\frac12} \mathfrak{E}_{k}(B_\omega)^{\frac12} 
    + \left(\Vert {\Riem_g} \Vert_{C^k} + \Vert \sff \Vert_{C^k}^2 + \left\Vert \frac{\nabla\sff}{\diff\tau} \right\Vert_{C^k} \right) \mathfrak{E}_{k}(E_\omega)^{\frac12}
    \\[0.1cm]\nonumber
    &\quad         
    + \Vert \sff \Vert_{C^{k+a}}  \mathfrak{E}_{k+a}(B_\omega)^{\frac12}
    + \mathfrak{E}_{k}(E_\omega)^{\frac12} \mathfrak{E}_{k}(\Phi)
    + \mathfrak{E}_{k+1}(\Phi)
    + \Vert \sff \Vert_{C^{k-1}}^k \mathfrak{E}_{k+a}(\Psi)
\end{align*}
and hence
\begin{align}\label{eq-E-total-box}
    \Vert \Box_\omega E_\omega \Vert_{H^{k+a-1}} \mathfrak{E}_{k+a}(E_\omega)^{\frac12}
    \lesssim&\,
    \mathfrak{E}_{k+a}^+(F_\omega)^{\frac32}
    + C(k,g,\sff) \, \mathfrak{E}_{k+a}^+(F_\omega)
    + \mathfrak{E}_{k+a}^+(F_\omega) \mathfrak{E}_{k}(\Phi)
    \\[0.1cm] \nonumber
    &
    + \mathfrak{E}_{k+a}^+(F_\omega)^{\frac12} \mathfrak{E}_{k+1}(\Phi)
    + C(k,\sff) \, \mathfrak{E}_{k+a}^+(F_\omega)^{\frac12}  \mathfrak{E}^+_{k+a}(\Psi).
\end{align}

Similarly, for the equation for $\Box_\omega B_\omega$ we find
\begin{equation*}
    \Vert B_\omega \ast B_\omega + E_\omega \ast E_\omega \Vert_{H^{k+a-1}}
    \lesssim \Vert B_\omega \Vert_{H^k}^2 + \Vert E_\omega \Vert_{H^k}^2,
\end{equation*}
while for the linear terms
\begin{align*}
    \left\Vert {\Riem_g} \ast B_\omega + \frac{\nabla\sff}{\diff\tau} \ast B_\omega + \sff \ast \sff \ast B_\omega \right\Vert_{H^{k+a-1}}
    &\lesssim \left(\Vert {\Riem_g} \Vert_{C^k} + \left\Vert \frac{\nabla\sff}{\diff\tau} \right\Vert_{C^k} + \Vert \sff \Vert_{C^k}^2\right) \Vert B_\omega \Vert_{H^{k}}, 
    \\[0.2cm]
    \Vert \sff \ast D_\omega E_\omega + D\sff \ast E_\omega \Vert_{H^{k+a-1}}
    &\lesssim  \Vert \sff \Vert_{C^{k+a}}  \Vert E_\omega \Vert_{H^{k+a}},
\end{align*}
and for the currents
\begin{align*}
    \Vert \Re \langle \rho_*(B_\omega)\Phi,\, \rho_*\Phi \rangle\Vert_{H^{k+a-1}}
    &\lesssim 
    \Vert B_\omega \Vert_{H^k} \Vert \Phi \Vert_{H^k}^2,
    \\[0.2cm]
    \Vert \Re \langle D_\omega \Phi \wedge \rho_*D_\omega \Phi  \rangle\Vert_{H^{k+a-1}}
    &\lesssim 
    \Vert \Phi \Vert_{H^{k+1}}^2,
    \\[0.2cm]
    \Vert \Im \langle \id^{T\Sigma} \cdot \Psi \wedge \chi_*D_\omega\Psi \rangle \Vert_{H^{k+a-1}}
    &\lesssim 
    \Vert \sff \Vert_{C^{k-1}}^k \Psi \Vert_{H^{k+a}}^2.
\end{align*}
Combining all the estimates, we get 
\begin{align}\label{eq-B-wave-estimate}
    \Vert \Box_\omega B_\omega \Vert_{H^{k+a-1}} 
    &\lesssim
    \mathfrak{E}_{k}(B_\omega) + \mathfrak{E}_{k}(E_\omega) 
    + \left(\Vert {\Riem_g} \Vert_{C^k} + \left\Vert \frac{\nabla\sff}{\diff\tau} \right\Vert_{C^k} + \Vert \sff \Vert_{C^k}^2\right) \mathfrak{E}_{k}(B_\omega)^{\frac12}
    \\[0.1cm]\nonumber
    &\quad 
    + \Vert \sff \Vert_{C^{k+a}} \mathfrak{E}_{k+a}(E_\omega)^{\frac12} 
    + \mathfrak{E}_{k}(B_\omega)^{\frac12} \mathfrak{E}_{k}(\Phi)
    + \mathfrak{E}_{k+1}(\Phi)
    + \Vert \sff \Vert_{C^{k-1}}^k \mathfrak{E}_{k+a}(\Psi),
\end{align}
and hence $B_\omega$ satisfies essentially the same estimate as $E_\omega$, cf.\ (\ref{eq-E-total-box}),
\begin{align}\label{eq-B-total-box}
    \Vert \Box_\omega B_\omega \Vert_{H^{k+a-1}} \mathfrak{E}_{k+a}(B_\omega)^{\frac12}
    \lesssim&\,
    \mathfrak{E}_{k+a}^+(F_\omega)^{\frac32}
    + C(k,g,\sff) \, \mathfrak{E}_{k+a}^+(F_\omega)
    + \mathfrak{E}_{k+a}^+(F_\omega) \mathfrak{E}_{k}(\Phi)
    \\[0.1cm] \nonumber
    &
    + \mathfrak{E}_{k+a}^+(F_\omega)^{\frac12} \mathfrak{E}_{k+1}(\Phi)
    + C(k,\sff) \, \mathfrak{E}_{k+a}^+(F_\omega)^{\frac12}  \mathfrak{E}^+_{k+a}(\Psi).
\end{align}

\subsubsection{Contribution from the bundle curvatures}
\label{sec-ym-bundle-curvatures}

For the bundle curvature terms we calculate, using Lemma \ref{lem-curvature-generalized-cyl}, the curvature tensor with respect to the connection $\nabla^{\pi^\ast T\Sigma}$ as
\begin{align*}
    h(R^{\pi^\ast T\Sigma}(\partial_\tau, e_i)e_j, e_k) &= \tensor{R}{_{0ijk}} = (D\sff)_{kji} - (D\sff)_{jik}, \\[0.1cm]
    h(R^{\pi^\ast T\Sigma}(e_i, e_j)e_k, e_\ell) &= \tensor{R(g)}{_{ijk\ell}}.
\end{align*}
Thus, for any $\xi \in \Gamma(\pi^\ast T^\ast \Sigma \otimes \Ad P)$, we find
\begin{align*}
    R^{\pi^\ast T^\ast \Sigma \otimes \Ad P}(\partial_\tau, e_k) \xi
    &= \left(\tensor{(D\sff)}{_j_k^i}-\tensor{(D\sff)}{^i_k_j}\right) e^j \otimes \xi_i + e^i \otimes [E_k, \xi_i]
    = D\sff \ast \xi + E_\omega \ast \xi,\\[0.2cm]
    R^{\pi^\ast T^\ast \Sigma \otimes \Ad P}(e_k, e_\ell) \xi
    &= \tensor{R(g)}{_k_\ell_j^i}\,  e^j \otimes \xi_i + e^i \otimes [B_{k\ell}, \xi_i]
    = {\Riem_g} \ast \xi + B_\omega \ast \xi.
\end{align*}
This implies that
\begin{align*}
    \left\Vert R^{\pi^\ast T^\ast \Sigma  \otimes \Ad P}(\partial_\tau, \cdot) \frac{\nabla_\omega E_\omega}{\diff\tau} \right\Vert_{H^{k+a-2}} 
    &\lesssim \Vert \sff \Vert_{C^k} \left\Vert \frac{\nabla_\omega E_\omega}{\diff\tau} \right\Vert_{H^{k-1}} + \Vert E_\omega \Vert_{H^{k}} \left\Vert \frac{\nabla_\omega E_\omega}{\diff\tau} \right\Vert_{H^{k-1}},
    \\[0.2cm]
    \Vert R^{\pi^\ast T^\ast \Sigma  \otimes \Ad P} E_\omega \Vert_{H^{k+a-1}}
    &\lesssim 
    \Vert {\Riem_g}\Vert_{C^k} \Vert E_\omega \Vert_{H^{k}} + \Vert B_\omega \Vert_{H^{k}} \Vert E_\omega \Vert_{H^k},
    \\[0.2cm]
    \Vert R^{\pi^\ast T^\ast \Sigma  \otimes \Ad P}(\partial_\tau, \cdot) E_\omega \Vert_{H^{k+a-1}}
    &\lesssim 
    \Vert \sff \Vert_{C^{k+1}} \Vert E_\omega \Vert_{H^k} + \Vert E_\omega \Vert_{H^k}^2,
\end{align*}
so that
\begin{equation}\label{eq-E-bundle-curvature}
    \mathfrak{R}_{k+a-1}^{\pi^\ast T\Sigma \otimes \Ad P} (E_\omega) \, \mathfrak{E}_{k+a}(E_\omega)^{\frac12}
    \lesssim
    C(k,g,\sff) \, \mathfrak{E}_{k+a}^+(F_\omega) + \mathfrak{E}_{k+a}^+(F_\omega)^{\frac32},
\end{equation}
and in an analogous way

\begin{equation}\label{eq-B-bundle-curvature}
    \mathfrak{R}_{k+a-1}^{\pi^\ast \Lambda^2\Sigma \otimes \Ad P} (B_\omega) \, \mathfrak{E}_{k+a}(B_\omega)^{\frac12}
    \lesssim
    C(k,g,\sff) \, \mathfrak{E}_{k+a}^+(F_\omega) + \mathfrak{E}_{k+a}^+(F_\omega)^{\frac32}.
\end{equation}

A combination of (\ref{eq-E-total-box}, \ref{eq-E-bundle-curvature}) and (\ref{eq-B-total-box}, \ref{eq-B-bundle-curvature}) gives (\ref{eq-energy-estimate-yang-mills-total}).

\subsection{Dirac energy estimate}
\label{sec-dirac-energy-estimates}

Finally, we treat the Dirac energy $\mathfrak{E}^+_k(\Psi) = \mathfrak{E}^+_k(\Psi; \omega)$.

\subsubsection{Contribution from the Euler-Lagrange equations}
We need to estimate $\Vert \Box_\omega \Psi \Vert_{H^{k+c-1}}$.
We recall from Theorem \ref{thm-wave-eqs} that $\Psi$ satisfies the wave equation
\begin{equation*}
    \Box_\omega \Psi =\tfrac14 \Scal \Psi + \chi_*(F_\omega) \cdot \Psi + \Y_{\nabla_\omega\Phi}\cdot\Psi - \Y_\Phi^2 \Psi.
\end{equation*}
Here, the term involving $\nabla_\omega\Phi$ requires special attention.
We can expand
\begin{equation*}
    \Y_{\nabla_\omega \Phi} \cdot \Psi = -e_0 \cdot \Y_{\frac{\nabla_\omega\Phi}{\diff\tau}}\Psi + \Y_{D_\omega\Phi} \cdot \Psi
\end{equation*}
and estimate
\begin{align*}
    &\Vert \Y_{\nabla_\omega \Phi} \cdot \Psi \Vert_{H^{k+c-1}} 
    \\[0.1cm]
    &= \left\Vert -e_0 \cdot \Y_{\frac{\nabla_\omega\Phi}{\diff\tau}}\Psi + \Y_{D_\omega\Phi} \cdot \Psi  \right\Vert_{H^{k+c-1}}
    \\[0.1cm]
    &\lesssim 
    \Vert e_0 \Vert_{C^{k+c-1}(TM)} \left\Vert \frac{\nabla_\omega\Phi}{\diff\tau} \right\Vert_{H^{k+c-1}} \Vert \Psi \Vert_{H^{k}}
    +
    \Vert \id^{T\Sigma} \Vert_{C^{k+c-1}(TM)} \left\Vert D_\omega\Phi \right\Vert_{H^{k+c-1}} \Vert \Psi \Vert_{H^{k}}
    \\[0.1cm]
    &\lesssim
    \Vert \sff \Vert_{C^{k-1}}(1+ \Vert \sff \Vert^{k-1}_{C^{k-1}}) \left( \left\Vert \frac{\nabla_\omega\Phi}{\diff\tau} \right\Vert_{H^{k+c-1}} + \Vert \Phi \Vert_{H^{k+c}} \right) \Vert\Psi\Vert_{H^k},
\end{align*}
cf.\ \S \ref{sec-estimates-clifford}.
It follows that
\begin{equation*}
    \Vert \Y_{\nabla_\omega \Phi} \cdot \Psi \Vert_{H^{k+c-1}} \lesssim
    \Vert \sff \Vert_{C^{k-1}}^{k} \left( \left\Vert \frac{\nabla_\omega\Phi}{\diff\tau} \right\Vert_{H^{k+c-1}} + \Vert \Phi \Vert_{H^{k+c}} \right) \Vert\Psi\Vert_{H^k}.
\end{equation*}

The other terms can be estimated as follows.
Using $(\ref{eq-scalar-curvature})$, we get
\begin{equation*}
    \Vert {\Scal_h} \Psi \Vert_{H^{k+c-1}}
    \lesssim
    \left( \left\Vert \tfrac{\nabla\sff}{\diff\tau} \right\Vert_{C^{k}} + \Vert \sff \Vert_{C^{k}}^2 + \Vert {\Scal_g} \Vert_{C^{k}} \right) \Vert \Psi \Vert_{H^{k+c-1}}.
\end{equation*}
On the other hand, by (\ref{eq-dirac-yukawa-commutator}), Lemma \ref{lem-sobolev} and the fact that $\Y$ is real linear in $\Phi$, we have
\begin{align*}
    \Vert \Y_\Phi^2\Psi \Vert_{H^{k+c-1}}
    &\lesssim
    \sum_{|\ell|\leq k+c-1} \Vert \Y_{D_\omega^{\ell_1}\Phi} \Y_{D_\omega^{\ell_2}\Phi} D_\omega^{\ell_3}\Psi \Vert_{L^2} \\
    &\lesssim \sum_{|\ell| \leq k+c-1} \Vert D_\omega^{\ell_1}\Phi \ast D_\omega^{\ell_2}\Phi \ast D_\omega^{\ell_3}\Psi \Vert_{L^2}  
    \lesssim \Vert \Phi \Vert_{H^k}^2 \Vert \Psi \Vert_{H^{k+c-1}}.
\end{align*}
Finally, from Lemma \ref{lem-curvature-splitting}, we find
\begin{align*}
    D_\omega (\chi_*(F_\omega) \cdot \Psi) &= \chi_*(D_\omega^{TM} F_\omega) \cdot \Psi + \chi_*(F_\omega) \cdot D_\omega \Psi \\
    &= - \, e^k \otimes e_0 \cdot e_i \cdot \left( E^i (D_\omega\Psi)_k + \left(\tensor{(D_\omega E_\omega)}{_k^i} + \tensor{\sff}{_k_\ell}B^{\ell i}\right)\Psi \right)\\
    &\; +\tfrac12 \, e^k \otimes e_i \cdot e_j \cdot \left( B^{ij} (D_\omega\Psi)_k + \left( \tensor{(D_\omega B_\omega)}{_k^i^j} - \tensor{\sff}{_k^j}E^i + \tensor{\sff}{_k^i}E^j \right)\Psi \right).
\end{align*}

Iterating, we get
\begin{equation*}
    \Vert \chi_*(F_\omega) \cdot \Psi \Vert_{H^{k+c-1}} 
    \lesssim
    \Vert \sff \Vert_{C^{k-1}}^2(1+ \Vert \sff \Vert^{k-1}_{C^{k-1}})^2 \left( \Vert E_\omega \Vert_{H^{k+c-1}} + \Vert B_\omega \Vert_{H^{k+c-1}} \right) \Vert \Psi \Vert_{H^{k}},
\end{equation*}
cf.\ \S \ref{sec-estimates-clifford}.

Combining all the estimates, we get that 
\begin{align}\label{eq-psi-wave-estimate}
    \Vert \Box_\omega\Psi \Vert_{H^{k+c-1}} 
    &\lesssim
    \left( \left\Vert \frac{\nabla\sff}{\diff\tau} \right\Vert_{C^{k}} + \Vert \sff \Vert_{C^{k}}^2 + \Vert {\Scal_g} \Vert_{C^{k}} \right) \mathfrak{E}^+_{k+c-1}(\Psi)^{\frac12} + \mathfrak{E}_{k}(\Phi) \mathfrak{E}^+_{k+c-1}(\Psi)^{\frac12} 
    \\\nonumber
    &\quad + \Vert \sff \Vert_{C^{k-1}}^{2k} \left( \mathfrak{E}_{k+c}(\Phi)^{\frac12} + \mathfrak{E}^+_{k+c-1}(F_\omega)^{\frac12} \right) \mathfrak{E}^+_{k}(\Psi)^{\frac12}.
\end{align}

\subsubsection{Contribution from the bundle curvatures}
\label{sec-dirac-bundle-curvatures}
We first note that
\begin{align*}
    R^{\Fer}(X,Y)\Psi &= R^{\Sigma M}(X,Y)\Psi + R^{\mathscr{S}}(X,Y)\Psi\\[0.1cm]
    &=  -\frac14 \Riem(X,Y) e_\mu \cdot e^\mu \cdot \Psi + \chi_*(F_\omega(X,Y))\Psi,
\end{align*}
where in the last equality we use the well-known relation (\ref{eq-riemann-vs-spinor}) between the curvature tensor of the spinor bundle and the Riemann tensor.
The $R^{\mathscr{S}}$ terms can be estimated as (cf.\ \S \ref{sec-higgs-bundle-curvatures})
\begin{align*}
    \left\Vert R^{\mathscr{S}}(\partial_\tau, \cdot) \frac{\nabla_\omega\Psi}{\diff\tau} \right\Vert_{H^{k+c-2}} 
    &\lesssim \Vert E_\omega \Vert_{H^{k+c-1}} \left\Vert \frac{\nabla_\omega\Psi}{\diff\tau}\right\Vert_{H^{k-1}},
    \\[0.2cm]
    \Vert R^{\mathscr{S}}\Phi \Vert_{H^{k+c-1}}
    &\lesssim \Vert B_\omega \Vert_{H^{k+c-1}} \Vert \Psi \Vert_{H^k},
    \\[0.2cm]
    \Vert R^{\mathscr{S}}(\partial_\tau, \cdot) \Psi \Vert_{H^{k+c-1}}
    &\lesssim \Vert E_\omega \Vert_{H^{k+c-1}} \Vert \Psi \Vert_{H^{k}}.
\end{align*}
The spinorial curvature terms are more involved, however. Using Lemma \ref{lem-curvature-generalized-cyl}, we find
\begin{align*}
    R^{\Sigma M}(\partial_\tau,e_i)\Psi &= \frac14 \Riem(e_0, e_i) e_0 \cdot e_0 \cdot \Psi - \frac14 \Riem(e_0,e_i) e_k \cdot e^k \cdot \Psi\\
    &= \frac12 \left( \tensor{\left(\tfrac{\nabla\sff}{\diff\tau}\right)}{_i^k} - \sff_{i\ell} \sff^{\ell k} \right) e_0 \cdot e_k \cdot \Psi
    - \frac14 \left( \tensor{(D\sff)}{^\ell^k_i} - \tensor{(D\sff)}{^k^\ell_i} \right) e_\ell \cdot e_k \cdot \Psi,
    \\[0.2cm]
    R^{\Sigma M}(e_i, e_j) \Psi &= \frac14 \Riem(e_i, e_j) e_0 \cdot e_0 \cdot \Psi - \frac14 \Riem(e_i,e_j) e_k \cdot e^k \cdot \Psi\\
    &= \frac12 \left( \tensor{(D\sff)}{_i_j^k} - \tensor{(D\sff)}{_j_i^k} \right) e_0 \cdot e_k \cdot \Psi
    - \frac14 \left( \tensor{R(g)}{_i_j^k^\ell} + \tensor{\sff}{_i^k} \tensor{\sff}{_j^\ell} - \tensor{\sff}{_j^k} \tensor{\sff}{_i^\ell} \right) e_\ell \cdot e_k \cdot \Psi.
\end{align*}
Then as before
\begin{align*}
    \left\Vert R^{\Sigma M}(\partial_\tau, \cdot)\frac{\nabla_\omega\Psi}{\diff\tau} \right\Vert_{H^{k+c-2}}
    &\lesssim 
    \left( \left\Vert \tfrac{\nabla\sff}{\diff\tau} \right\Vert_{C^{k-1}} + \Vert \sff \Vert_{C^{k-1}}^2 + \Vert \sff \Vert_{C^{k}} \right) \Vert \sff \Vert_{C^{k-2}}^2(1+ \Vert \sff \Vert^{k-2}_{C^{k-2}})^2 \Vert \Psi \Vert_{H^{k+c-1}},
    \\[0.2cm]
    \Vert R^{\Sigma M}(\partial_\tau, \cdot)\Psi \Vert_{H^{k+c-1}}
    &\lesssim 
    \left( \left\Vert \tfrac{\nabla\sff}{\diff\tau} \right\Vert_{C^k} + \Vert \sff \Vert_{C^k}^2 + \Vert \sff \Vert_{C^{k+1}} \right) \Vert \sff \Vert_{C^{k-1}}^2(1+ \Vert \sff \Vert^{k-1}_{C^{k-1}})^2 \Vert \Psi \Vert_{H^{k+c-1}},
    \\[0.2cm]
    \Vert R^{\Sigma M}(\cdot, \cdot) \Psi \Vert_{H^{k+c-1}}
    &\lesssim
    \left( \Vert \sff \Vert_{C^{k+1}} + \Vert \sff \Vert_{C^k}^2 + \Vert {\Riem_g} \Vert_{C^k} \right) \Vert \sff \Vert_{C^{k-1}}^2(1+ \Vert \sff \Vert^{k-1}_{C^{k-1}})^2 \Vert \Psi \Vert_{H^{k+c-1}},
\end{align*}
and it follows that
\begin{equation*}
    \mathfrak{R}_{k+c-1}^{\Sigma M \otimes \Ad P} (\Psi) \, \mathfrak{E}^+_{k+c}(\Psi)^{\frac12}
    \lesssim 
    C(k,\sff,g) \,\mathfrak{E}^+_{k+c}(\Psi) + \mathfrak{E}^+_{k+c-1}(F_\omega)^{\frac12} \mathfrak{E}^+_{k+c}(\Psi).
\end{equation*}

\section{Proof of the main theorem}
\label{sec-main-result}

The goal of this section is to prove Theorem \ref{thm-main}.
To do this, we will construct a conformal Standard Model on the Gaussian foliation $[0,T) \times \Sigma$, and then conformally transform it back to the original spacetime, where it will be future global, cf.\ \S \ref{sec-global-original}.

\subsection{Local existence}

In view of the definition of the space $C^0(I, \mathcal{H}^{k+1}(P^\Sigma))$ of Sobolev connections, it is natural to choose a connection $\omega^0 = \Pi^\ast \sigma^0$ for some connection $\sigma^0$ on $P^\Sigma$ of bounded geometry of order $k+1$, and perturb it by a spatial $\Ad P$-valued form $\eta$, having $\eta^0$ as its initial value, i.e.\ we consider connections of the form
\begin{equation*}
    \omega = \omega^0 + \eta,
\end{equation*}
and view $\eta$ as the unknown.
This approach has two main benefits.
Firstly, having a fixed reference connection will allow us to write down the Yang-Mills equation as an equation for $\eta$, and this equation will be tensorial.
Secondly, the connection $\omega^0$ has bounded geometry of order $k+1$ uniformly for $\tau \in I$,  which will also ensure that the necessary Sobolev machinery holds for the the corresponding Sobolev spaces defined in terms of $\omega^0$.

Observe that by (\ref{eq-curvature-sum-intro}), the curvature form satisfies
\begin{equation}\label{eq-curvature-sum}
    F_\omega = F_{\omega^0} + \diff_{\omega^0} \eta + \frac12 [\eta \wedge \eta],
\end{equation}
where we have $F_{\omega^0} = \pi^\ast F_{\sigma^0}$ since $\omega^0 = \Pi^\ast \sigma^0$.
In particular, we also see that $\partial_\tau \iprod F_{\omega^0} = 0$.

The idea now is to view $F_\omega$ as an independent quantity $F \in \Omega^2(M,\Ad P)$, and then use (\ref{eq-curvature-sum}) to set up the equations for $\eta$, resp.\ the Bianchi identity and the Yang-Mills equation to set up the equations for $F$. 
It will then eventually follow that $F$ are indeed the components of the curvature form $F_\omega$.
This idea (used in a local temporal gauge) is due to \cite{MR1604914, MR649158, MR638511}.
This strategy seems to be well-suited for our global gauge invariant approach, although we would like to point out that the local existence result can also be obtained by working in a local Lorenz gauge, along the lines of \cite[\S 4]{MR654209}.
More precisely, the idea is to solve the PDE system%
\footnote{
We wish to separate all the derivatives in the normal (i.e.\ temporal) direction from the ones tangential to $\Sigma$. 
For $\diff_{\omega^0} \sigma$ and $\diff_\omega F$, the normal (resp.\ tangential) derivatives all reside in terms which contain (resp.\ do not contain) $\diff\tau$, while dually for $\diff_\omega^* F$ these are all the terms that do not contain (resp.\ contain) $\diff\tau$.
} 
\begin{subequations}
\begin{empheq}[left=\empheqlbrace]{align}
    \label{eq-local-sys-connection}
    \partial_\tau \iprod \left(\diff_{\omega^0} \eta - F + F_{\omega^0} + \frac12 [\eta \wedge \eta]\right) = 0,\\[0.1cm]
    \label{eq-local-sys-bianchi}
    \partial_\tau \iprod \diff_\omega F = 0,\\[0.1cm]
    \label{eq-local-sys-ym}
    \diff\tau \wedge (\diff_\omega^* F - \mathfrak J[\omega,\Phi,\Psi]) = 0,\\[0.1cm]
    \label{eq-local-sys-higgs}
    \Box_{\omega} \Phi - \tfrac16 \Scal_h \Phi - \lambda|\Phi|^2 \Phi - \langle \Psi, i\Y^-  \Psi \rangle = 0, \\[0.1cm]
    \label{eq-local-sys-dirac}
    \Box_\omega \Psi - \tfrac14 \Scal_h \Psi - \chi_*(F) \cdot \Psi - \Y_{\nabla_\omega\Phi} \cdot \Psi + \Y_\Phi^2 \Psi = 0,
\end{empheq}
\end{subequations}
subject to the constraints
\begin{subequations}
\begin{empheq}[left=\empheqlbrace]{align}
    \label{eq-local-sys-connection-constraint}
    \diff\tau \wedge \left(\diff_{\omega^0} \eta - F + F_{\omega^0} + \frac12 [\eta \wedge \eta]\right) = 0, \\[0.1cm]
    \label{eq-local-sys-bianchi-constraint}
    \diff\tau \wedge \diff_\omega F = 0,\\[0.1cm]
    \label{eq-local-sys-ym-constraint}
    \partial_\tau \iprod (\diff_\omega^* F - \mathfrak J[\omega,\Phi,\Psi]) = 0,\\[0.1cm]
    \label{eq-local-sys-dirac-constraint}
    \dirac_{\omega} \Psi + \Y_\Phi \Psi = 0,
\end{empheq}
\end{subequations}
where we recall that
\begin{equation*}
    \mathfrak J[\omega,\Phi,\Psi] = -\Re \langle \nabla_\omega\Phi \otimes \rho_*\Phi \rangle + \tfrac12 \Im \langle \id\cdot\Psi \otimes \chi_*\Psi \rangle.
\end{equation*}
In the following sections, we will see that this is a semilinear symmetric hyperbolic system, so that it can be treated using standard methods at a local level, and then the energy estimate derived in \S \ref{sec-energy-estimates} will give the desired global result.

\subsubsection{Evolution equations}

Let us write the system (\ref{eq-local-sys-connection}--\ref{eq-local-sys-dirac}) with respect to an adapted frame%
\footnote{We do the computations in an adapted frame for simplicity, but the formulae are tensorial so they also hold in any other frame.}
but in a (principal bundle) gauge invariant manner.
Observe that
\begin{equation*}
    \frac{\nabla_\omega \xi}{\diff\tau} = \frac{\nabla_{\omega^0} \xi}{\diff\tau},\qquad 
    D_\omega \xi = D_{\omega^0 + \eta} \, \xi = D_{\omega^0}\, \xi + [\eta \otimes \xi],
\end{equation*}
for any spatial $\Ad P$-valued form $\xi$, since $\eta$ is also purely spatial.
As before, we will decompose $F = \diff\tau \wedge E + B$.
For convenience, we also define the (spatial) Hodge dual of $B$ via
\begin{equation*}
    Q = \star B = e^i \otimes \frac12 \tensor{\epsilon}{_i^{jk}} B_{jk},
\end{equation*}
where $\star = \star_{g_t}$ here and throughout the rest of the section denotes the Hodge star of the spatial slice $(\Sigma, g_\tau)$, and $\epsilon$ is the Levi-Civita symbol.

A calculation using \S \ref{sec-adapted-frame} gives
\begin{equation*}
    (\diff_{\omega^0}\eta)_{0i} = (\nabla_{\omega^0}\eta)_{0i} - (\nabla_{\omega^0}\eta)_{i0} = \left(\frac{\nabla_{\omega^0} \eta}{\diff\tau}\right)_i - \tensor{\sff}{_i^k} \eta_k,
\end{equation*}
and thus if we evaluate (\ref{eq-curvature-sum}) at $0i$, which is equivalent to (\ref{eq-local-sys-connection}), we see that $\eta$ satisfies the evolution equation
\begin{equation}\label{eq-eta-evolution}
    \left(\frac{\nabla_{\omega^0} \eta}{\diff\tau}\right)_i = \tensor{\sff}{_i^k} \eta_k + E_i.
\end{equation}
In a similar way, using Lemma \ref{lem-curvature-splitting} and the formulae
\begin{equation}\label{eq-diff-codiff-components}
    (\diff_\omega F_\omega)_{\mu\nu\lambda} = (\nabla_\omega F_\omega)_{\mu\nu\lambda}
    + (\nabla_\omega F_\omega)_{\nu\lambda\mu}
    + (\nabla_\omega F_\omega)_{\lambda\mu\nu}, \qquad
    (\diff_\omega^\ast F_\omega)_\nu = - \tensor{(\nabla_\omega F_\omega)}{^\mu_\mu_\nu},
\end{equation}
we see that the Bianchi and Yang-Mills equations (\ref{eq-local-sys-bianchi}, \ref{eq-local-sys-ym}) are equivalent to 
\begin{subequations}
    \begin{empheq}[left=\empheqlbrace]{align}
    \label{eq-Q-evolution}
    \left( \frac{\nabla_{\omega^0} Q}{\diff\tau} \right)_i &= \tensor{\epsilon}{_i^{jk}} \left( (D_{\omega^0}E)_{jk} + [\eta_j, E_k] \right) + 3 H Q_i - \tensor{\sff}{_i^k} Q_i, \\[0.1cm]
    \label{eq-E-evolution-intermsofQ}
    \left(\frac{\nabla_{\omega^0} E}{\diff\tau}\right)_i &= \tensor{(D_{\omega^0} \star Q)}{^k_{ki}} + \tensor{\epsilon}{^j_{ki}}\, [\eta^k, Q_j] + 3H E_i - \tensor{\sff}{_i^k} E_k + \mathfrak J_i,  
    \end{empheq}
\end{subequations}
where we recall that $H = \tfrac13\tr_g(\sff)$ is the mean curvature and 
\begin{equation*}
    \mathfrak J_i = -\Re \langle (D_\omega\Phi)_i \otimes \rho_*\Phi \rangle + \tfrac12 \Im \langle e_i \cdot\Psi \otimes \chi_*\Psi \rangle.
\end{equation*}

The equations (\ref{eq-local-sys-higgs}, \ref{eq-local-sys-dirac}) are wave equations, for which it is well-known that they can be rewritten as a symmetric hyperbolic system.
Indeed, if we put
\begin{equation*}
    \dot\Phi := \frac{\nabla_{\omega^0} \Phi}{\diff\tau}, \qquad Z := D_\omega\Phi, \qquad \dot\Psi := \frac{\nabla_{\omega^0} \Psi}{\diff\tau}, \qquad S := D_\omega\Psi
\end{equation*}
then using (\ref{eq-commutator-temporal-spatial-explicit}) to commute the temporal and spatial derivatives, we find that the Higgs equation (\ref{eq-local-sys-higgs}) is equivalent to 
\begin{subequations}
\begin{empheq}[left=\empheqlbrace]{align}
    \label{eq-Phi-evolution}
    \frac{\nabla_{\omega^0} \Phi}{\diff\tau} &= \dot\Phi\\[0.1cm]
    \frac{\nabla_{\omega^0} \dot\Phi}{\diff\tau} &= - D_{\omega^0}^\ast Z + \rho_*(\eta^k) Z_k + 3H \dot\Phi - \tfrac16 \Scal_h \Phi - \lambda|\Phi|^2 \Phi - \langle \Psi, i\Y^-  \Psi \rangle\\[0.1cm]
    \label{eq-R-evolution}
    \left(\frac{\nabla_{\omega^0} Z}{\diff\tau}\right)_i &= (D_{\omega^0}\dot\Phi)_i + \rho_*(\eta_i)\dot\Phi + \rho_*(E_i) \Phi +
    \tensor{\sff}{_i^k} Z_k,
\end{empheq}
\end{subequations}
and the Dirac equation (\ref{eq-local-sys-dirac}) is equivalent to
\begin{subequations}
\begin{empheq}[left=\empheqlbrace]{align}
    \label{eq-Psi-evolution}
    \frac{\nabla_{\omega^0} \Psi}{\diff\tau} &= \dot\Psi\\[0.1cm]
    \frac{\nabla_{\omega^0} \dot\Psi}{\diff\tau} &= - D_{\omega^0}^\ast S + \chi_*(\eta^k)S_k + 3H \dot\Psi  - \tfrac14 \Scal_h \Psi  \\\nonumber
    & \quad  + e_0 \cdot \chi_*(E) \cdot \Psi - \chi_*(\star \, Q) \cdot \Psi + \partial_\tau \cdot \Y_{\dot\Phi} \Psi - \Y_{Z} \cdot \Psi + \Y_\Phi^2 \Psi \\[0.1cm]
    \label{eq-S-evolution}
    \left(\frac{\nabla_{\omega^0} S}{\diff\tau}\right)_i &= (D_{\omega^0}\dot\Psi)_i + \chi_*(\eta_i)\dot\Psi  + \frac12 \left( \tensor{\left(\frac{\nabla\sff}{\diff\tau}\right)}{_i^k} - \sff_{i\ell} \sff^{\ell k} \right) \partial_\tau \cdot e_k \cdot \Psi  \\\nonumber
    & \quad 
    - \frac14 \left( \tensor{(D\sff)}{^\ell^k_i} - \tensor{(D\sff)}{^k^\ell_i} \right) e_\ell \cdot e_k \cdot \Psi + \chi_*(E_i)\Psi + \tensor{\sff}{_i^k}S_k,
\end{empheq}
\end{subequations}
where in the last equation we also use the formulas from \S \ref{sec-dirac-bundle-curvatures}.
We also note that the current in the Yang-Mills equation (\ref{eq-E-evolution-intermsofQ}) is given by
\begin{equation*}
    \mathfrak J_i = -\Re \langle Z_i, \rho_*\Phi \rangle + \tfrac12 \Im \langle e_i \cdot\Psi, \chi_*\Psi \rangle.
\end{equation*}

Now let
\begin{equation*}
    \mathscr{U} = (\pi^\ast T^\ast \Sigma \otimes \Ad P)^3 \oplus \Hig^2 \oplus (\pi^\ast T^\ast \Sigma \otimes \Hig) \oplus \Fer_+^2 \oplus (\pi^\ast T^\ast \Sigma \otimes \Fer_+)
\end{equation*}
be the total bundle with sections $u = (\eta,Q,E,\Phi,\dot\Phi,Z,\Psi,\dot\Psi,S) \in \Gamma(\mathscr{U})$.
Define 
\begin{equation*}
    X : \Gamma(\mathscr{U}) \to \Gamma(\mathscr{U})
\end{equation*}
as the (non-linear) first-order differential operator associated to the equations (\ref{eq-eta-evolution}, \ref{eq-E-evolution-intermsofQ}, \ref{eq-Q-evolution}, \ref{eq-Phi-evolution}--\ref{eq-R-evolution}, \ref{eq-Psi-evolution}--\ref{eq-S-evolution}).
A quick inspection of the equations shows that we can decompose the operator as $X=L+N$, where $L$ is the linear first order operator containing all the derivative terms, and the remainder $N$ is a zeroth order nonlinear operator.
Let $\sigma_{L}$ be the principal symbol of the linear part, which we recall is the map taking one-forms $\theta \in \Gamma(T^\ast M)$ to maps of sections $\sigma_{L}(\theta) : \Gamma(\mathscr{U}) \to \Gamma(\mathscr{U})$, characterized by
$L(f u) = f L(u) + \sigma_{L}(\diff f)(u)$,
where $u \in \Gamma(\mathscr{U})$ and $f \in C^\infty(M)$, cf.\ \cite[\S 2]{MR3846239} and \cite[Definition 5.1]{MR3302643}.
It is easy to see that $\sigma_{L}$ is in this setting given by the vector bundle endomorphism
\begin{equation*}
    \sigma_{L}(\diff\tau) = \id^{\mathscr{X}}, 
    \qquad 
    \sigma_{L}(\xi) = -\left[
    \begin{array}{ccc|ccc|ccc}
        0 & & & & & & & &\\
        & 0 & -\xi^\sharp \iprod \star  & & & & & & \\
        & \xi^\sharp \iprod \star & 0 & & & & & & \\ \hline
        & & & 0 & & & & &\\
        & & & & 0 & \xi^\sharp \iprod & & & \\
        & & & & \xi \,\wedge & 0 & & & \\ \hline
        & & & & & & 0 & & \\
        & & & & & & & 0 & \xi^\sharp \iprod \\
        & & & & & & & \xi \,\wedge & 0
   \end{array}\right],
\end{equation*}
where $\xi \in \Gamma(\pi^\ast T^\ast \Sigma)$ is spatial and $\star = \star_{g_\tau}$ denotes the Hodge star on the spatial slice $(\Sigma, g_\tau)$.%
\footnote{
If we interpret one-forms as vectors (i.e.\ by viewing the frame elements $e^j$ as a basis), we can also write the spatial part of $\sigma_L$ in matrix form as
\begin{equation*}
    \sigma_{L}(e^j) = -
    \left[
    \begin{array}{ccc|ccc|ccc}
        0 & & & & & & & &\\
        & 0 & \epsilon^j & & & & & & \\
        & -\epsilon^j & 0 & & & & & & \\ \hline
        & & & 0 & & & & &\\
        & & & & 0 & (\delta^j)^\top & & & \\
        & & & & \delta^j & 0 & & & \\ \hline
        & & & & & & 0 & & \\
        & & & & & & & 0 & (\delta^j)^\top\\
        & & & & & & & \delta^j & 0
   \end{array}\right],
\end{equation*}
where $\delta^j$ is the vector with $(\delta^j)_k = \tensor{\delta}{^j_k}$ (i.e.\ has 1 in the $j$-th position and zeros elsewhere), and $\epsilon^j$ are the skew-symmetric matrices with components $\tensor{(\epsilon^j)}{_i^k} = \tensor{\epsilon}{_i^{jk}}$.
}
In particular, we see that $\sigma_{L}(\diff\tau)$ is positive-definite and $\sigma_{L}(\theta)$ is symmetric (with respect to the inner product on $\mathscr{U}$ induced by the natural inner products on the factors) for each $\theta \in T^\ast M$, and hence $X$ is a semilinear symmetric hyperbolic operator (cf.\ \cite[\S 2]{MR3846239} and \cite[Definition 5.1]{MR3302643}). 

Using the standard theory of semilinear symmetric hyperbolic systems, we therefore see that given initial data $u^0 \in H^k_{\sigma^0}(\iota_0^\ast \mathscr{U})$, the Cauchy problem 
\begin{equation}\label{eq-sym-hyp-cauchy}
    X(u) =0, \qquad \iota_0^\ast u = u^0
\end{equation}
has a unique local solution $u \in C\left([0,\varepsilon), \; H^k_{\omega^0}(\mathscr{U})\right)$ for $k\geq 2$.

\begin{Bem}
    The standard theory as presented in \cite{MR4703941} actually only concerns functions on $\R^n$ rather than sections of bundles, so one should technically take an atlas of trivializations of the manifold $M$ and the bundle $\mathscr{U}$, apply the classical theory within each trivialization, and then patch them together. The exact choice of atlas is immaterial here provided that it preserves the regularity of the coefficients sufficiently well (the assumption of bounded geometry allows for coordinates with good bounds on the metric coefficients and connection symbols, cf.\ \cite{MR3126616, MR2343536}), although in the present setting the most convenient choice is probably an atlas of adapted gauges (cf.\ \S \ref{sec-spatial-principal-bundle})), since the temporal derivative on a vector bundle associated to $P$ via a representation $\rho$ is in any adapted gauge just the tensor derivative $\frac{\nabla_{\omega^0}}{\diff\tau} = \frac{\nabla}{\diff\tau}$.
    Finally, we would like to note that the theory from  \cite{MR4703941} also only produces solutions with regularity $k > 3/2 + 1$ (i.e.\ $k \geq 3$) since the most general quasilinear symmetric hyperbolic systems are considered, but it is well-known (see for example \cite{MR2183957})
  that the regularity can be lowered to $k > 3/2$ (i.e.\ $k\geq 2$) for semilinear symmetric hyperbolic systems by a simple modification of the argument.
\end{Bem}

\subsubsection{Initial data and constraints}

If the local solution produced in the previous section is to be consistent with the original theory, the initial data $u^0 \in \Gamma(\iota_0^\ast \mathscr{X})$ cannot be arbitrary and we need to also take the constraints (\ref{eq-local-sys-connection-constraint}--\ref{eq-local-sys-dirac-constraint}) into account.
Evaluating (\ref{eq-curvature-sum}) at $ij$, we get the curvature constraint 
\begin{equation}\label{eq-eta-constraint}
    (\star Q)_{ij} - (\star Q_{\omega^0})_{ij} - (D_{\omega^0}\eta)_{ij} + (D_{\omega^0}\eta)_{ji} - [\eta_i, \eta_j] = 0,
\end{equation}
which is equivalent to (\ref{eq-local-sys-connection-constraint}) and determines the initial datum $Q^0$ in terms of $\eta^0$ (and $\sigma^0$).
Similarly, using (\ref{eq-diff-codiff-components}) and Lemma \ref{lem-curvature-splitting}, the Bianchi constraint (\ref{eq-local-sys-bianchi-constraint}) is%
\footnote{Note that the terms of the form $\sff_{ij}E_k$ all cancel out upon skew-symmetrizing.
This constraint holds automatically at the initial hypersurface since (\ref{eq-eta-constraint}) essentially says that $B=\star Q$ is the curvature form of $\sigma^0+\eta^0$.}
\begin{equation}\label{eq-B-constraint}
    (D_{\omega^0} \star Q)_{ijk} + (D_{\omega^0} \star Q)_{jki} + (D_{\omega^0} \star Q)_{kij} + [\eta_{i}, (\star Q)_{jk}] + [\eta_{j}, (\star Q)_{ki}] + [\eta_{k}, (\star Q)_{ij}] = 0,
\end{equation}
and the Yang-Mills constraint (\ref{eq-local-sys-ym-constraint}) is
\begin{equation}\label{eq-E-constraint}
    \tensor{(D_{\omega^0} E)}{^k_k} + [\eta^k, E_k] + \Re \left\langle \frac{\nabla_{\omega^0}\Phi}{\diff\tau}, \rho_*\Phi \right\rangle - \frac12 \Im \langle \partial_\tau \cdot\Psi, \chi_*\Psi \rangle = 0.
\end{equation}
The Dirac constraint (\ref{eq-local-sys-dirac-constraint}) can be Clifford multiplied by $\partial_\tau$ to obtain
\begin{equation}\label{eq-dirac-constraint}
    \frac{\nabla_{\omega^0}\Psi}{\diff\tau} = \partial_\tau \cdot \left( e^k \cdot (D_{\omega^0} \Psi)(e_k) + e^k \cdot \chi_*(\eta_k)\Psi + \Y_\Phi \Psi \right),
\end{equation}
which determines the initial datum $\dot\Psi^0$ in terms of $\eta^0, \Phi^0, \Psi^0$ (and $\sigma^0$).
Finally, $Z=D_\omega\Phi$ and $S=D_\omega\Psi$ are predetermined by the choice of $\Phi^0$ and $\Psi^0$ (as they only involve derivatives tangential to $\Sigma$), explicitly $Z^0=D_{\sigma^0}\Phi^0 + \rho_*(\eta^0)\Phi^0$ and $S^0= D_{\sigma^0}\Psi^0 + \chi_*(\eta^0)\Psi^0$.
Thus the initial data vector $u^0$ is parametrized by the tuple $(\eta^0, E^0, \Phi^0, \dot\Phi^0, \Psi^0)$  where $E^0$ satisfies (\ref{eq-E-constraint}), as assumed in the statement of Theorem \ref{thm-main}.

\subsubsection{Propagation of constraints}

Now we need to also show that the constraints (\ref{eq-eta-constraint}--\ref{eq-dirac-constraint}) are preserved by the evolution.

First, define the tensor $G \in \Gamma(\pi^\ast \Lambda^2 \Sigma \otimes \Ad P)$ with components equal to the left-hand side of the curvature constraint (\ref{eq-eta-constraint}).
A simple computation shows that if (\ref{eq-eta-evolution}, \ref{eq-Q-evolution}) are satisfied, then
\begin{equation}\label{eq-curvature-constraint-propagation}
    \left(\frac{\nabla_{\omega^0} G}{\diff\tau}\right)_{ij} = \tensor{\sff}{_i^\ell} G_{\ell j} +  \tensor{\sff}{_j^\ell} G_{i \ell}. 
\end{equation}
This is a linear first order ODE for $G$, and consequently it has a unique solution.
It follows that the constraint (\ref{eq-eta-constraint}) is preserved, if satisfied initially.
The equations (\ref{eq-eta-evolution}, \ref{eq-eta-constraint}) together are equivalent to (\ref{eq-curvature-sum}), and we therefore see that $F$ coincides with the curvature form of $\omega$  (i.e.\ we have $E_\omega = E$ and $B_\omega = B$).
Consequently, the constraint (\ref{eq-B-constraint}) will also be satisfied automatically, being the spatial part of the Bianchi identity (\ref{eq-local-sys-bianchi-constraint}).
Both of these facts are true for any $E, \Phi, \Psi$, not necessarily satisfying the Yang-Mills, Higgs and Dirac equations.

Now consider the constraint (\ref{eq-dirac-constraint}), which is just the original Dirac equation.
Note that (\ref{eq-local-sys-dirac}) is equivalent to 
\begin{equation*}
    \dirac_{\omega} \Theta = 0, \quad\text{where}\quad \Theta = \dirac_{\omega} \Psi + \Y_\Phi \Psi,
\end{equation*}
cf.\ Theorem \ref{thm-wave-eqs} and recall that we have already shown that $F_\omega = F$ above.
Applying $\dirac_{\omega}$ one more time, we see that $\Theta$ satisfies the linear wave equation
\begin{equation}\label{eq-dirac-constraint-propagation}
    \Box_{\omega} \Theta = \frac14 \Scal_h\Theta + \chi_*(F_\omega) \cdot \Theta,
\end{equation}
and thus $\Theta \equiv 0$, since it vanishes initially.

Now it only remains to show that the Yang-Mills constraint (\ref{eq-E-constraint}) is preserved.
To see this, define the one-form $C$ by
\begin{equation*}
    C  = \diff \tau \otimes C_0 + \iota^\ast C = \diff_\omega^\ast F - \mathfrak J.
\end{equation*}

Then $\iota^\ast C = 0$ is equivalent to the Yang-Mills evolution equations (\ref{eq-E-evolution-intermsofQ}) while the vanishing of $C_0$ is the Yang-Mills constraint (\ref{eq-E-constraint}).
By (\ref{eq-curvature-divergence-free}) and the fact that $F = F_\omega$, we get that
\begin{equation*}
    - \diff_\omega^\ast \mathfrak J = \diff_\omega^\ast C = -\tensor{(\nabla_\omega C)}{^\mu_\mu} = \frac{\nabla_\omega C_0}{\diff\tau } - 3HC_0 + D_\omega^\ast (\iota^\ast C),
\end{equation*}
and since $\iota^\ast C = 0$ by construction, we see that $C_0$ satisfies
\begin{equation}\label{eq-ym-constraint-propagation}
    \frac{\nabla_{\omega^0} C_0}{\diff\tau} = 3HC_0 - \diff_\omega^\ast \mathfrak J = 3HC_0,
\end{equation}
where the twisted divergence $\diff_\omega^\ast\,\mathfrak{J}$ vanishes when the Higgs and the (original) Dirac equations are satisfied, by Proposition \ref{prop-currents-div-free}.
Consequently the Yang-Mills constraint (\ref{eq-E-constraint}) will be preserved, if satisfied initially.

\begin{Bem}
    There is a small subtlety involving the propagation of constraints, namely that the equations, in particular the Dirac constraint wave equation (\ref{eq-dirac-constraint-propagation}), have coefficients that are only Sobolev-regular rather than smooth if one directly starts by producing the solution $u$ as in the previous section and only then considers these equations.
    This can be circumvented either by smoothening the rough coefficients using a mollifier, or more simply by instead simultaneously solving the Cauchy problem (\ref{eq-sym-hyp-cauchy}) for $u$ and the constraint problem
    \begin{equation*}
        K(v) = 0, \qquad v^0 = 0,
    \end{equation*}
    where $K$ is the differential operator associated to (\ref{eq-curvature-constraint-propagation}, \ref{eq-ym-constraint-propagation}, \ref{eq-dirac-constraint-propagation}) with respect to $v = (G, C_0, \Theta, \dot\Theta, D_\omega\Theta)$, and the wave equation (\ref{eq-dirac-constraint-propagation}) is reduced to a first order system as before.
    The previously rough coefficients are then considered as unknown variables of the system.
    The combined operator $X \oplus K$ is still semilinear and symmetric hyperbolic, giving a unique local solution $u \oplus v$.
    In particular, $u$ matches the previously constructed solution, while $v \equiv 0$ clearly solves the supplemented constraint problem.
    Since the tuple $(G, C_0, \Theta, \dot\Theta, D_\omega\Theta)$ defined in terms of the solution $u$ solves the same equations as $v$, the constraints get propagated by uniqueness. 
\end{Bem}   

This finishes the proof of the local existence and uniqueness of a solution $(\eta, \Phi, \Psi)$.

\subsection{Global existence}

Now consider the total energy of the Standard Model, cf.\ \S \ref{sec-sm-energy}.
Note however that this notion of energy depends on the choice of connection $\omega$, and since $\omega = \omega^0 + \eta$ is one of the variables of the system, we cannot use the proven energy estimate directly.
We therefore modify the argument by instead considering the energy with respect to the connection induced by $\omega^0$, and we also add the norm of $\eta$, i.e.\
\begin{equation*}
    \mathfrak{E}_{k+1}(\omega^0) := \mathfrak{E}_{k+1}(\eta; \omega^0) + \mathfrak{E}^+_{k+1}(F; \omega^0) + \mathfrak{E}_{k+1}(\Phi; \omega^0) + \mathfrak{E}^+_{k+1}(\Psi; \omega^0),
\end{equation*}
where
\begin{equation*}
    \mathfrak{E}_{k+1}(\eta; \omega^0) = \Vert \eta \Vert_{H^{k+1}_{\omega^0}}^2. 
\end{equation*}
We note that we can make the energy $\mathfrak{E}_{k+1}(\omega^0)$ as small as we like by choosing the initial data (in particular $\epsilon$ as in the statement of Theorem \ref{thm-main}) sufficiently small.
Indeed, the $H_{\omega^0}^{k+1}$-norms of 
\begin{equation}\label{eq-H-k+1-vars}
\eta, \quad  E, \quad B=\star\, Q, \quad \Phi, \quad \ \Psi    
\end{equation}
are then clearly small initially by assumption (note that $\iota_0^\ast B = F_{\sigma^0+\eta^0}$), while the temporal derivatives
\begin{equation}\label{eq-H-k-vars}
    \frac{\nabla_{\omega^0} E}{\diff\tau}, \quad \frac{\nabla_{\omega^0} B}{\diff\tau}, \quad \frac{\nabla_{\omega^0} \Phi}{\diff\tau},  \quad \frac{\nabla_{\omega^0} \Psi}{\diff\tau}
\end{equation}
have small $H_{\omega^0}^{k}$-norms by assumption and the identities (\ref{eq-E-evolution-intermsofQ}, \ref{eq-Q-evolution}, \ref{eq-dirac-constraint}).
We now wish to show that this energy stays uniformly bounded for sufficiently small $\epsilon$, along the lines of Theorem \ref{thm-total-energy-estimate}. 

To estimate the energy of $\eta$, we use Lemma \ref{lem-Hk-norm-evolution} and the evolution equation (\ref{eq-eta-evolution}) for $\eta$ to find that, for $k\geq 2$,
\begin{align}\nonumber
    \frac{\diff}{\diff\tau} \mathfrak{E}_{k+1}(\eta; \omega^0) 
    &\lesssim \left(  \left\Vert \frac{\nabla_{\omega^0}\eta}{\diff \tau} \right\Vert_{H^{k+1}_{\omega^0}} + \Vert R^{\pi^\ast T^\ast \Sigma \otimes \Ad P}_{\omega^0}(\partial_\tau, \cdot)\eta \Vert_{H^{k}_{\omega^0}} + \Vert \sff \Vert_{C^{k+1}} \Vert \eta \Vert_{H^{k+1}_{\omega^0}} \right) \Vert \eta \Vert_{H^{k+1}_{\omega^0}}\\ \label{eq-eta-Hk-estimate}
    &\lesssim \mathfrak{E}_{k+1}(E; \omega^0)^{\frac12} \mathfrak{E}_{k+1}(\eta; \omega^0)^{\frac12} + \Vert \sff \Vert_{C^{k+1}}  \mathfrak{E}_{k+1}(\eta; \omega^0),
\end{align}
where $R^{\pi^\ast T^\ast \Sigma \otimes \Ad P}_{\omega^0}(\partial_\tau,\cdot)\eta = D\sff \ast \eta$ by the results of \S \ref{sec-ym-bundle-curvatures} (recall that $E_{\omega^0} = 0$).

To estimate the other terms, we note that if $\xi$ is a section of a bundle associated to $P$ via a representation $\rho$ (or more generally a spatial form or spinor with values in such a bundle), then
\begin{equation*}
    \Box_\omega \xi = \Box_{\omega^0} \xi - \rho_*(D_{\omega^0}^\ast \eta) \xi + \rho_*(\eta^k)(D_{\omega^0}\xi)_k + \rho_*(\eta^k)\rho_*(\eta_k)\xi. 
\end{equation*}
In particular, if $k\geq 2$, we have the bound
\begin{equation*}
    \Vert \Box_{\omega^0} \xi \Vert_{H^{k}_{\omega^0}} \lesssim \Vert \Box_\omega \xi \Vert_{H^{k}_{\omega^0}} + \Vert \eta \Vert_{H^{k+1}_{\omega^0}}^2 \Vert \xi \Vert_{H^{k+1}_{\omega^0}},
\end{equation*}
and therefore (\ref{eq-energy-estimate-general}) gives
\begin{align*}
    \frac{\diff}{\diff\tau} \mathfrak{E}_{k+1}(\xi;\omega^0) 
    \lesssim&\, \left( 1 + \Vert {\Riem_g} \Vert_{C^{k}} + \Vert \sff \Vert_{C^{k+1}} \right)\mathfrak{E}_{k+1}(\xi;\omega^0)\\
    &\, + \left(\Vert \Box_{\omega} \xi \Vert_{H^{k}_{\omega^0}} 
    + \mathfrak{R}_{k}(\xi;\omega^0)
    \right) \mathfrak{E}_{k+1}(\xi;\omega^0)^{\frac12}
    + \mathfrak{E}_{k+1}(\eta;\omega^0)\mathfrak{E}_{k+1}(\xi;\omega^0).
\end{align*}
The estimates of the wave-terms from \S \ref{sec-higgs-energy-estimates}--\ref{sec-dirac-energy-estimates} carry over, word for word, to this setting, where the $H^k$-norms with respect to $\omega$ are simply replaced by ones with respect to $\omega^0$.
The estimates of the bundle curvature terms need only a simple modification, as the terms coming from $\mathfrak R_{k-1}$ will now depend on $F_{\omega^0} = \pi^\ast F_{\sigma^0}$ rather than $F_\omega$, and we can more simply estimate them by the $C^k_{\omega^0}$-norm. 
In summary, we get that the total energy with respect to $\omega^0$ satisfies an inequality of the form
\begin{align*}
    \frac{\diff}{\diff\tau} \mathfrak{E}_{k+1}(\omega^0) 
    \lesssim& 
    \left( 1 + \Vert {\Riem_g} \Vert_{C^{k}} + \Vert \sff \Vert_{C^{k+1}} + \Vert F_{\omega^0} \Vert_{C^{k}_{\omega^0}} \right)\mathfrak{E}_{k+1}(\omega^0)\\
    &\, + C(k,\sff) \mathfrak{E}_{k+1}(\omega^0)^{\frac32} + C(k,\sff) \mathfrak{E}_{k+1}(\omega^0)^2.
\end{align*}
Analogously as in the proof of Theorem \ref{thm-total-energy-estimate}, we see that the total energy stays uniformly bounded provided that it is sufficiently small initially, and therefore the solution exists future globally on $[0,T) \times \Sigma$ by a standard continuation argument.
In particular, the variables (\ref{eq-H-k+1-vars}) persist in $H^{k+1}_{\omega^0}$ while the variables (\ref{eq-H-k-vars}) persist in $H^{k}_{\omega^0}$.
This completes the global existence part of Theorem \ref{thm-main}.

\begin{Bem}
    At this point we could conclude also that the total energy with respect to $\omega = \omega^0+\eta$ (rather than $\omega^0$) is also uniformly bounded, but this is not needed for the proof.
\end{Bem}

\begin{Bem}
    \label{rem-optimal-regularities}
    Throughout the proof of Theorem \ref{thm-main}, we need to assume that the connection and the curvature have the same Sobolev regularity.
    In contrast, the energy estimate from Theorem \ref{thm-total-energy-estimate} holds even if the curvature has regularity lowered by one.
    One would perhaps expect that the latter interplay between regularities is more natural since curvature is in some sense the derivative of the connection.
    However, using our techniques here we cannot conclude that $\eta$ has regularity $H^{k+1}_{\omega^0}$ if the curvature has regularity $H^k_{\omega^0}$. 
    In particular, the evolution equation (\ref{eq-eta-evolution}) only allows us to conclude that $\eta$ has the same regularity as $E=E_\omega$.
    On the other hand, the formula (\ref{eq-B-constraint}) only allows us to uniformly bound the $H^k_{\sigma^0}$-norm of the exterior derivative $\diff_{\omega^0}\eta$ (i.e.\ the components $(D_{\omega^0}\eta)_{ij}-(D_{\omega^0}\eta)_{ji})$, while the formula (\ref{eq-E-constraint}) together with the evolution equation (\ref{eq-eta-evolution}) gives a uniform estimate of the $H^k_{\omega^0}$-norm of the codifferential $D_{\sigma^0}^\ast\eta$. These alone are however not sufficient to get a uniform bound for the $H^{k}_{\omega^0}$-norm of the covariant derivative $D_{\omega^0}\eta$. 
    In the seminal work on Yang-Mills theory on flat Minkowski space, Eardley and Moncrief \cite{MR649158} circumvent this issue by considering the decomposition of (Sobolev regular) vector fields into transversal and longitudinal components, and separately analyzing the equations for each of the parts.
    However, the proof of the well-definedness of such decompositions relies on the geometry of Minkowski space and does not carry over to globally hyperbolic manifolds in general. Nevertheless, it may still be possible to optimize the regularity in our setting if one could show the existence of such decompositions and mimic the techniques from \cite{MR649158}.
\end{Bem}

\subsection{Construction of the bundle automorphism}

To complete the proof of Theorem \ref{thm-main} we need to show the existence and uniqueness of the desired bundle automorphism $f : P \to P$, such that the transformed conformal Standard Model triplet $(f\omega, f\Phi, f\Psi)$ has a prescribed temporal component, i.e.\
\begin{equation*}
    (f\omega)(\partial_\tau) := \alpha \in C^0(I, H^{k+2}_{\omega^0}(\Ad P)).
\end{equation*}
Note that such $\alpha$ is in $C^0(I, C^1_{\omega^0}(\Ad P))$ by the Sobolev embedding and in particular continuous.%
\footnote{Note that the metric $g$ and the connection $\omega^0$ have bounded geometry only up to order $k+1$, so we can in particular only use Sobolev embeddings up to order $k+1$, and we cannot (nor need to) conclude that $\alpha\in C^0(I, C^2_{\omega^0}(\Ad P))$}
We can decompose the variation of a bundle automorphism $f : P \to P$ (cf.\ \S \ref{sec-bundle-aut}) as
\begin{equation}\label{eq-aut-derivative-decomposition}
    f\omega - \omega =: \autdiff_\omega f = \diff t \otimes \frac{\autdiff_\omega f}{\diff t} + \autdiff^\Sigma_\omega f 
    = \diff t \otimes \frac{\autdiff_{\omega^0} f}{\diff t} + \autdiff^\Sigma_{\omega^0} f + \left(f-\Id\right)\eta,
\end{equation}
where $\autdiff^\Sigma_\omega f = \autdiff_\omega f |_{TP^\Sigma}$ is the spatial part.
Our goal is therefore to construct a bundle automorphism $f$ such that
\begin{equation*}
    \frac{\autdiff_{\omega^0} f}{\diff t} = \alpha, \qquad  
    f\eta+\autdiff_{\omega^0}^\Sigma f \in C^0(I, H_{\omega^0}^{k+1}(\pi^\ast T^\ast \Sigma \otimes \Ad P)),
\end{equation*}
since then $f^\ast \omega$ will have the desired temporal component, while its spatial part will have the same regularity as $\omega$, i.e.\ $(f\omega)|_{TP^\Sigma} \in C^0([0,T), \mathcal{H}^{k+1}(P^\Sigma))$ (cf.\ also \cite[Proposition 3.1]{MR1458665}).

As in \S \ref{sec-bundle-aut}, we denote by $g : P \to G$ the map associated to the automorphism $f$ via $f(p) = p\cdot g(p)$.
By (\ref{eq-aut-derivative-formula}), we have 
\begin{align*}
    \frac{\autdiff_{\omega^0} f}{\diff t}
    = -\diff R_{g^{-1}} \circ \diff g(\partial_\tau),
\end{align*}
and therefore the map $g : P \to G$ must satisfy the linear ODE
\begin{equation*}
     \dot g{(t,p)} = -\diff R_{g(t,p)}(\alpha(t,p)), \qquad g(0,p) = e, 
\end{equation*}
for each fixed $p \in P^\Sigma$, where we interpret $\alpha$ at bundle level as a map $P \to \g$ of $\Ad$-type.
Since $\alpha$ is continuous, this ODE has a unique solution by standard theory (cf. \cite[Satz 1.10]{BAUM}), and one can easily show that the map $g$ satisfies 
\begin{equation*}
    g(t,p\cdot h) = h^{-1}\cdot g(t,p)\cdot h, \qquad h \in G,    
\end{equation*}
by differentiating with respect to $t$ for fixed $p \in P^\Sigma$ and using that the identity holds at $t=0$, so that $g$ indeed induces a bundle automorphism $f$.

Now we want to show show that $\autdiff_{\omega^0}^\Sigma f$ has the correct regularity (i.e. that it is uniformly in $H^{k+1}_{\omega^0}$).
Specializing (\ref{eq-aut-derivative-commutator}), we get that
\begin{align}\nonumber
    \left(\frac{\nabla_{\omega^0}}{\diff \tau} \autdiff^\Sigma_{\omega^0} f\right)(e_i)
    &= (\nabla_{\omega^0} \, \autdiff_{\omega^0} f) (\partial_\tau, e_i) \\\nonumber
    &= (\nabla_{\omega^0} \, \autdiff_{\omega^0} f) (e_i, \partial_\tau) + (f - \Id) (\partial_\tau \iprod F_{\omega^0})_i - \left[ \tfrac{\autdiff_{\omega^0} f}{\diff \tau}, (\autdiff^\Sigma_{\omega^0} f)_i \right]\\\label{eq-evolution-spatial-aut-derivative}
    &= (D_{\omega^0} \alpha)_i + \tensor{\sff}{^k_i} (\autdiff^\Sigma_{\omega^0} f)_k - \left[ \alpha, (\autdiff^\Sigma_{\omega^0} f)_i \right] 
\end{align}
where the curvature term is zero since $\omega^0 = \Pi^\ast \sigma^0$ is independent of $\tau$ and purely spatial.
By Lemma \ref{lem-Hk-norm-evolution}, we then see that
\begin{align*}
    \frac{\partial}{\partial \tau} \Vert \autdiff_{\omega^0}^\Sigma f \Vert_{H^{k+1}_{\omega^0}}
    &\lesssim \Vert \alpha \Vert_{H^{k+2}_{\omega^0}} + \left(\Vert \sff \Vert_{C^{k+1}} + \Vert \alpha \Vert_{H^{k+1}_{\omega^0}} \right)\Vert \autdiff_{\omega^0}^\Sigma f \Vert_{H^{k+1}_{\omega^0}} \\
    &\lesssim C(\alpha,\sff) \left( 1 + \Vert \autdiff_{\omega^0}^\Sigma f \Vert_{H^{k+1}_{\omega^0}} \right),
\end{align*}
and therefore, since $\iota_0^\ast (\autdiff_{\sigma^0}^\Sigma f) = 0$ by construction,
\begin{equation*}
    \sup_{0\leq\tau\leq T}\Vert \autdiff_{\omega^0}^\Sigma f \Vert_{H^{k+1}_{\omega^0}} \lesssim e^{C(\alpha,\sff)T} - 1.
\end{equation*}

Now finally, we need to show that $f\eta \in C^0([0,T), H_{\omega^0}^{k+1}(\pi^\ast T^\ast \Sigma \otimes \Ad P))$.
To this end we note that if $\xi$ is a section of an associated bundle $P \times_\rho W$, then
\begin{equation*}
    \Vert f \xi \Vert_{H^{k+1}_{\omega^0}} \lesssim \Vert \xi \Vert_{H^{k+1}_{\omega^0}} + \Vert \autdiff_{\omega^0}^\Sigma f \Vert_{H^{k}_{\omega^0}}^{k+1} \Vert \xi \Vert_{H^{k}_{\omega^0}},
\end{equation*}
which follows by iterating (\ref{eq-section-aut-derivative}).
Thus the uniform boundedness of the $H_{\omega^0}^{k+1}$-norm of $f\eta$ follows from the boundedness of the $H_{\omega^0}^{k+1}$-norm of $\autdiff_{\omega^0} f$.
This also shows that $f\Phi$ and $f\Psi$ belong to the same Sobolev space as $\Phi$ and $\Psi$.
Hence, $f$ is the desired bundle automorphism.

\begin{Bem}
    \label{rem-temporal-order-one-higher}
    Note that we needed to assume that the temporal coefficient $\alpha$ has one higher order of Sobolev regularity than the spatial part of the connection, but (\ref{eq-evolution-spatial-aut-derivative}) indicates that one cannot do better in general, at least with the techniques used here, cf.\ also Remark \ref{rem-optimal-regularities}.
\end{Bem}

\subsection{Conformal transformation back to the original spacetime}
\label{sec-global-original}

To wrap up, we would like to say a few words about the conformal transformation back to the original spacetime.
Let us, as in Definition \ref{def-expanding-spacetime}, denote the original spacetime by $(M,h)$ and the conformally transformed metric by $\tilde h = (Ns)^{-2} h$.
We recall from Proposition \ref{prop-conformal-trans} that the conformal Standard Model triplet should, after the conformal transformation (with $\Omega = (Ns)^{-1}$), be rescaled as
\begin{equation}\label{eq-conformal-rescaling}
    \widetilde\omega = \omega, \qquad \widetilde\Phi = Ns \, \Phi, \qquad \widetilde \Psi = (Ns)^{\frac32} \, \Psi.
\end{equation}
Note that we have also used the Gaussian temporal coordinate $\tau$ rather than the original coordinate $t$ in the last two sections. The two are related via $\diff \tau = s^{-1} \diff t$, so certain quantities that we have considered will also pick up an additional factor of $s$ when expressed in terms of $t$ (e.g.\ the temporal coefficients/derivatives).
In particular, we note that
\begin{equation}\label{eq-E-tilde-vs-E}
    \widetilde E_{\omega} = s E_\omega, \qquad \widetilde B_{\omega} = B_\omega
\end{equation}
where $\widetilde F_{\omega} = F_{\widetilde\omega} = \diff \tau \wedge \widetilde E_{\omega} + \widetilde B_{\omega}$ and $F_\omega = \diff t \wedge E_\omega + B_\omega$.

We note that, by the Sobolev embedding (which applies since $\tilde g$ and $\omega^0$ have uniformly bounded geometry \cite{MR1205820}), the solution $(\eta, \widetilde\Phi, \widetilde\Psi)$ has uniformly bounded supremum norm.
In particular, since $N$ is uniformly lower and upper bounded it follows from (\ref{eq-conformal-rescaling}) that $\Phi$ decays uniformly with rate $s^{-1}$ while $\Psi$ decays uniformly with rate $s^{-\frac32}$ for $t\to\infty$, and $E_\omega$ decays uniformly with rate $s^{-1}$ for $t\to\infty$ by (\ref{eq-E-tilde-vs-E}).

Finally, let us briefly analyze the behaviour of the $H^k$-norms of the different quantities under the conformal transformation.
\begin{itemize}
    \item If $\xi$ is a section of a bundle $E = P \times_\rho W$ associated to $P$, then since the natural bundle metric on such bundles does not depend on $g$ or $\tilde g$, the $L^2$-norm transforms as
    \begin{equation*}
        s(t)^{-\frac32} \Vert \xi \Vert_{L^2(g_t)} \lesssim \Vert \xi \Vert_{L^2(\tilde g_t)} \lesssim s(t)^{-\frac32} \Vert \xi \Vert_{L^2(g_t)},
    \end{equation*}
    since $\dv_{\tilde{g}_t} = (Ns)^{-3} \dv_{g_t}$
    and $N$ is uniformly upper and lower bounded.
    The same holds if $\xi$ is a twisted spinor, due to our identification of spin structures (cf.\ Table \ref{tab-conformal-transform}).
    \item If $\theta \in \Omega^p(\Sigma)$, then $|\theta|_{\tilde g}^2 = (Ns)^{2p} |\theta|_{g}^2$, and hence
    \begin{equation*}
         s(t)^{p-\frac32} \Vert \theta \Vert_{L^2(g_t)} \lesssim \Vert \theta \Vert_{L^2(\tilde g_t)} \lesssim s(t)^{p-\frac32} \Vert \theta \Vert_{L^2(g_t)}.
    \end{equation*}
    \item If $\xi$ is a $p$-covariant spatial tensor with values in some tensor product of associated bundles to $P$ and/or the spinor bundle, then the covariant derivatives transform as (cf.\ Table \ref{tab-conformal-transform})
    \begin{equation*}
        \widetilde{D}^\ell \xi = D^\ell \xi + \sum_{i,j=0}^{\ell-1} (D^{i+1} \log N) \ast D^j\xi,
    \end{equation*}   
    and since the derivatives of $N$ have uniformly bounded $C^k(g)$-norm, we have
    \begin{equation*}
         \Vert \xi \Vert_{H^k(\tilde g_t)} \lesssim s(t)^{p-\frac32} \sum_{\ell=0}^k  s(t)^\ell \Vert \xi \Vert_{H^\ell(g_t)} \lesssim s(t)^{k+p-\frac32} \Vert \xi \Vert_{H^k(g_t)},
    \end{equation*} 
    where we use that  $s$ is lower bounded in the final estimate.
    Replacing $\tilde g$ by $g$ (and hence $s$ by $s^{-1}$) in the first inequality above, one similarly also gets a lower bound for the norm with respect to $\tilde g$ where one instead estimates $\sum_{\ell=0}^k s^{-\ell} \lesssim 1$, and we conclude that
    \begin{equation*}
        s(t)^{p-\frac32} \Vert \xi \Vert_{H^k(g_t)} \lesssim \Vert \xi \Vert_{H^k(\tilde g_t)} \lesssim s(t)^{k+p-\frac32} \Vert \xi \Vert_{H^k(g_t)}.
    \end{equation*}
    More generally, if we put $\widetilde\xi = (Ns)^{\gamma} \xi$ for some $\gamma \in \R$, then by similar arguments 
    \begin{equation*}
        s(t)^{\gamma+p-\frac32} \Vert \xi \Vert_{H^k(g_t)}
        \lesssim
        \Vert \widetilde{\xi} \Vert_{H^k(\tilde g_t)} 
        \lesssim 
        s(t)^{\gamma+k+p-\frac32} \Vert \xi \Vert_{H^k(g_t)}.
    \end{equation*}
    \item Finally, temporal derivatives transform as
    \begin{equation*}
        \frac{\widetilde\nabla\xi}{\diff \tau} = s \left(\frac{\nabla \xi}{\diff t} + \left(\partial_t(\log N) + s^{-1}\dot s\right) \xi \right).
    \end{equation*}
    Hence, assuming that $s^{-1}\dot s$ is bounded,
    \begin{equation*}
        s(t)^{\gamma+p-\frac12} \left( \left\Vert \frac{\nabla \xi}{\diff t} \right\Vert_{H^k(g_t)} - \Vert \xi \Vert_{H^k(g_t)}  \right)
         \lesssim
        \left\Vert \frac{\widetilde\nabla \widetilde\xi}{\diff\tau} \right\Vert_{H^k(\tilde g_t)}
        \lesssim s(t)^{\gamma+k+p-\frac12} \left( \left\Vert \frac{\nabla \xi}{\diff t} \right\Vert_{H^k(g_t)} + \Vert \xi \Vert_{H^k(g_t)}  \right).
    \end{equation*}
\end{itemize}

The $H^k_{\sigma^0}$-norms at the initial hypersurface $t = 0$ (equivalent to $\tau = 0$) are thus equivalent with respect to $g$ and $\tilde g$, depending only on the bounds for $N$ and the fixed positive number $s(0)$. In particular, initial data on $(M,h)$ induce initial data on $(M,\tilde h)$, and their smallness in the appropriate $H^k_{\sigma^0}$-norms (without weights) is equivalent up to a constant.
Asymptotically, the norms differ however, and we can only conclude that the global solution elements satisfy
\begin{align*}
    \eta  &\in C^0\left([0,\infty), s^{-\frac12} H^{k+1}_{g,\omega^0}(\pi^\ast T^\ast \Sigma \otimes \Ad P)\right),\\
    E_\omega &\in C^0\left([0,\infty), s^{\frac12} H^{k+1}_{g,\omega^0}(\pi^\ast T^\ast \Sigma \otimes \Ad P)\right),\\
    B_\omega &\in C^0\left([0,\infty), s^{\frac12} H^{k+1}_{g,\omega^0}(\pi^\ast \Lambda^2 \Sigma \otimes \Ad P)\right),\\
    \Phi &\in C^0\left([0,\infty), s^{-\frac12} H^{k+1}_{g,\omega^0}(\Hig)\right),\\
    \Psi &\in C^0\left([0,\infty), H^{k+1}_{g,\omega^0}(\Fer_+)\right),
\end{align*}
while if $s^{-1}\dot s$ is also assumed to be bounded, then the above also holds with $C^0$ replaced by $C^1$ and $H^{k+1}$ replaced by $H^k$.

\begin{Bem}
    We would like to point out that, other than using an energy method as highlighted in this paper, one could likely also reach a global existence result in weighted Sobolev spaces using the conformal extension technique \`a la Ginoux and M\"uller \cite{MR3846239} after sorting out some technical details. To do this, one would have to extend the Gaussian foliation across $\tau = T$, which should be possible in view of our assumptions of bounded geometry and the Hamilton compactness theorem (cf.\ \cite[\S 7]{MR2265040}). Our technique has the benefit of not requiring any such extensions, and gives a more explicit and global description of the Sobolev spaces involved.
\end{Bem}

\begin{appendix}

\section{Euler-Lagrange equations}
\label{appendix-el-eqs-derivation}

For completeness, we prove here Theorem \ref{thm-el-eqs}, i.e.\ we derive the Euler-Lagrange equations corresponding to the Lagrangian 
\begin{equation*}
    L: (\omega,\Phi,\Psi) \mapsto 
    -\left(\vert F_\omega \vert^2  + \vert \nabla_{\omega} \Phi \vert^2 + U(\Phi) + 
    \Re\langle \Psi, i(\dirac_{\omega} + \Y_\Phi) \Psi \rangle\right).
\end{equation*}

More precisely, fix a semi-Riemannian manifold $(M,g)$ and look for $(\omega,\Phi,\Psi)$ such that for any compact $K\subset M$ and variation
\begin{equation*}
    (\omega(s), \Phi(s), \Psi(s)) = (\omega,\Phi,\Psi) + s(\xi, \alpha, \beta)
\end{equation*}
such that $(\xi,\alpha,\beta)$ has support contained in $K$, we have
\begin{equation*}
    \frac{\diff }{\diff s}\bigg|_{s=0} \int_K L(\omega(s),\Phi(s),\Psi(s)) \,\dv_g = 0.
\end{equation*}

First off all, we observe that (keeping $\omega$ and $\Phi$ fixed)
\begin{align*}
    \frac{\diff }{\diff s}\bigg|_{s=0} \int_K  \Re\langle \Psi(s), i(\dirac_\omega + \Y_\Phi)\Psi(s) \rangle \dv_g
    =  2\Re \int_K \left\langle \beta, i(\dirac_\omega + \Y_\Phi)\Psi \right\rangle\dv_g
\end{align*}
since both $i\dirac_\omega$ and $i\Y_\Phi$ are self-adjoint with respect to the $L^2$ inner product.
Next, keeping also $\omega$ fixed, we have 
\begin{align*}
    &\frac{\diff }{\diff s}\bigg|_{s=0} \int_K \left( |\nabla_{\omega} \Phi(s)|^2 + U(\Phi(s)) + \langle \Psi, i\Y_{\Phi(s)}\Psi \rangle \right) \dv_g\\[0.2cm]
    =& \, \int_K \left( 2\Re\langle \nabla_{\omega} \alpha, \nabla_{\omega} \Psi \rangle + \diff U_\Phi(\alpha) + \langle \Psi, i\Y_\alpha\Psi \rangle \right) \dv_g\\[0.2cm]
    =& \, -2\Re \int_K \left\langle \alpha, \Box_{\omega} \Phi - \tfrac12\grad U_\Phi - \langle \Psi, i\Y^- \Psi \rangle \right\rangle \dv_g,
\end{align*}
where 
\begin{enumerate}
    \item $\Box_{\omega} = -(\nabla_{\omega})^\ast \nabla_{\omega}$ is the connection d'Alembertian on $\Hig$,
    \item $\Y^-$ is the complex antilinear part of $\Y^-$, i.e.\ the corresponding term is given in an orthonormal frame $W_K$ of $\Hig$ by
    \begin{equation*}
        \langle \Psi, \Y^- \Psi \rangle = \sum_{k=1}^{\dim W} \langle \Psi, \tfrac12(\Y_{W_k} + i\Y_{iW_k})\Psi \rangle_{\Fer}\, W_k,
    \end{equation*}
    \item $\grad U_\Phi \in \Gamma(\Hig)$ is the "bundle gradient" of $U$ evaluated at $\Phi$, satisfying 
    \begin{equation*}
        \Re \langle \alpha, \grad U_\Phi \rangle = \diff U_\Phi(\alpha),
    \end{equation*}
    for all $\alpha \in \Gamma(\Hig)$ (note that we are using the canonical tangent bundle decomposition $T\Hig = TM \oplus \Hig$ here). 
\end{enumerate}
Now it only remains to calculate the variation with respect to $\omega$.
The difference of two connections can be viewed as an element $\xi$ of $\Omega^1(M,\Ad(P))$. 
Then we can write
\begin{align*}
    F_{\omega(s)} &= F_\omega + s\,\diff _\omega \xi + \frac{s^2}{2} [\xi \wedge \xi], \\[0.2cm]
    \nabla_{\omega(s)} \Phi &= \nabla_{\omega} \Phi + s\, \rho_*(\xi) \Phi, \\[0.2cm]
    \dirac_{\omega(s)}\Psi &= \dirac_\omega \Psi + s\, \chi_*(\xi) \cdot \Psi,
\end{align*}
and hence
\begin{align*}
    &\frac{\diff }{\diff s}\bigg|_{s=0} \int_K |F_{\omega(s)}|^2 \,\dv_g 
    = 2 \Re\int_K \langle \diff _\omega \xi, F_\omega \rangle \,\dv_g 
    = 2 \Re\int_K \langle \xi, \diff _\omega^\ast  F_\omega \rangle \,\dv_g, 
    \\[0.2cm]
    &\frac{\diff }{\diff s}\bigg|_{s=0} \int_K |\nabla_{\omega(s)}\Phi|^2 \,\dv_g 
    = 2 \Re\int_K \langle \nabla_{\omega} \Phi, \rho_*(\xi) \Phi \rangle\,\dv_g
    = 2 \int_K \langle \xi, \Re\langle \nabla_{\omega} \Phi \otimes \rho_* \Phi \rangle \rangle\,\dv_g,
    \\[0.2cm]
    &\frac{\diff }{\diff s}\bigg|_{s=0} \int_K \Re\langle \Psi, i\dirac_{\omega(s)}\Psi \rangle \,\dv_g
    = -\Im \int_K \langle \Psi, \chi_*(\xi) \cdot \Psi \rangle \, \dv_g
    = -\int_K \langle \xi, \Im\langle \id \cdot \Psi \otimes \chi_*\Psi \rangle \rangle \,\dv_g,
\end{align*}
where the currents are given explicitly by
\begin{align*}
    \Re \langle \nabla_{\omega}\Phi \otimes \rho_*\Phi \rangle
        &= \Re \langle \nabla_{\omega}\Phi(e_\mu), \, \rho_*(\xi^a)\Phi \rangle_\Hig \, e^\mu \otimes \xi_a,\\[0.1cm]
        \Im \langle \id\cdot\Psi \otimes \chi_*\Psi \rangle
        &= \Im \langle e_\mu \cdot \Psi, \, \chi_*(\xi^a)\Psi \rangle_{\Fer}\, e^\mu \otimes \xi_a,
\end{align*}
for an orthonormal frame $\xi_a \in \Ad P$, and we also use the symmetry of Clifford multiplication in the last equality. Putting everything together gives Theorem \ref{thm-el-eqs}.

\section{Lichnerowicz-Weitzenböck formulae}

The goal of this section is to state some identities used throughout the article.
These are by no means new, but seem to be difficult to find in the literature.

In the following results, given a vector bundle $E$ with a connection $\nabla^E$, we denote by $R^E$ the \textit{curvature endomorphism} of $\nabla^E$, i.e.\ the map
\begin{equation*}
    R^E : \Gamma(TM \otimes TM) \to \Gamma(\End(E)), \quad R^E(X,Y)\xi = \nabla^E_X \nabla^E_Y \xi - \nabla^E_Y \nabla^E_X \xi - \nabla^E_{[X,Y]}\xi.
\end{equation*}
We can also view it as a mapping taking sections of $E$ to $E$-valued two-forms
\begin{equation*}
    R^E\xi = e^\mu \wedge e^\nu \otimes \tfrac12 R^E(e_\mu,e_\nu)\xi.
\end{equation*}
For the particular case $E=TM$ equipped with the Levi-Civita connection, we index the Riemann tensor as described in (\ref{eq-riemann-convention}).

\begin{Prop}[Lichnerowicz-Weitzenböck formula for bundle-valued spinors]\label{prop-dirac-lichnerowicz-weitzenbock}
    Let $(M,g)$ be a spin Lorentzian manifold and $E$ a vector bundle with a connection $\nabla^E$.
    Then the twisted Dirac operator $\dirac^E$ on $\Sigma M \otimes E$ satisfies
    \begin{equation*}
        (\dirac^E)^2 = - \Box^{\Sigma M \otimes E} + \tfrac14 \Scal^M + \slashed{\mathcal{R}}{}^E, 
    \end{equation*}
    where $\Scal^M$ is the scalar curvature of $(M,g)$, and $\slashed{\mathcal{R}}^E$ is the operator on $\Sigma M \otimes E$ satisfying
    \begin{equation*}
        \slashed{\mathcal{R}}{}^E (\sigma \otimes \xi) = \frac12 (e_\mu \cdot e_\nu \cdot \sigma) \otimes (R^E\xi)^{\mu\nu},
    \end{equation*}
    where $R^E$ is the curvature endomorphism of $E$.
\end{Prop}

\begin{Prop}[Weitzenböck formula for bundle-valued forms]\label{prop-twisted-weitzenbock}
    Let $E$ be a vector bundle with a connection $\nabla^E$. Then, on $\Omega^k(M,E)$, we have
    \begin{equation*}
        (\diff^E)^* \diff^E + \diff^E (\diff^E)^*
        = -\Box^{\Lambda^k\hspace{-0.04cm}M\otimes E} + \mathcal{R}^{\Lambda^k\hspace{-0.04cm}M\otimes E},
    \end{equation*}
    where the twisted Weitzenböck curvature operator is given by
    \begin{equation*}
        \mathcal{R}^{\Lambda^k\hspace{-0.04cm}M\otimes E}\xi = e^\mu \wedge e_\nu \iprod \tensor{(R^{\Lambda^k\hspace{-0.04cm}M \otimes E}\xi)}{^\nu_\mu}.
    \end{equation*}
\end{Prop}

We also derive some formulas for the Lichnerowicz-Weitzenböck curvature operators in specific cases.

\begin{Lem}\label{lem-weitzenbock-curvature-explicit}
\leavevmode
\begin{enumerate}
    \item The twisted Weitzenböck curvature operator satisfies
    \begin{align}
        \mathcal R^{\Lambda^{k+\ell}M} (\theta \wedge \eta) 
        =&\, \mathcal R^{\Lambda^kM} \theta \wedge \eta + \theta \wedge \mathcal R^{\Lambda^\ell M}\eta \nonumber\\[0.1cm]
        &+ \tensor{(R\theta)}{^\nu_\mu} \wedge e^\mu \wedge (e_\nu \iprod \eta) + e^\mu \wedge (e_\nu \iprod \theta) \wedge \tensor{(R\eta)}{^\nu_\mu}, \label{eq-twisted-weitzenbock-leibniz}
    \end{align}
    for $\theta \in \Omega^k(M), \eta \in \Omega^\ell(M)$.
    In particular, we have
    \begin{align*}
        \mathcal R^{T^\ast M} e^\alpha &= \tensor{R}{^\alpha_\mu}  e^\mu,\\[0.1cm]
        \mathcal R^{\Lambda^2M}(e^\alpha \wedge e^\beta) &= (\delta^\beta_\nu\tensor{R}{^\alpha_\mu} + \delta^\alpha_\mu \tensor{R}{^\beta_\nu}- 2\tensor{R}{^\alpha_{\mu\nu}^\beta}) \, e^\mu \wedge e^\nu.
    \end{align*}
    \item Assume $E = P \times_\rho V$ is the associated vector bundle for some principal $G$-bundle $P$ and representation $\rho : G \to \mathbf{GL}(V)$, and let $\omega \in \Omega^1(P,\g)$ be a connection on $P$. 
    If
    \begin{equation*}
        \xi = e^{\mu_1}\wedge\cdots\wedge e^{\mu_k} \otimes \tfrac{1}{k!}\,\xi_{\mu_1\ldots\mu_k} \in \Omega^k(M,E),
    \end{equation*}
    then
    \begin{align*}
        \slashed{\mathcal{R}}{}^E\xi &= \rho_*(F_\omega) \cdot \xi = \tfrac12\, e^\mu \cdot e^\nu \cdot \rho_*(F_{\mu\nu}) \xi, \\[0.2cm]
        \mathcal{R}^E\xi &=  e^{\mu_1}\wedge\cdots\wedge e^{\mu_k} \otimes \tfrac{1}{(k-1)!}\,\rho_*(\tensor{F}{^\lambda_{\mu_1}})\,\xi_{\lambda \mu_2\ldots \mu_k}.
    \end{align*}
\end{enumerate}
\end{Lem}

\section{Differentials of currents}
\label{sec-differentials-currents}

We provide a result that allows us to show the identities (\ref{eq-higgs-current-divergence},\ref{eq-dirac-current-divergence}), as well as some parts of Theorem \ref{thm-wave-eqs}.

\begin{Lem}\label{lemma-current-differential}
    Let $G$ be a Lie group equipped with an Ad-invariant metric.
    Let $G \to P \xrightarrow{\pi} M$ be a principal bundle with a connection $\omega$.
    Let $\rho : G \to \mathbf{GL}(W)$ be a representation of $G$ and consider the associated vector bundle $E = P \times_\rho W$, equipped with the metric connection $\nabla_\omega$ induced by $\omega$.
    \begin{enumerate}
        \item If $\sigma \in \Gamma(E)$, then we have
        \begin{equation*}
            \nabla_\omega (\rho_* \sigma) = \rho_* \nabla_\omega \sigma.
        \end{equation*}
        \item
        If $\theta \in \Omega^1(M, E)$ and $\sigma \in \Gamma(E)$, then we have
        \begin{align*}
            \diff_\omega^* \langle \theta \otimes \rho_* \sigma \rangle &= -\langle \tr\nabla_\omega \theta, \rho_*\sigma \rangle_E - \langle \theta, \rho_*\nabla_\omega \sigma \rangle_{T^\ast M \otimes E},\\[0.1cm]
            \diff_\omega \langle \theta \otimes \rho_* \sigma \rangle &= \langle \diff_\omega \theta, \rho_*\sigma \rangle_E - \langle \theta \wedge \rho_* \nabla_\omega \sigma \rangle_E.
        \end{align*} 
    \end{enumerate}
\end{Lem}

\begin{Bem}
As before, here we use the notations
\begin{align*}
    \rho_*(\theta \otimes \sigma) &:= \theta \otimes \rho_*(\xi^a) \sigma \otimes \xi_a \in \Gamma(T^\ast M \otimes E \otimes \Ad P)\\
    \langle \theta \otimes \rho_* \sigma \rangle &:= \langle \theta_\mu, \rho_*(\xi^a) \sigma \rangle_E \, e^\mu \otimes \xi_a \in \Omega^1(M, \Ad P),
\end{align*}
and we also introduce
\begin{equation*}
    \langle \theta \wedge \rho_* \nabla_\omega \sigma \rangle_E
    := \langle \theta_\mu, \rho_*(\xi^a) (\nabla_\omega \sigma)_\nu \rangle_E \; e^\mu \wedge e^\nu \otimes \xi_a.
\end{equation*}
\end{Bem}

\begin{proof}
    We only verify (i), as the other statements follow easily from it after some simple calculations.
    In a local gauge, we have that
    \begin{equation*}
        \nabla_\omega (\rho_*(\xi^a) \sigma) = \rho_*(\nabla_\omega\xi^a) \sigma + \rho_*(\xi^a) \nabla_\omega\sigma,
    \end{equation*}
    since $\rho_*$ is a constant linear map, and $\rho_*[X,Y] = [\rho_*(X),\rho(Y)]$ for all $X,Y \in \g$.
    Therefore
    \begin{equation*}
        \nabla_\omega (\rho_* \sigma) = \rho_* \nabla_\omega \sigma + \widetilde\sigma,
    \end{equation*}
    where 
    \begin{equation}\label{eq-divergence-residue}
        \tilde\sigma = \rho_*(\nabla_\omega\xi^a) \sigma \otimes \xi_a + \rho_*(\xi^a)\sigma \otimes (\nabla_\omega\xi_a).
    \end{equation}
    We wish to show that $\tilde\sigma = 0$.
    To see this, consider the connection symbols of $\nabla_\omega^{\Ad P}$ in the frame $\xi_a$, defined by
    \begin{equation*}
        \nabla_\omega \xi_a = \tensor{\gamma}{_a^b} \otimes \xi_b.
    \end{equation*}
    Since the $\xi_a$ are orthonormal and the connection on $\Ad P$ is metric, we see that the Christoffel symbols are skew-symmetric, so that $\tensor{\gamma}{_a^b} + \tensor{\gamma}{^b_a} = 0$.
    But then if we expand (\ref{eq-divergence-residue}) in terms of the symbols, we see directly that $\tilde \sigma \equiv 0$ by this skew-symmetry.
\end{proof}

\section{Splitting of the curvature form}
\label{appendix-curvature-splitting}

The goal of this section is to prove Proposition \ref{prop-E-B-wave-eqs}.

\begin{Lem}\label{lem-curvature-splitting}
    In terms of the decomposition $F_\omega = e^0 \wedge E_\omega + B_\omega$ as in \S \ref{sec-yang-mills-energy-definition}, we have
    \begin{align*}
        \frac{\nabla_\omega F_\omega}{\diff\tau} &= e^0 \wedge \frac{\nabla_\omega E_\omega}{\diff\tau} + \frac{\nabla_\omega B_\omega}{\diff\tau}, \\[0.2cm]
        (\nabla_\omega F_\omega)_k &= e^0 \wedge e^i \otimes \left( (D_\omega E_\omega)_{ki} + \tensor{\sff}{_k^\ell} B_{\ell i} \right)\\
        &+  e^i \wedge e^j \otimes \frac12 \left( (D_\omega B_\omega)_{kij} - \tensor{\sff}{_k_j} E_i + \tensor{\sff}{_k_i} E_j \right),
        \\[0.2cm]
        \Box_\omega F_\omega &= e^0 \wedge e^i \otimes \Big( (\Box_\omega E_\omega)_i  + 2|\sff|^2 E_i - \sff_{ik}\sff^{k\ell}E_\ell
        + 2\tensor{\sff}{^k^\ell} (D_\omega B_\omega)_{k\ell i} - \tensor{(D^\ast \sff)}{^k} B_{ki}  \Big)
        \\[0.1cm]
        &+ e^i \wedge e^j \otimes \frac12 \Big((\Box_\omega B_\omega)_{ij} + 2\tensor{\sff}{_i^k} (D_\omega E_\omega)_{kj} - 2\tensor{\sff}{_j^k} (D_\omega E_\omega)_{ki}
        \\[0.1cm]
        &\hspace{2.5cm} - (D^\ast \sff)_i E_j + (D^\ast \sff)_j E_i + \sff_{ik} \sff^{k\ell}B_{\ell j} - \sff_{jk} \sff^{k\ell}B_{\ell i} \Big).
    \end{align*}
\end{Lem}

\begin{proof}
    Throughout the proof for $F_\omega$ we use the connection on $\Lambda^2 M \otimes \Ad P$, while for the spatial forms we use the connections on $\pi^\ast T^\ast \Sigma \otimes \Ad P$ and $\pi^\ast \Lambda^2 \Sigma \otimes \Ad P$ respectively.
    For notational ease, we will drop the subscript $\omega$ from the twisted connections where appropriate.
    
    The first formula is trivial since $\nabla_{e_0} e_\mu = 0$ in the adapted frame (cf.\ \S \ref{sec-adapted-frame}). 
    The second formula follows from the calculations
    \begin{align*}
        (\nabla_\omega F_\omega)_{k0i}
        &= \nabla_{e_k}^{\Ad P} E_i - F_\omega (\nabla_{e_k} e_0, e_i) - F_\omega(e_0, \nabla_{e_k} e_i)
        \\[0.1cm]
        &= \nabla_{e_k}^{\Ad P} E_i + \tensor{\sff}{_k^\ell} B_{\ell i} + \sff_{ki} F_{00} - E_\omega( D_{e_k} e_i)
        \\[0.1cm]
        &= (D_\omega E_\omega)_{ki} + \tensor{\sff}{_k^\ell} B_{\ell i}
    \end{align*}
    and
    \begin{align*}
        (\nabla_\omega F_\omega)_{kij}
        &= \nabla_{e_k}^{\Ad P} B_{ij} - F_\omega (\nabla_{e_k} e_i, e_j) - F_\omega(e_i, \nabla_{e_k} e_j)
        \\[0.1cm]
        &= \nabla_{e_k}^{\Ad P} B_{ij} + \sff_{ki} E_j - F_\omega (D_{e_k} e_i, e_j) - \sff_{kj} E_i - F_\omega(e_i, D_{e_k} e_j)
        \\[0.1cm]
        &= (D_\omega B_\omega)_{kij} + \sff_{ki} E_j - \sff_{kj} E_i.
    \end{align*}
    On the other hand,
    \begin{align}
        \nonumber
        (\Box_\omega F_\omega)_{0i}
        &= \eta^{\mu\nu} (\nabla_\omega^2 F_\omega)(e_\mu, e_\nu, e_0, e_i)
        \\[0.1cm] \nonumber
        &= \eta^{\mu\nu} \Big( \nabla_{e_\mu}^{\Ad P} (\nabla_\omega F_\omega)_{\nu 0 i} - (\nabla_\omega F_\omega)(\nabla_{e_\mu}e_\nu, e_0, e_i)
        \\\nonumber
        &\hspace{1cm} - (\nabla_\omega F_\omega)(e_\nu, \nabla_{e_\mu}e_0, e_i)  - (\nabla_\omega F_\omega)(e_\nu, e_0, \nabla_{e_\mu}e_i)  \Big)
        \\ \nonumber
        &= -\left(\tfrac{\nabla_\omega^2 E_\omega}{\diff\tau^2}\right)_i + g^{k\ell} \Big( \nabla_{e_k}^{\Ad P} (\nabla_\omega F_\omega)_{\ell 0 i} - (\nabla_\omega F_\omega)(\nabla_{e_k}e_\ell, e_0, e_i)\\ \label{eq-box-F0i}
        &\hspace{1cm} - (\nabla_\omega F_\omega)(e_\ell, \nabla_{e_k}e_0, e_i)  - (\nabla_\omega F_\omega)(e_\ell, e_0, \nabla_{e_k}e_i)  \Big).
    \end{align}
    Using the first part of the result, we can write
    \begin{align*}
        \nabla_{e_k}^{\Ad P} (\nabla_\omega F_\omega)_{\ell 0 i}
        &= \nabla_{e_k}^{\Ad P} (D_\omega E_\omega)_{\ell i} + e_k(\tensor{\sff}{_\ell^a}) B_{ai} + \tensor{\sff}{_\ell^a}\,\nabla_{e_k}^{\Ad P} B_{a i}\\[0.1cm]
        &= (D_\omega^2 E_\omega)_{k \ell i} + (D_\omega E_\omega)(D_{e_k} e_\ell, e_i) + (D_\omega E_\omega)(e_\ell, D_{e_k} e_i)
        \\[0.1cm]
        &\quad + \tensor{(D\sff)}{_k_\ell^a} B_{ai}  +   g^{ab}\, \sff(D_{e_k} e_\ell, e_a) B_{bi}  + \tensor{\sff}{_\ell^a}(D_\omega B_\omega)_{kai} + g^{ab}\, \sff_{\ell a} B(e_b, D_{e_k} e_i),
    \end{align*}
    and
    \begin{align*}
        (\nabla_\omega F_\omega)(\nabla_{e_k}e_\ell, e_0, e_i)
        &= -\sff_{k\ell} \left(\tfrac{\nabla_\omega E_\omega}{\diff\tau}\right)_i + (D_\omega E_\omega)(D_{e_k} e_\ell, e_i) + g^{ab}\, \sff(D_{e_k} e_\ell, e_a) B_{bi} ,
        \\[0.1cm]
        (\nabla_\omega F_\omega)(e_\ell, \nabla_{e_k}e_0, e_i)
        &= -\tensor{\sff}{_k^m} (D_\omega B_\omega)_{\ell m i} +\tensor{\sff}{_k^m}\sff_{\ell m} E_i -\tensor{\sff}{_k^m}\sff_{\ell i} E_m,
        \\[0.1cm]
        (\nabla_\omega F_\omega)(e_\ell, e_0, \nabla_{e_k}e_i)
        &=
        (D_\omega E_\omega)(e_\ell, D_{e_k} e_i) + g^{ab}\, \sff_{\ell a} B(e_b, D_{e_k} e_i). 
    \end{align*}
    Combining these calculations, we get
    \begin{align*}
        \text{(\ref{eq-box-F0i})}
        &= -\left(\tfrac{\nabla_\omega^2 E_\omega}{\diff\tau^2}\right)_i + \tr(\sff) \left(\tfrac{\nabla_\omega E_\omega}{\diff\tau}\right)_i 
        + \tensor{(D_\omega^2 E_\omega)}{^\ell_\ell_i}  \\[0.1cm]
        &\quad + \tensor{(D\sff)}{^\ell_\ell^a} B_{ai} + 2\sff^{\ell a} (D_\omega B_\omega)_{\ell ai}
        + \sff^{k\ell}\sff_{k\ell} E_i - \sff^{\ell m}\sff_{\ell i} E_m.
        \\[0.1cm]
        &= (\Box_\omega E_\omega)_i - \tensor{(D^\ast \sff)}{^\ell} B_{\ell i} + 2\sff^{k\ell} (D_\omega B_\omega)_{k\ell i}
        + 2|\sff|^2 E_i - \sff^{k\ell}\sff_{ki} E_\ell.
    \end{align*}
    The formula for $(\Box_\omega F_\omega)_{ij}$ is obtained in an entirely similar fashion.
\end{proof}

Next, we decompose all the terms on the right-hand side of (\ref{eq-yang-mills-wave}).
First of all, it is trivial to see that 
\begin{equation*}
    \langle [F_\omega \wedge F_\omega] \rangle 
    = e^0 \wedge e^i \otimes 2[E^k, B_{ik}] + e^i \wedge e^j \otimes (-[E_i, E_j] + [\tensor{B}{^k_i}, B_{kj}]).
\end{equation*}

Next we calculate $\mathcal R F_\omega$.
If we first decompose $\mathcal R F_\omega$ into temporal and spatial components, 
and then apply the identities from Lemma \ref{lem-curvature-generalized-cyl}, we get
\begin{align*}
    \mathcal{R} F_\omega 
    &= e^0 \wedge e^i \otimes (- 2 R_{k0i0} E^k  + R_{ik} E^k  - R_{00} E_i  - 2 R_{k0i\ell} B^{k\ell} + R_{0k} \tensor{B}{^k_i})\\
    &+ e^i \wedge e^j \otimes \frac12 \left( 2(R_{0ijk}-R_{kij0})E^k - R_{i0}E_j + R_{j0}E_i - 2R_{kij\ell} B^{k\ell} + R_{ik} \tensor{B}{^k_j} + R_{jk} \tensor{B}{_i^k} \right)
    \\[0.2cm]
    &= e^0 \wedge e^i \otimes \Big(R(g)_{ik} E^k - 3H \sff_{ik}E^k + \left(\tfrac{\nabla\sff}{\diff\tau}\right)_{ik} E^k - \tr\left(\tfrac{\nabla\sff}{\diff\tau}\right) E_i + 2|\sff|^2 E_i \\
    &\hspace{2cm} + 2 (D\sff)_{\ell ki}B^{k\ell} + 3(\diff H)_k \tensor{B}{^k_i} + (D^\ast \sff)_k \tensor{B}{^k_i} \Big)\\
    &+ e^i \wedge e^j \otimes \frac12 \Big(R(g)_{ik}\tensor{B}{^k_j} + R(g)_{jk}\tensor{B}{_i^k} - 2R(g)_{kij\ell} B^{k\ell} - \left(\tfrac{\nabla\sff}{\diff\tau}\right)_{ik}\tensor{B}{^k_j} -  \left(\tfrac{\nabla\sff}{\diff\tau}\right)_{jk}\tensor{B}{_i^k} \\
    &\hspace{2cm}- 3H (\sff_{ik} \tensor{B}{^k_j} + \sff_{jk} \tensor{B}{_i^k}) -2 \sff_{kj} \sff_{i\ell} B^{k\ell} + 2\tensor{\sff}{_i^\ell} \sff_{\ell k}\tensor{B}{^k_j} + 2\tensor{\sff}{_j^\ell} \sff_{\ell k} \tensor{B}{_i^k}  \\[0.1cm]
    &\hspace{2cm}+ 2\left( (D\sff)_{ijk} - (D\sff)_{jik} \right) E^k - 3(\diff H)_i E_j + 3(\diff H)_j E_i - (D^\ast \sff)_i E_j + (D^\ast \sff)_j E_i\Big).
\end{align*}

Inserting these calculations into (\ref{eq-yang-mills-wave}) proves Proposition \ref{prop-E-B-wave-eqs}.

\end{appendix}

\bibliographystyle{plain}
\bibliography{mybib}

\end{document}